\RequirePackage{fix-cm}
\documentclass[smallextended,envcountsame,numbook]{svjour3-vfinal}      
\smartqed  
\usepackage{tikz-cd}
\usepackage{cancel}
\usepackage{wrapfig}

\usepackage{tikz}
\usetikzlibrary{arrows.meta, positioning, calc}
\usetikzlibrary{shapes.misc}

\tikzset{
  nodebox/.style={draw, rounded corners=2pt, inner sep=3pt, minimum size=6mm}
}

\usepackage{mathtools}
\usepackage{booktabs}
\usepackage{array}

\usepackage{appendix}

\usepackage{graphicx}
\usepackage{nicematrix}
\usepackage{wrapfig}
\usepackage{svg}
\usepackage{dirtytalk}
\usepackage{vwcol}  

\usepackage{multicol}

 \usepackage{hyperref}
 \hypersetup{
     colorlinks=true,
     linkcolor=blue,
     filecolor=magenta,      
     urlcolor=cyan,
     citecolor=red,
}
\usepackage{amsfonts, amssymb, amsbsy, amsmath}
 \usepackage{blkarray}
 \usepackage{graphicx}
 \usepackage{float}
 \usepackage{epstopdf}
 \usepackage{algorithmic}
 \ifpdf
   \DeclareGraphicsExtensions{.eps,.pdf,.png,.jpg}
 \else
   \DeclareGraphicsExtensions{.eps}
 \fi
 \usepackage{mathrsfs}
 \usepackage{amscd,empheq} 
 \usepackage{yfonts,bm}

 \usepackage[pdf]{graphviz}

 \usepackage[all,2cell]{xy}
 \usepackage{pb-diagram}
 \usepackage{pb-xy}

 \usepackage{tikz}
 \usepackage{tikz-cd}
  \tikzcdset{
  cells={font=\everymath\expandafter{\the\everymath\displaystyle}},
}
 \usetikzlibrary{patterns}
 \usetikzlibrary{positioning}
 \usetikzlibrary{shapes,decorations.markings,arrows.meta}
 \usetikzlibrary{matrix,arrows}
 \tikzset{hfit/.style={rounded rectangle, inner xsep=0pt},
           vfit/.style={rounded corners}}

 \usepackage{morewrites}
 \usepackage{scalerel}
 \usepackage{mathabx}
 \usepackage{palatino}
 
\usepackage{enumitem}
\setlist[enumerate]{leftmargin=.5in}
\setlist[itemize]{leftmargin=.5in}

\usepackage{amsopn}

\newcommand{\rmOmega}{\mathrm{\Omega}}

\newcommand{\dff}{{\mathrm d}}

\newcommand{\sAtt}{{\mathsf{ Att}}}
\newcommand{\sANbhd}{{\mathsf{ ANbhd}}}

\newcommand{\sRep}{{\mathsf{ Rep}}}
\newcommand{\sMorse}{{\mathsf{ Morse}}}

\newcommand{\sMTile}{{\mathsf{MTile}}}
\newcommand{\sInvset}{{\mathsf{Invset}}}
\newcommand{\sCo}{{\mathsf{Co}}}

\newcommand{\scrT}{\mathscr{T}}

\newcommand{\pred}{\blacktriangleleft}

\newcommand{\Ext}{\mathrm{Ext}}

\renewcommand{\subset}{\subseteq}
\newcommand{\alphaOg}{\alpha_{\rm o}(\gamma_x^-)}

\definecolor{Blue}{RGB}{0,162,255}
\definecolor{Orange}{RGB}{243,144,25}
\definecolor{Red}{RGB}{236,93,87}

\newcommand{\bR}{{R}}

\newcommand{\K}{{\mathbb{K}}}
\newcommand{\N}{{\mathbb{N}}}
\newcommand{\R}{{\mathbb{R}}}

\newcommand{\Z}{{\mathbb{Z}}}

\newcommand{\cE}{{\mathcal S}}

\newcommand{\sA}{{\mathsf A}}

\newcommand{\sF}{{\mathsf F}}

\newcommand{\sJ}{{\mathsf J}}
\newcommand{\sO}{{\mathsf O}}
\newcommand{\sK}{{\mathsf K}}
\newcommand{\sL}{{\mathsf L}}
\newcommand{\sM}{{\mathsf M}}
\newcommand{\sN}{{\mathsf N}}
\newcommand{\sP}{{\mathsf P}}
\newcommand{\sQ}{{\mathsf Q}}

\newcommand{\sS}{{\mathsf S}}
\newcommand{\sT}{{\mathsf T}}
\newcommand{\sU}{{\mathsf U}}

\renewcommand{\imath}{{\mathsf i}}
\renewcommand{\jmath}{{\mathsf j}}

\newcommand{\sSet}{{\mathsf{Set}}}

\newcommand{\Hom}{\mathrm{Hom}}

\newcommand{\id}{\mathrm{id}}
\newcommand{\coker}{\mathrm{coker}}

\DeclareMathOperator{\Inv}{\mathrm{Inv}}
\DeclareMathOperator{\cl}{\mathrm{cl}}
\DeclareMathOperator{\Int}{\mathrm{int}}

\newcommand{\HC}{\mathrm{HC}}

\newcommand{\image}{\mathrm{im}}
\newcommand{\coimage}{\mathrm{coim}}

 \newcommand{\bt}{{\bm 0}}
\newcommand{\smbt}{{\scalebox{0.6}{$\bm{0}$}}}
\newcommand{\tp}{{\bm 1}}
\newcommand{\smtp}{{\scalebox{0.6}{$\bm{1}$}}}
\newcommand{\covv}{{\prec }}

\newcommand{\bPoset}{\text{{\bf Poset}}}
\newcommand{\bSet}{\text{{\bf Set}}}

\newcommand{\sRmodR}{\mathbf{Mod}_{R}}
\newcommand{\sRmodZ}{\mathbf{Mod}_{\Z}}
\newcommand{\sKVec}{\text{{\bf Vect}}_\K}

\newcommand{\Ch}{\mathbf{Ch}}
\newcommand{\Kh}{\mathbf{Kh}}
\newcommand{\Dh}{\mathbf{Dh}}

\newcommand{\bfE}{\text{\bf E}}

\newcommand{\bfD}{\text{\bf D}}

\newcommand{\bfS}{\text{\bf S}}

{\begin{snugshade}\begin{quote}}
{\hfill\end{quote}\end{snugshade}}

\definecolor{shadecolor}{rgb}{0.8,0.8,0.8}

\begin{document}

\setlength{\abovedisplayskip}{4pt}
\setlength{\belowdisplayskip}{3pt}

\title{Lattice-ordered Spectral Systems for  Morse Representations
}

\titlerunning{Lattice-ordered Spectral Systems for  Morse Representations}

\author{Samuel Castelo-Mourelle \and Robert Vandervorst 
}


\institute{S. Castelo-Mourelle \at
              VU University Amsterdam\\
              Amsterdam, The Netherlands\\
              \email{s.castelomourelle@vu.nl}          
           \and
R.C.A.M. Vandervorst \at
              VU University Amsterdam\\
              Amsterdam, The Netherlands\\
              \email{r.c.a.m.vander.vorst@vu.nl}
}

\date{Version date: \today}

\maketitle

\begin{abstract}
We show that 
  Morse representations can be 
algebraically described by spectral systems, a notion introduced by B.\ Matschke for poset-filtered chain complexes. This provides an algebraic alternative to connection matrices and yields algebraic invariants via $\bR$-modules  over an arbitrary ring $\bR$. This approach applies in particular to   Conley index theory for semi-flows. We discuss the theory of spectral systems and, as an application to dynamical systems, explain
 the relation to 
 R.\ Franzosa's theory of connection matrices.
 In this context we introduce the concept of a spectral net.
 In contrast to traditional spectral sequences, spectral nets over a lattice display a rich network of (higher) differentials. In the application to dynamical systems these differentials provide information about connecting orbits.

\end{abstract}

\tableofcontents

\section{Introduction}
\label{prelude}

In this paper three key ingredients come together. Chronologically, the first is the foundational work by R.\ Franzosa on connection matrices in the setting of dynamical systems, v.\ \cite{fran2,fran}. This provides an algebraic theory for Morse representations   with information about connecting orbits and bifurcations in flows. The second is the work by J.\ Robbin \& D.\ Salamon on connection matrix theory for discrete dynamical systems, v.\ \cite{robbin:salamon2}. In this work they utilize the novel observation that the attractors in a dynamical system display a lattice structure. The third is the introduction of spectral systems for filtered chain complexes by B.\ Matschke, which allow filterings over any poset, v.\ \cite{Matschke}. In this paper we introduce spectral systems in an abstract setting beyond filtered chain complexes. This abstract setting is applied to the theory of Morse representations and
extends the approach via spectral sequences by P.\ Bart{\l}omiejczyk, v.\ \cite{Bartl3,Bartl1}; cf.\ \cite{CorneaRez}.

The decomposition of a dynamical system $\phi$ into a hierarchical structure of invariant sets is a cornerstone of the global analysis of dynamical systems.\footnote{In this paper a dynamical system is a semi-flow on a compact topological space.} Central to this decomposition is the concept of \emph{Morse representation}, defined as a finite partially ordered set consisting of pairwise disjoint, compact, non-empty invariant sets --- \emph{Morse sets} --- outside of which $\phi$ has gradient-like behavior, v.\ Sect.\ \ref{tessMD} and \cite{conley:cbms,fran,lsa3,KVPriest}.
Formally, a Morse representation is equivalent to a finite sublattice of attractors $\sA$ of $\phi$.  To each Morse set $M$, i.e.\ the intersection of an attractor $A\in \sA$ and the dual repeller $A'^*$ of an attractor $A'\in \sA$, with $A'\subset A$, we assign the \emph{Conley index} $\HC(M)$, an algebraic-topological invariant defined as the relative homology $H(U,U')$, where $U,U'$, with $U'\subset U$, are trapping regions for $A$ and $A'$ respectively.
The interplay between the Conley indices is captured by a \emph{Cartan--Eilenberg system} $\bfE(\sM)$ (or \emph{$E$-system}), v.\ \cite{SpV4,HelleRognes}. This system provides an exhaustive algebraic framework consisting of $\bR$-modules $\HC(M)=H(U,U')$ for every pair $(A, A') \in \sA \times \sA$, with $A' \subset A$. Within this system, every nested pair of attractors $(A, A')$  yields an exact triangle of the form:
\begin{equation}
    \label{exacttrinr12}
\begin{tikzcd}
\cdots \arrow[r] & H(U') \arrow[r, "i"] & H(U) \arrow[r, "j"] & {H(U,U')} \arrow[r, "k"] & H(U') \arrow[r] & \cdots
\end{tikzcd}
\end{equation}
An $E$-system serves as a \emph{maximal system of invariants}, encoding the totality of algebraic (homological) data regarding Morse sets and the   connecting orbits between them through a network of commutative squares and exact triangles,
v.\ Sect. \ref{specMorserepr}. In \cite{SpV4} it is shown that $E$-systems for a Morse representation are equivalent to Franzosa's homology braids.

An equivalent formulation of this information is via a \emph{spectral system}\footnote{The term spectral system was coined in \cite{Matschke}, in the setting of filtered chain complexes.} $\bfS(\sM)$ (or \emph{$S$-system}), which represents the data as a collection of $\bR$-modules and differentials. If $(\sQ,\le)$ is a poset and $\bm{\sQ_n}$  is the poset of ordered $n$-tuples $\alpha_1\le\cdots\le\alpha_n$ in $\sQ$, ordered componentwise and regarded as a thin category, then a spectral system $\bfS$ is defined by a functor $\sS\colon \bm{\sQ_4}\to \sRmodR$, $(\alpha,\beta,\gamma,\delta)\mapsto S^{\gamma\delta}_{\alpha\beta}$, together with differentials $\dff\colon S^{\epsilon\zeta}_{\gamma\delta}\to S^{\gamma\delta}_{\alpha\beta}$, $(\alpha,\beta,\gamma,\delta,\epsilon,\zeta)\in\bm{\sQ_6}$, which satisfy a list of axioms, v.\ Defn.\ \ref{specsysdefn}.
In our setting $\sQ=\sO(\sM)$, the lattice of down-sets in $(\sM,\le)$, and
\begin{equation*}
S_{\alpha\beta}^{\gamma\delta}=\displaystyle  \frac{\ker \colon E^{\gamma}_{\beta} \xrightarrow{k} E^{\beta}_{\alpha}}{\image \colon E^{\delta}_{\gamma} \xrightarrow{k} E^{\gamma}_{\beta}},\quad \forall (\alpha,\beta,\gamma,\delta)\in \bm{\sO(\sM)_4},
\end{equation*}
where $E^\beta_\alpha = H(F_\beta U,F_\alpha U)$ are the $E$-terms for $\bfE(\sM)$ --- Conley indices --- and the sets $F_\alpha U$ define a lattice $\sU$ of trapping regions isomorphic to $\sA$, v.\ Sect.~\ref{tessMD}. The differentials $k\colon E^\gamma_\beta\to E^\beta_\alpha$ are the connecting homomorphisms  of $\bfE(\sM)$. Spectral systems indexed by a finite distributive lattice, such as $\sO(\sM)$, will be called \emph{lattice-ordered}.

\subsection{The spectral net approach}
\label{specnetapp}
In the tradition of Franzosa, chain-generated $E$-systems, v.\ Sect.~\ref{specMorserepr}, can be encoded by an \emph{$\sM$-graded differential module} $(C, \dff)$, also referred to as a \emph{connection matrix} (or \emph{spectral matrix} in the context of Cartan--Eilenberg systems). We emphasize that a connection matrix is not a reduction of the data but an alternative, highly structured encoding of the system. Knowing $(C, \dff)$ determines all $E$-terms, $S$-terms, and differentials. Connection matrices are not unique in general. The existence of connection matrices is proved by Franzosa \cite{fran}, but not over arbitrary rings. For the extension to arbitrary rings, v.\ \cite{Vander2}.
Because computing the full connection matrix requires global, chain-level information, which is often not available, we introduce the concept of a \emph{spectral net}. Analogous to spectral sequences in homological algebra, spectral nets provide a framework to deduce structural information about $E$-systems and $S$-systems from partial, localized data by exploiting the underlying spectral structure.

Spectral systems offer several key advantages over traditional connection matrices. First, they operate over an arbitrary ring $\bR$, providing broad algebraic generality. Second, their differentials are uniquely determined by the $E$-system, which makes them robust under perturbations of $\phi$, v.~\cite{fran3,DKV}, and non-trivial differentials carry direct information about the underlying connecting orbit structure, v.\ Theorems~A and~B below. In contrast, connection matrix entries are generally non-unique and lack this straightforward geometric interpretation. Finally, spectral systems make information more accessible: while computing a connection matrix requires complete, global knowledge of the entire system up front, spectral systems can be analyzed using localized or partial data via  spectral nets.
To exploit this accessibility,
spectral nets enable an iterative computational approach. Rather than relying on a complete system description, spectral nets incrementally recover information about the spectral system --- higher $S$-terms and higher differentials --- step-by-step from limited initial data: the $E$-terms $E^\beta_\alpha$ with $\beta\setminus\alpha$ a singleton, and the first-order differentials between them.

\subsection{An example}
\label{theN5}
We consider an $E$-system over the lattice of down-sets $\sO(\sP)$ of the poset $\sP=\sN_5$. Figure \ref{fig:downset-lattice} gives the Hasse diagrams of both $\sN_5$ and $\sO(\sN_5)$.
\begin{figure}[h] 
  \centering
\begin{tikzpicture}[scale=0.8, every node/.style={circle, draw, minimum size=2mm, inner sep=2pt, font=\small}]
  \node (1) at (0,0) {{\small 1}};         
  \node (2) at (-2,1.5) {{\small 2}};
  \node (3) at (2,2.25) {{\small 3}};
  \node (4) at (-2,3) {{\small 4}};
  \node (5) at (0,4.5) {{\small 5}};       

  \draw (1) -- (2);
  \draw (1) -- (3);
  \draw (2) -- (4);
  \draw (3) -- (5);
  \draw (4) -- (5);
\end{tikzpicture}
\quad\quad\quad
\begin{tikzpicture}[
  scale=0.8,
  every node/.style={draw, rounded corners, fill=white, inner sep=2pt, font=\small}
]

  \node (D0) at (0,0) {$\bt=\emptyset$};

  \node (D1) at (0,1.2) {$\alpha=\{1\}$};

  \node (D2) at (-2,2.4) {$\beta=\{1,2\}$};
  \node (D3) at (2,2.4) {$\gamma=\{1,3\}$};

  \node (D4) at (0,3.6) {$\delta=\{1,2,3\}$};
  \node (D5) at (-4,3.6) {$\epsilon=\{1,2,4\}$};

  \node (D6) at (-2.0,4.8) {$\zeta=\{1,2,3,4\}$};

  \node (D7) at (-2,6) {$\tp=\{1,2,3,4,5\}$};

  \draw (D0) -- (D1);
  \draw (D1) -- (D2);
  \draw (D1) -- (D3);
  \draw (D2) -- (D4);
  \draw (D2) -- (D5);
  \draw (D3) -- (D4);
  \draw (D4) -- (D6);
  \draw (D5) -- (D6);
  \draw (D6) -- (D7);

\end{tikzpicture}
 \caption{The poset $\sP = \sN_5= \{1,2,3,4,5\}$ with covering relations
 $1\!\prec\! 2\!\prec\! 4\!\prec\! 5$ and $1\!\prec \!3\!\prec\! 5$
 [left] and its lattice of down-sets $\sO(\sP)$ [right].}
 \vspace{-.5cm}
  \label{fig:downset-lattice}
\end{figure}
The lattice $\sO(\sN_5)$ has $3$ maximal chains with $3$  associated linear extensions of $\sN_5$:
\[
\begin{aligned}
    \sK_1\colon~~& \bt \subset\alpha\subset \beta \subset \epsilon\subset \zeta\subset \tp,\quad \implies \quad 1 \covv 2 \covv 4 \covv 3 \covv 5;\\
    \sK_2\colon~~& \bt \subset\alpha\subset \beta \subset \delta\subset \zeta\subset \tp,\quad \implies \quad 1\covv  2 \covv  3 \covv  4 \covv  5;\\
    \sK_3\colon ~~& \bt \subset\alpha\subset \gamma \subset \delta\subset \zeta\subset \tp,\quad \implies \quad 1\covv  3 \covv  2 \covv  4 \covv  5.
\end{aligned}
\]
From the structure of the Cartan--Eilenberg system we define page-$1$   of the \emph{spectral net} by considering $E$-terms corresponding to consecutive down-sets in the maximal chains.
We assume that the $E$-system  is \emph{excisive}, in particular $E^\beta_\alpha\cong E^{\beta'}_{\alpha'}$ for all $\beta\setminus\alpha=\beta'\setminus\alpha'$.
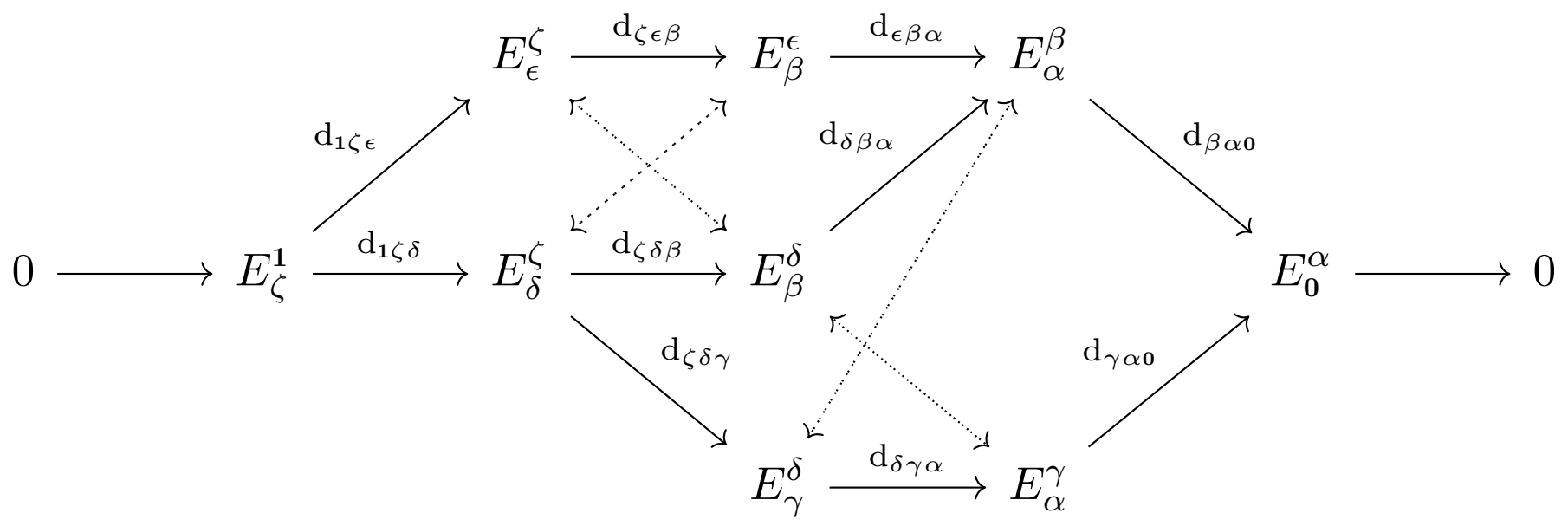
\[
\begin{tikzcd}[row sep=0.9cm, column sep =1.2cm]
            &                                                                                            & E^\zeta_\epsilon \arrow[r, "\dff_{\zeta\epsilon\beta}"] \arrow[rd, dotted]                                    & E^\epsilon_\beta \arrow[r, "\dff_{\epsilon\beta\alpha}"] \arrow[ld, dotted]                 & E^\beta_\alpha \arrow[rd, "\dff_{\beta\alpha\smbt}"] \arrow[ldd, dotted]  &                          &   \\
0 \arrow[r] & E^\smtp_\zeta \arrow[ru, "\dff_{\smtp\zeta\epsilon}"] \arrow[r, "\dff_{\smtp\zeta\delta}"] & E^\zeta_\delta \arrow[r, "\dff_{\zeta\delta\beta}"] \arrow[rd, "\dff_{\zeta\delta\gamma}"] \arrow[ru, dotted] & E^\delta_\beta \arrow[ru, "\dff_{\delta\beta\alpha}"] \arrow[lu, dotted] \arrow[rd, dotted] &                                                                           & E^\alpha_\smbt \arrow[r] & 0 \\
            &                                                                                            &                                                                                                               & E^\delta_\gamma \arrow[r, "\dff_{\delta\gamma\alpha}"] \arrow[ruu, dotted]                  & E^\gamma_\alpha \arrow[ru, "\dff_{\gamma\alpha\smbt}"] \arrow[lu, dotted] &                          &
\end{tikzcd}
\]
If we display each maximal chain we obtain the following diagram:
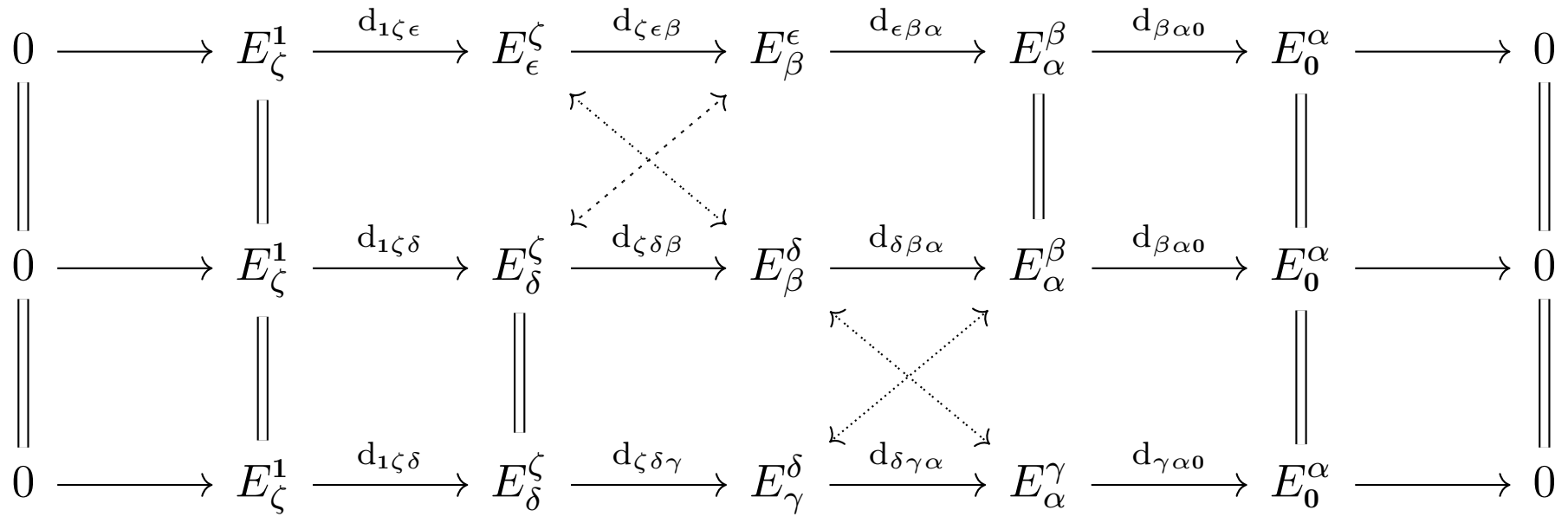
\[
\begin{tikzcd}[row sep=0.9cm, column sep =1.2cm]
0 \arrow[d, equal] \arrow[r] & E^\smtp_\zeta \arrow[d, equal] \arrow[r, "\dff_{\smtp\zeta\epsilon}"] & E^\zeta_\epsilon \arrow[r, "\dff_{\zeta\epsilon\beta}"] \arrow[rd, dotted]                & E^\epsilon_\beta \arrow[r, "\dff_{\epsilon\beta\alpha}"] \arrow[ld, dotted]                & E^\beta_\alpha \arrow[d, equal] \arrow[r, "\dff_{\beta\alpha\smbt}"]   & E^\alpha_\smbt \arrow[d, equal] \arrow[r] & 0 \arrow[d, equal] \\
0 \arrow[r] \arrow[d, equal] & E^\smtp_\zeta \arrow[r, "\dff_{\smtp\zeta\delta}"] \arrow[d, equal]   & E^\zeta_\delta \arrow[r, "\dff_{\zeta\delta\beta}"] \arrow[ru, dotted] \arrow[d, equal] & E^\delta_\beta \arrow[lu, dotted] \arrow[rd, dotted] \arrow[r, "\dff_{\delta\beta\alpha}"] & E^\beta_\alpha \arrow[r, "\dff_{\beta\alpha\smbt}"] \arrow[ld, dotted]   & E^\alpha_\smbt \arrow[r] \arrow[d, equal] & 0 \arrow[d, equal] \\
0 \arrow[r]                    & E^\smtp_\zeta \arrow[r, "\dff_{\smtp\zeta\delta}"]                      & E^\zeta_\delta \arrow[r, "\dff_{\zeta\delta\gamma}"]                                      & E^\delta_\gamma \arrow[r, "\dff_{\delta\gamma\alpha}"] \arrow[ru, dotted]                  & E^\gamma_\alpha \arrow[lu, dotted] \arrow[r, "\dff_{\gamma\alpha\smbt}"] & E^\alpha_\smbt \arrow[r]                    & 0
\end{tikzcd}
\]
where the dotted arrows indicate the isomorphisms between $E$-terms due to the excisive property. The latter diagram makes it easier to display the higher pages of the spectral net.
By the excisive property we choose $\bR$-modules $G_pC$: $G_1C\cong E^{\alpha}_{\smbt}$, $G_2C\cong E^\beta_\alpha\cong E^\delta_\gamma$, $G_3C \cong E^\zeta_\epsilon\cong E^\delta_\beta\cong E^\gamma_\alpha$, $G_4C \cong E^\zeta_\delta\cong E^\epsilon_\beta$, and $G_5C\cong E^\smtp_{\zeta}$,
whose identifications are indicated with dotted arrows in the diagram above.
If we identify the isomorphic $E$-terms with $G_1C,\ldots,G_5C$, then we obtain the following diagram:
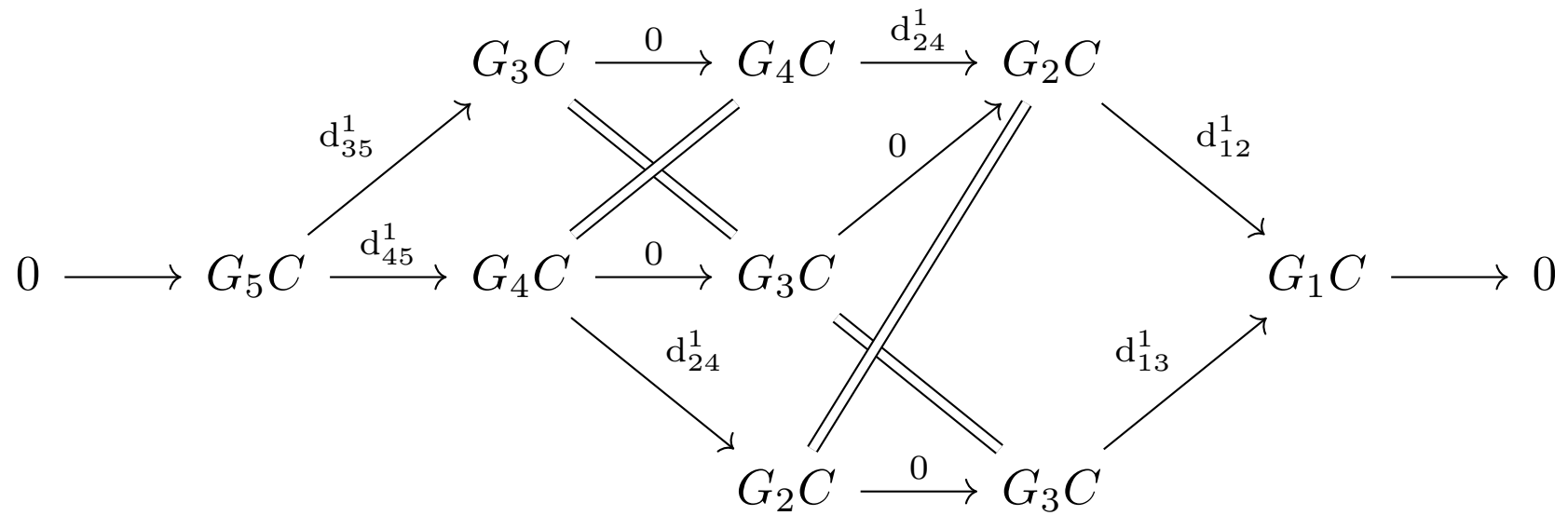
\begin{equation}
    \label{page1forC}
    \begin{tikzcd}[row sep=large]
            &                                                         & G_3C \arrow[rd, equal] \arrow[r, "0"]                           & G_4C \arrow[r, "\dff^1_{24}"] & G_2C \arrow[rd, "\dff^1_{12}"] &                &   \\
0 \arrow[r] & G_5C \arrow[ru, "\dff^1_{35}"] \arrow[r, "\dff^1_{45}"] & G_4C \arrow[ru,equal] \arrow[r, "0"] \arrow[rd, "\dff^1_{24}"] & G_3C \arrow[rd, equal] \arrow[ru, "0"]          &                                & G_1C \arrow[r] & 0 \\
            &                                                         &                                               & G_2C \arrow[ruu, equal] \arrow[r, "0"]           & G_3C \arrow[ru, "\dff^1_{13}"] &                &
\end{tikzcd}
\end{equation}
Repeated application of homology yields the \say{higher} $S$-terms and associated differentials.
The following diagram displays  page-2:
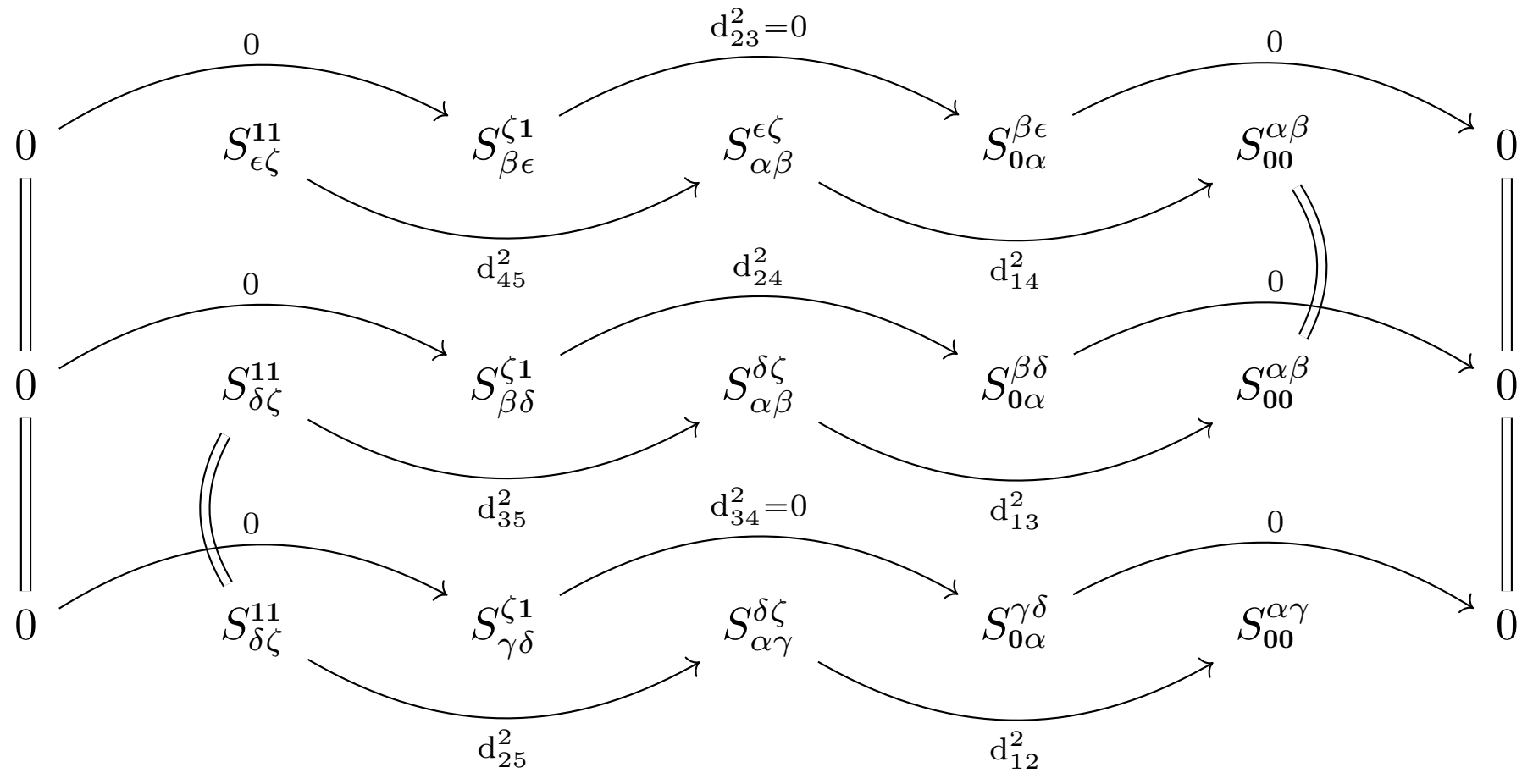
\[
\begin{tikzcd}[row sep=1.1cm, column sep=1.1cm]
0 \arrow[d, equal] \arrow[rr, "0", bend left] & S^{\smtp\smtp}_{\epsilon\zeta} \arrow[rr, "\dff^2_{45}"', bend right]                          & S^{\zeta\smtp}_{\beta\epsilon} \arrow[rr, "\dff^2_{23}=0", bend left] & S^{\epsilon\zeta}_{\alpha\beta} \arrow[rr, "\dff^2_{14}"', bend right] & S^{\beta\epsilon}_{\smbt\alpha} \arrow[rr, "0", bend left] & S^{\alpha\beta}_{\smbt\smbt} \arrow[d, equal, bend left] & 0 \arrow[d, equal] \\
0 \arrow[d, equal] \arrow[rr, "0", bend left] & S^{\smtp\smtp}_{\delta\zeta} \arrow[d, equal,  bend right] \arrow[rr, "\dff^2_{35}"', bend right] & S^{\zeta\smtp}_{\beta\delta} \arrow[rr, "\dff^2_{24}", bend left]   & S^{\delta\zeta}_{\alpha\beta} \arrow[rr, "\dff^2_{13}"', bend right]   & S^{\beta\delta}_{\smbt\alpha} \arrow[rr, "0", bend left]   & S^{\alpha\beta}_{\smbt\smbt}                            & 0 \arrow[d, equal] \\
0 \arrow[rr, "0", bend left]                    & S^{\smtp\smtp}_{\delta\zeta} \arrow[rr, "\dff^2_{25}"', bend right]                            & S^{\zeta\smtp}_{\gamma\delta} \arrow[rr, "\dff^2_{34}=0", bend left]  & S^{\delta\zeta}_{\alpha\gamma} \arrow[rr, "\dff^2_{12}"', bend right]  & S^{\gamma\delta}_{\smbt\alpha} \arrow[rr, "0", bend left]  & S^{\alpha\gamma}_{\smbt\smbt}                           & 0
\end{tikzcd}
\]
Spectral nets are based on partial information of the spectral system. At page-1 we start with the modules $G_pC$ and differentials defined between consecutive (in the linear extensions of $\sP$) $G$-terms given by the diagram in \eqref{page1forC}. Each path of solid arrows from $G_5C$ to $G_1C$ in \eqref{page1forC} represents a spectral sequence associated to a maximal chain in $\sO(\sP)$. The \say{bridges} (equal signs between $G$-terms) allow for additional paths and therefore new relations on the differentials, v.\ Thm.\ \ref{firstantocomm}.
In Figures \ref{perspec1} and \ref{perspec2} additional relations are displayed.
In Example \ref{franexample} we give another example of the complexity of relations between differentials in a spectral net applied to a problem in \cite[Ex.\ 6.2]{fran}.

\subsection{Main results}
\label{mainresults}
In Section \ref{tessMD} we start with the dynamical context of this paper. A \emph{semi-flow} $\phi$ on a compact topological space $(X,\scrT)$  is given by
 a continuous function $\phi\colon\R^+\times X \to X$, with $\R^+=[0,\infty)$, that satisfies: (i) $\phi^0(x)=x$ for all $x\in X$, (ii) $\phi^{t+s}=\phi^t\circ\phi^s$ for all $s,t\in\R^+$ (semi-group property), and (iii) $\phi^t$ is proper for all $t\in\R^+$, where $\phi^t:=\phi(t,\cdot)$, v.\ \cite{KVPriest}. We summarize the duality between Morse representations and finite sublattices of attractors for $\phi$.
In Section \ref{specMorserepr} we define the local invariants for Morse sets, \emph{the Conley index}, and we introduce an algebraic structure that combines the \say{local} Conley indices into a \say{global} structure via  a Cartan--Eilenberg system, v.\ \cite{SpV4}. In \cite{fran2,fran} the global structure is given by the equivalent structure called a homology braid,\footnote{In Franzosa's algebraic theory of Morse representations the notion of homology braid is introduced, v.\ \cite{fran2,fran}. In \cite{SpV4} it is proved that homology braids and excisive $E$-systems are equivalent. The disadvantage of the former is that it only records the convex sets of Morse sets, whereas $E$-systems keep track of the defining pairs of attractors, which makes properties and results   more accessible. $E$-systems are also well-defined without the excisive property.} v.\ \cite{SpV4}.
In Section \ref{specsys} we introduce the notion of a spectral system motivated by Matschke's notion of spectral systems for poset-filtered chain complexes, v.\ \cite{Matschke}. We give a  new axiomatic definition of spectral system and we prove that spectral systems are equivalent to Cartan--Eilenberg systems, v.\ Thm. \ref{duality}.
In Section \ref{spectralnets} we give a formal definition of the fundamental spectral net for an excisive $S$-system. The first page consists of the $E$-terms $E^\beta_\alpha$, with $\beta\setminus\alpha=\{p\}$ for all $p\in \sP$, and the higher pages are inductively defined by successively applying homology. In Section \ref{theN5} we already explained the basics via an example. Other spectral nets such as subnets and nets of $E$-terms of convex sets can be defined to obtain different information about the spectral system and its differentials.
In Section \ref{properties} we study properties of the differentials in an $S$-system and in particular the higher differentials.
Section \ref{firstordinfo} starts with the notion of \emph{perspective pairs} of down-sets, which provide additional relations on the first-order differentials of a spectral net besides the usual relations coming from the spectral sequences, v.\ Thm.\ \ref{firstantocomm}.
In Section \ref{deriveddiff} we show that first-order differentials allow \emph{derived differentials} that carry information about extensions, v.\ Thm.\ \ref{exactfortr}, and in Section \ref{highdiffspec} various properties of higher differentials are discussed, v.\ Thm. \ref{propo0} and the subsequent results.
In Section \ref{connorbits} we apply the properties of differentials
to Morse representations, which yields direct information about the structure of connecting orbits.
For flow-orders, v.\ Defn.\ \ref{flowdefn} and \cite{fran}, we obtain the following result for connecting orbits, where $C(M',M)$ denotes the set of connecting orbits from $M'$ to $M$, and $\bt=\varnothing$ and $\tp=\sM$ are the bottom and top of $\sO(\sM)$:
\vspace{.2em}

\noindent \textbf{Theorem A [v.\ Thm.\ \ref{conn2aa}]}
\emph{
Let $(\sM,\le)$ be a flow-order. Suppose the differential
    \[
    \dff\colon S^{\delta\smtp}_{\beta\gamma}\rightarrow S^{\beta\gamma}_{\smbt\alpha}
    \]
    is non-trivial, with $\beta\setminus\alpha=\{M\}$ and $\delta\setminus\gamma = \{M'\}$. Then, there exists a directed path $\{M_k\}_{k=0}^{n}\subset \delta\setminus\alpha$
    from $M'=M_n$ to $M=M_0$ (a path of covering relations, v.\ Cor.~\ref{propo1}), such that for all $1\le k\le n$ either
    $C(M_k,M_{k-1})\neq \emptyset$, or $C(M_{k-1},M_{k})\neq \emptyset$, depending on whether $M_k\succ M_{k-1}$, or $M_k\prec M_{k-1}$ respectively.
}
\vspace{.2em}

Other results are described in Section \ref{connorbits}, such as
Theorem \ref{conn2}, which extend results in \cite{fran} and  \cite{Bartl1}.
In contrast, the interpretation of connection matrices is less direct as their higher-order entries depend on extension data.
\vspace{.2em}

\noindent \textbf{Theorem B [v.\ Thm.\ \ref{connviahigherdiff}]}
\emph{
Let $(\sM,\le)$ be a Morse representation and
suppose the differential
    \[
    \dff\colon S^{\epsilon\zeta}_{\gamma\delta}\rightarrow S^{\gamma\delta}_{\alpha\beta}
    \]
    is non-trivial, with $\gamma\setminus\beta=\{M\}$ and $\epsilon\setminus\delta = \{M'\}$. Then,
    \begin{enumerate}[label=\textup{(\roman*)},leftmargin=2.2em,itemsep=0.1em,topsep=0pt]
        \item [{\rm (i)}] there exists a Morse set $M''\in \delta\setminus\beta$, $M''<M'$,
   such that $C(M',M'')\neq \emptyset$;
    \item[{\rm (ii)}] there exists a Morse set $M'''\in \epsilon\setminus\gamma$, $M<M'''$, such that
    $C(M''',M)\neq \emptyset$.
    \end{enumerate}
    }
\vspace{.2em}

Theorem~B is a key result: in contrast to most previous results, it extracts information about connecting orbits from (higher) differentials alone, and it applies to arbitrary Morse representations, not only to flow-orders.
In Section \ref{examplesofconn} we illustrate the results on spectral nets and the interpretation of higher differentials by some basic examples, including the examples in \cite{fran}.

\section{Morse representations}
\label{tessMD}
Let $\phi$ be a semi-flow on a compact topological space $X$, v.\ Sect.~\ref{mainresults}.
A \emph{trapping region} $U \subset X$ is a forward invariant attracting neighborhood, i.e.\ $\phi^t(U)\subset U$ for all $t\ge 0$ and $\phi^\tau(\cl U)\subset \Int U$ for some $\tau>0$.  Every  trapping region $U$ defines an \emph{attractor} $A= \omega(U)$, the $\omega$-limit set of $U$, which is non-empty  and invariant for $\phi^t$ if $U$ is non-empty, v.\ \cite{KVPriest}. The trapping regions  $\sANbhd_+(\phi)$ (with $\cup$ and $\cap$) and the attractors  $\sAtt(\phi)$ (with $A\vee A'=A\cup A'$ and $A\wedge A'=\omega(A\cap A')$) form bounded distributive lattices, v.\ \cite{lsa3}. All sublattices of $\sANbhd_+(\phi)$ and $\sAtt(\phi)$ are assumed to contain the bottom and the top element, i.e.\ $\varnothing$ and $X$, resp.\ $\varnothing$ and $\omega(X)$. The relationship between attractors and trapping regions is captured by the surjective lattice homomorphism\footnote{If there is no ambiguity about the dynamical system we write
$\omega_\phi(U)=\omega(U)$.}
\begin{equation}
    \label{omegamap}
\omega_\phi\colon \sANbhd_+(\phi) \twoheadrightarrow \sAtt(\phi),\quad U \mapsto \omega_\phi(U)=A.
\end{equation}
Every attractor has a well-defined \emph{dual repeller} $A^* = \alpha(U^c)$, the $\alpha$-limit set of $U^c=X\setminus U$, where $U$ is a trapping region with $\omega(U)=A$; the set $U^c$ is called a \emph{repelling region}, v.\ \cite{lsa,lsa3,KVPriest}. The lattice of repellers is denoted by $\sRep(\phi)$. A pair of attractors $A,A'\in \sAtt(\phi)$  defines a \emph{Morse set} $M=A-A':=A\cap A'^*$,\footnote{On every bounded, distributive lattice there exists a unique (up to isomorphism) monotone, distributive lattice form which satisfies  the axioms and properties of set-difference in a Boolean algebra. In the context of dynamics the lattice form is referred to as the Conley form, v.\ \cite{lsa3}.} and a pair   of trapping regions $U,U'\in \sANbhd_+(\phi)$ defines a \emph{Morse tile} $T=U\setminus U':=U\cap U'^c$. This defines a homomorphism $\Inv\colon \sMTile(\phi)\twoheadrightarrow \sMorse(\phi)$, via $\Inv(T)= \omega(U)-\omega(U')$, with $T=U\setminus U'$, which is the maximal invariant set in $T$ (so that it does not depend on the choice of $U$ and $U'$), and where $\sMTile(\phi)$ and $\sMorse(\phi)$ denote the semilattices of Morse tiles and Morse sets respectively, v.\ \cite{lsa3}.

For every finite sublattice $\sU\subset \sANbhd_+(\phi)$ of trapping regions the homomorphism in \eqref{omegamap} yields a finite sublattice $\sA\subset \sAtt(\phi)$ of attractors via  $\sU\ni U\mapsto \omega(U)=A\in \sA$, so that the diagram
\[
\begin{tikzcd}[row sep=large, column sep=large]
\sANbhd_+(\phi) \arrow[r, "\omega", two heads]                   & \sAtt(\phi)                      \\
\sU \arrow[u, "\subseteq"',"i", tail] \arrow[r, "\omega", two heads] & \sA \arrow[u, "\subseteq","i"', tail]
\end{tikzcd}
\]
commutes.
By Birkhoff's representation theorem the surjective homomorphism $\omega\colon \sU\to \sA$ dualizes to an order-embedding $\sJ(\omega)\colon\sJ(\sA) \hookrightarrow \sJ(\sU)$, where $\sJ$ is the contravariant functor that assigns to a finite, distributive lattice its poset of join-irreducible elements.\footnote{Let $\sL$ be a finite, distributive lattice. An element $a\in \sL$ is \emph{join-irreducible} if it has a unique immediate predecessor $a^\pred$, which is equivalent to $a\neq 0$ and: $a=b\vee c$ implies $b=a$ or $c=a$.}
Define $\sJ(\sA)\ni A \mapsto M= A-A^\pred$, a Morse set, and
$\sJ(\sU)\ni U \mapsto T=U\setminus U^\pred$, a Morse tile. We obtain the order-embedding\footnote{The map $\iota$ can be defined explicitly via Birkhoff duality, v.\ \cite{davey:priestley}: $\iota(A-A^\pred)=U\setminus U^\pred$, where $U=\sJ(\omega)(A)$ is the smallest element of $\sU$ with $\omega(U)=A$. However, the left inverse $\Inv$ determines the map implicitly: since the tiles in $\sT(\sU)$ are pairwise disjoint and the Morse sets in $\sM(\sA)$ are non-empty, $\iota(M)$ is the unique $T\in\sT(\sU)$ with $\Inv(T)=M$.}
\begin{equation}
    \label{defnTMD}
    \iota\colon\sM(\sA)\hookrightarrow \sT(\sU),\quad \Inv(\iota(M)) =M,
\end{equation}
where $\sM(\sA)=\bigl\{M=A-A^\pred\mid A\in \sJ(\sA)\bigr\}$, with $M\le M'$ if
and only if $A\subset A'$, $A,A'\in \sJ(\sA)$ (where $M'=A'-A'^\pred$), $\sT(\sU)=\bigl\{T=U\setminus U^\pred\mid U\in \sJ(\sU)\bigr\}$, with $T\le T'$ if
and only if $U\subset U'$, $U,U'\in \sJ(\sU)$ (where $T'=U'\setminus U'^\pred$), and $\Inv(\iota(M))$ is the maximal invariant set in $T=\iota(M)$. The poset $\sM(\sA)$ is called a \emph{Morse representation}, the poset $\sT(\sU)$ is called a \emph{Morse tessellation} and the embedding in \eqref{defnTMD} is called a \emph{(tessellated) Morse decomposition}, v.\ \cite{lsa3} for details.

Notationally, the poset $\sT(\sU)$ is denoted by $\sP$ and the isomorphism $\sO(\sP)\cong \sU$, where $\sO(\sP)$ denotes the lattice of down-sets in $\sP$, v.\ \cite{davey:priestley}, is expressed as the lattice filtering
\[
\sO(\sP) \xlongrightarrow{\cong} \sU,\quad
\alpha\mapsto F_\alpha U\in \sU,
\]
v.\ \cite{SpV4}.
The elements of $\sT(\sU)$ are denoted by $G_p U$, $p\in\sP$, i.e.\ $G_pU=F_{\downarrow p}U\setminus F_{\downarrow p\setminus\{p\}}U$, where $\downarrow p$ is the down-set generated by $p$. For any convex set $\beta\setminus\alpha\in \sCo(\sP)$, where $\sCo(\sP)$ denotes the meet-semilattice of convex sets in $\sP$, $G_{\beta\setminus\alpha} U= F_\beta U\setminus F_\alpha U$.
This does not depend on the choice of $\alpha$ and $\beta$, since
\[
G_{\beta\setminus\alpha} U = \bigcup_{p\in \beta\setminus\alpha} G_pU,
\]
which defines the semilattice embedding $\sCo(\sP)\rightarrowtail\sMTile(\phi)$.
The following theorem provides an equivalent characterization of a Morse representation.

\begin{theorem}[\cite{lsa3}, Thm.\ 8.5]
\label{thm:mainMDthm}
Let $(\sM,\le)$ be a finite poset of non-empty, pairwise disjoint, compact, invariant sets\footnote{An invariant set  of $\phi$ is a subset $S\subset X$ that satisfies $\phi^t(S)=S$ for all $t\ge 0$. We denote the lattice of invariant sets by $\sInvset(\phi)$, v.\ \cite{lsa}.} in $X$. Then, $(\sM,\le)$
is a Morse representation subordinate to a finite sublattice $\sA(\sM)\subset \sAtt(\phi)$ if and only if
for every $x\in X$ there exists an $M\in\sM$ such that
$\omega(x)\subset M$, and for each
complete orbit $\gamma_x$\footnote{A complete orbit through $x$ is a map $\gamma_x\colon\R\to X$ with $\gamma_x(0)=x$ and $\phi^t(\gamma_x(s))=\gamma_x(s+t)$ for all $s\in\R$ and $t\in\R^+$. Its restriction $\gamma_x^-$ to $(-\infty,0]$ is a backward orbit through $x$. For the $\omega$-limit set $\omega(x)$ of $x$ and the $\alpha$-limit set $\alphaOg$ of $\gamma_x^-$, v.\ \cite{KVPriest}.} with
$x\not \in \bigcup_{M\in \sM} M$
  there exist $M<M'$ such that $\omega(x)\subset M$ and  $\alphaOg \subset M'.$
The associated
lattice
$\sA(\sM)$ is the image of the injective lattice homomorphism $\nu\colon\sO(\sM) \to \sInvset(\phi)$ given by
\begin{equation}\label{eq:numap}
\alpha \mapsto \nu(\alpha) =
\bigcup_{M\in \alpha} W^u(M) \in \sAtt(\phi),
\end{equation}
where $W^u(M) =\bigl\{x\in X\mid \text{there exists a backward orbit $\gamma_x^-$ with }\alphaOg\subset M\bigr\}$ is the unstable set of $M$. Moreover, $\sM = \sM\bigl(\sA(\sM)\bigr)$.
\end{theorem}

Let $\sF\colon \bPoset \to \bSet$ be the forgetful functor keeping the underlying set of a poset, and let $\sF(\sM)\subset \sSet(X)$, where $\sSet(X)$ denotes the set of all subsets of $X$, be the collection of Morse sets in a Morse representation $(\sM,\le)$.
\begin{definition}
    \label{flowdefn}
    A finite sublattice $\sA\subset \sAtt(\phi)$ is \emph{flow-maximal} if for every finite sublattice $\sA'\subset \sAtt(\phi)$ such that $\sF\bigl(\sM(\sA')\bigr)
    = \sF\bigl(\sM(\sA)\bigr)$, it holds that
    \[
    \sA'\subset \sA.
    \]
    The associated Morse representation $\bigl(\sM(\sA),\le\bigr)$ is called a \emph{flow-order} for $\phi$.\footnote{In \cite{fran} the order is referred to as the \emph{flow-defined order}, or \emph{flow-order} for $\phi$. In \cite{KVPriest} the latter is called an \emph{order model}.}
\end{definition}

A flow-order for $\phi$ is obtained from a Morse representation $(\sM,\le)$ by removing order relations (in finitely many steps) until the order is minimal with respect to the characterization in Theorem \ref{thm:mainMDthm}; since $\nu(\alpha)$ in \eqref{eq:numap} does not depend on the order, a smaller order, which has more down-sets, yields a larger sublattice $\sA(\sM)$. This implies that a flow-order $\bigl(\sM(\sA),\le\bigr)$ is characterized by the property: $M < M'$ if and only if there exists a chain $\{M_i\}_{i=0}^{n}\subset \sM(\sA)$, with $M_0=M$ and $M_n=M'$, such that
the sets of connecting orbits $C\bigl(M_i,M_{i-1}\bigr)$ are non-empty for all $1\le i\le n$; here $C(M',M)$ denotes the set of complete orbits $\gamma_x$ with $\alphaOg\subset M'$ and $\omega(x)\subset M$.

\section{Algebraic invariants}
\label{specMorserepr}
In order to organize the local data in a Morse representation, we start with a global structure reminiscent of homology braids as introduced in \cite{fran2,fran}, v.\ \cite{SpV4}.

\subsection{Cartan--Eilenberg systems}
\label{CAsystems}
Consider an arbitrary poset $\sQ$ and let $\bm{\sQ_1}$ be the associated posetal category. By $\bm{\sQ_2}$ we denote the arrow category of $\bm{\sQ_1}$, which consists of ordered pairs $(\alpha,\beta)$, $\alpha\le\beta$, ordered componentwise, i.e.\ $(\alpha,\beta)\le(\alpha',\beta')$ if and only if $\alpha\le\alpha'$ and $\beta\le\beta'$.
The (functor) categories $\bm{\sQ_n}$ consist of concatenations of $n-1$ arrows and are given by $n$-tuples
$\alpha_1\le\cdots\le\alpha_{n}$, again ordered componentwise, which may be regarded as a path of length $n-1$ in the directed graph induced by the poset $\sQ$.
Let $\sRmodR$ be the category of $\bR$-modules over a ring $\bR$.

\begin{definition}[cf.\ \cite{SpV4}]
    \label{CEsys}
A \emph{Cartan--Eilenberg system}, or \emph{$E$-system}, $\bfE$ over $\sQ$ consists of:
\begin{enumerate}[label=\textup{(\roman*)},leftmargin=2.2em,itemsep=0.1em,topsep=0pt]
    \item [(i)] $E$-terms $E_{\alpha}^{\beta}\in \sRmodR$, for all  pairs $(\alpha,\beta)\in \bm{\sQ_2}$;
    \item [(ii)] homomorphisms    $\ell\colon E^\beta_\alpha \to E^{\beta'}_{\alpha'}$ for all pairs $(\alpha,\beta)\le(\alpha',\beta')$, such that $\ell'\circ\ell=\ell''$, where $\ell'\colon E^{\beta'}_{\alpha'} \to E^{\beta''}_{\alpha''}$ and
$\ell''\colon E^\beta_\alpha \to E^{\beta''}_{\alpha''}$, for all triples $(\alpha,\beta)\le (\alpha',\beta')\le (\alpha'',\beta'')$,
 and such that  $(\alpha,\beta) = (\alpha',\beta')$ yields the identity map $\ell=\id$. Conditions (i) and (ii) define a functor $\bm{\sQ_2}\to\sRmodR$;
 \item[(iii)] connecting homomorphisms $k\colon E^\gamma_\beta\to E^\beta_\alpha$ and
$k'\colon E^{\gamma'}_{\beta'}\to E^{\beta'}_{\alpha'}$ such that the diagrams\footnote{In the special cases $(\alpha,\beta)\le (\alpha,\gamma)$ and $(\alpha,\gamma)\le (\beta,\gamma)$ the morphisms $\ell$ are denoted by $i$ and $j$ respectively.}
\begin{equation}
\label{unpackeddiag12}
\begin{tikzcd}[column sep=small]
E^\beta_\alpha \arrow[rr, "i"] &                   & E^\gamma_\alpha \arrow[ld, "j"] \\
                  & E^\gamma_\beta \arrow[lu, "k"] &
\end{tikzcd}
\qquad\qquad
\begin{tikzcd}
E^\gamma_\beta \arrow[r, "k"] \arrow[d, "\ell"'] & E^\beta_\alpha \arrow[d, "\ell"] \\
E^{\gamma'}_{\beta'} \arrow[r, "k'"]                 & E^{\beta'}_{\alpha'}
\end{tikzcd}
\end{equation}
are exact triangles and commutative squares respectively for all pairs of ordered triples $(\alpha,\beta,\gamma)\le (\alpha',\beta',\gamma')$.
\end{enumerate}
\end{definition}

By (ii) $\ell\colon E^\beta_\alpha \xrightarrow[\id]{\ell} E^\beta_\alpha$ is the identity homomorphism and the exactness of the left diagram in \eqref{unpackeddiag12} shows that $E^\alpha_\alpha =0$ for all $(\alpha,\alpha) \in \bm{\sQ_2}$, cf.\ \cite{HelleRognes}.
In many situations $\sQ$ is a modular\footnote{Distributive lattices, which are considered in this paper, are modular.} lattice, in which case one can define a notion of excision.
A Cartan--Eilenberg system is \emph{excisive} if the homomorphisms $\ell\colon$ $E^\alpha_{\alpha\wedge \beta} \to E^{\alpha\vee\beta}_\beta$ are isomorphisms for all $\alpha,\beta\in \sQ$.  In this paper we employ Cartan--Eilenberg systems over a finite, distributive lattice $\sQ=\sL$, represented as $\sL=\sO(\sP)$ for some finite poset $\sP$ via Birkhoff's representation theorem, v.\ \cite{davey:priestley}.
Since the elements in $\sO(\sP)$ are down-sets, the excisive property is equivalent to: $\ell\colon E^\beta_\alpha\to E^{\beta'}_{\alpha'}$ is an isomorphism for all $(\alpha,\beta)\le(\alpha',\beta')$ with $\beta\setminus\alpha=\beta'\setminus\alpha'$. In particular,
$E_\alpha^\beta\cong E_{\alpha'}^{\beta'}$ for all $\beta\setminus\alpha=\beta'\setminus\alpha'$, v.\ \cite{SpV4}.

A morphism between Cartan--Eilenberg systems $\bfE$ and $\bfE'$ over $\sQ$ is given by homomorphisms  $h^{\beta}_{\alpha}\colon E^\beta_\alpha \to E'^\beta_\alpha$ such that the diagram
\begin{equation}
\label{dia:exacttriangles}
    \begin{tikzcd}[row sep=large, column sep=large]
{} \arrow[r] & E_\alpha^\beta \arrow[r, "i"] \arrow[d, "h^\beta_\alpha"] & E^\gamma_\alpha \arrow[r, "j"] \arrow[d, "h^{\gamma}_{\alpha}"] & E^\gamma_\beta \arrow[r, "k"] \arrow[d, "h^{\gamma}_{\beta}"] & E_\alpha^\beta \arrow[r] \arrow[d, "h^{\beta}_{\alpha}"] & {} \\
{} \arrow[r] & E'^\beta_\alpha \arrow[r, "i'"]                           & E'^\gamma_\alpha \arrow[r, "j'"]                                & E'^\gamma_\beta \arrow[r, "k'"]                               & E'^\beta_\alpha \arrow[r]                                & {}
\end{tikzcd}
\end{equation}
commutes for any $\alpha\le\beta\le\gamma$. We write $h\colon \bfE\to \bfE'$, and $h$ is an \emph{isomorphism} if all $h^\beta_\alpha$ are isomorphisms.

\begin{remark}
\label{changeorder}
    In some settings it is necessary to consider Cartan--Eilenberg systems over different posets, which requires extending the notion of morphism between Cartan--Eilenberg systems. Let $\bfE$ and $\bfE'$ be Cartan--Eilenberg systems over $\sQ$ and $\sQ'$ respectively. Then, a morphism between $\bfE$ and $\bfE'$ is given by an order-preserving map $\psi\colon \sQ\to \sQ'$ and homomorphisms $h^{\beta}_{\alpha}\colon E^\beta_\alpha \to E'^{\psi(\beta)}_{\psi(\alpha)}$ such that the diagram in \eqref{dia:exacttriangles} is satisfied for all $\alpha\le\beta\le\gamma$, with the bottom row replaced by   $\bR$-modules $E'^{\psi(\beta)}_{\psi(\alpha)}$, $E'^{\psi(\gamma)}_{\psi(\alpha)}$ and $E'^{\psi(\gamma)}_{\psi(\beta)}$.
\end{remark}

In \cite[Lem.\ 4.8]{Matschke} it is shown that the information in a Cartan--Eilenberg system can be captured by only considering exact triangles and commutative squares of the form:\footnote{The bottom and top elements in $\sO(\sP)$ are  denoted by $\bt$ and $\tp$ respectively.}
\begin{equation}
\label{unpackeddiag344}
\begin{tikzcd}[column sep=small]
E^\alpha_\smbt \arrow[rr, "i"] &                   & E^\beta_\smbt \arrow[ld, "j"] \\
                  & E^\beta_\alpha \arrow[lu, "k"] &
\end{tikzcd}
\qquad\qquad
\begin{tikzcd}
E^\beta_\alpha \arrow[r, "k"] \arrow[d, "\ell"'] & E^\alpha_\smbt \arrow[d, "i"] \\
E^{\beta'}_{\alpha'} \arrow[r, "k'"]                 & E^{\alpha'}_\smbt
\end{tikzcd}
\end{equation}
   This reduced system is called an \emph{exact couple system} $\bfD$ and is reminiscent of Boardman's unrolled exact couple, cf.\ \cite{Boardman,HelleRognes}. An exact couple system fully reconstructs a Cartan--Eilenberg system.

\subsection{An invariant for  Morse representations}
\label{invforsect}
A Morse set $M$ is determined by a  pair of attractors. If $M=A-A'$, then
$A-A'=A-(A'\wedge A)$, which justifies considering only nested pairs $(A,A')$, $A'\subset A$.
A pair  $(U,U')$  is called an \emph{index pair}
for $M$ if (i) $U,U'$ are trapping regions with $U'\subset U$, and (ii) $M=A-A'$, with $A=\omega(U)$ and $A'=\omega(U')$. It holds that $M\subset \Int (U\setminus U')$. The term index pair can be used without   mention of the Morse set $M$.
A nested pair $U,U'$ of trapping regions with $U'\subsetneq U$ does not guarantee that $M$ is non-empty. Its relative homology yields a criterion to guarantee non-empty Morse sets $M$.
We define the \emph{Conley index} of $M$ by
\begin{equation}
    \label{defnCI}
     \HC(M) := H(U,U'),
\end{equation}
where $H$ denotes (co)homology with coefficients in $\bR$, e.g.\ singular homology or \v{C}ech cohomology. We deliberately ignore the grading and regard $H(U,U')$ as an $\bR$-module.
The index does not depend on the choice of index pair, i.e.\ if $(V,V')$ is an index pair with $M=\omega(V)-\omega(V')$, then $H(V,V')\cong H(U,U')$ via natural isomorphisms, which for nested index pairs are induced by inclusion, v.\ \cite{salamon}.
The Conley index contains information about Morse sets.
If $\HC(M)\neq 0$, then $M\neq \emptyset$, v.\ \cite{conley:cbms,salamon,Vander1}. The Conley index is invariant under continuation, v.\ \cite{conley:cbms,salamon}.

In this paper, however, we concentrate on Morse representations for a fixed dynamical system $\phi$.
For a Morse representation $(\sM,\le)$ for $\phi$ there exists a finite sublattice $\sU\subset\sANbhd_+(\phi)$ of trapping regions such that $\sM\cong \sT(\sU)$.\footnote{In the work of Franzosa such a lattice $\sU$ (more precisely, the filtering $\alpha\mapsto F_\alpha U$ below) is referred to as an \emph{index filtration}, v.\ \cite{fran2}.} In other words, since $\sM=\sM(\sA)$, with $\sA\subset \sAtt(\phi)$, $\sU$ is the image of a lift
$\sA\rightarrowtail \sANbhd_+(\phi)$, v.\ \cite{lsa}.
A lattice $\sU$ as indicated has the property that tiles in $\sT(\sU)$ correspond to (non-empty) Morse sets in $\sM$. As before, we represent $\sU$ by the lattice isomorphism
$\sO(\sM) \xlongrightarrow{\cong} \sU$ given by $\alpha\mapsto F_\alpha U\in \sU$.
The same applies to $\sO(\sM) \xlongrightarrow{\cong} \sA$ given by $\alpha\mapsto F_\alpha A\in \sA$.

Let $C(X)$ be an $\bR$-module of topological chains and let $\dff\colon C(X) \to C(X)$ be the differential making $\bigl(C(X),\dff\bigr)$ a differential module, for example the singular chain complex of $X$ with coefficients in $\bR$. The lattice filtering $\alpha\mapsto F_\alpha U$ yields a filtering of $C(X)$, i.e.\ $\alpha \mapsto F_\alpha C(X)$
   with $F_\alpha C(X) := C\bigl(F_\alpha U\bigr)$, where $C\bigl(F_\alpha U\bigr)$ are the topological chains in $C(X)$ restricted to $F_\alpha U$.
Observe that $F_{\alpha\cap \beta} C(X) = F_\alpha C(X) \cap F_\beta C(X)$.
   The same relation with  respect to union does \emph{not} hold in general. The differential satisfies $\dff F_\alpha C(X) \subset F_\alpha C(X)$ for all $\alpha\in \sO(\sM)$,  making $\alpha \mapsto F_\alpha C(X)$ a  meet semilattice filtered module.
   Define $C^\beta_\alpha = C(F_\beta U)/C(F_\alpha U)$ and $C^\alpha_\smbt = C(F_\alpha U)$, $\alpha\subset\beta$.
   We consider the following short exact sequences of differential modules.
Given a triple $\alpha\subset\beta\subset\gamma$ in $\sO(\sM)$, we have the short exact sequence:
\begin{equation}
    \label{shefilt1}
\begin{tikzcd}
0 \arrow[r] & C^\beta_\alpha \arrow[r, "i"] & C^\gamma_\alpha \arrow[r, "j"] & C^\gamma_\beta \arrow[r] & 0,
\end{tikzcd}
\end{equation}
and for pairs $\alpha\subset\beta$ and $\alpha'\subset\beta'$ in $\sO(\sM)$, with $\alpha\subset \alpha'$, $\beta\subset \beta'$, we define the compositions:
\begin{equation}
    \label{shefilt2}
\begin{tikzcd}
C^\beta_\alpha \arrow[r, "\subset"] & C^{\beta'}_{\alpha} \arrow[r, "\pi"] & C^{\beta'}_{\alpha'}.
\end{tikzcd}
\end{equation}
The homology  associated to the filtering of $\bigl(C(X),\dff\bigr)$ is denoted by $H(F_\alpha U) = H(C^\alpha_\smbt)$
and the   relative homology
is defined by
$H(F_\beta U,F_\alpha U) = H(C^\beta_\alpha)$.
The maps in \eqref{shefilt1} and \eqref{shefilt2}  yield the exact triangles:
\[
\begin{tikzcd}[row sep=large]
H(F_\beta U,F_\alpha U) \arrow[rr, "i"] &                   & H(F_\gamma U,F_\alpha U) \arrow[ld, "j"] \\
                  & H(F_\gamma U,F_\beta U) \arrow[lu, "k"] &
\end{tikzcd}
\]
and commutative squares:
\[
\begin{tikzcd}[row sep=large, column sep=large]
H(F_\gamma U,F_\beta U) \arrow[r, "k"] \arrow[d, "\ell"'] & H(F_\beta U,F_\alpha U) \arrow[d, "\ell"] \\
H(F_{\gamma'} U,F_{\beta'}U) \arrow[r, "k'"]                 & H(F_{\beta'}U,F_{\alpha'}U)
\end{tikzcd}
\]
where $k$ and $k'$ are the connecting homomorphisms, for all triples $(\alpha,\beta,\gamma)\le(\alpha',\beta',\gamma')$ in $\sO(\sM)$.
Via the above defined homology modules, the lattice $\sU$ defines a Cartan--Eilenberg system, which is denoted by $\bfE\bigl( \sT(\sU)\bigr)$.
The $E$-terms are given by
\[
E_\alpha^\beta := \HC(M)=H(F_\beta U,F_\alpha U),\quad \forall \alpha\subset \beta,~~\text{with}~~M=\omega(F_\beta U)-\omega(F_\alpha U),
\]
It is important to point out that $\bfE\bigl( \sT(\sU)\bigr)$ for a Morse representation is excisive even though $\bigl(C(X),\dff\bigr)$ is not necessarily lattice-filtered, i.e.\ $F_{\alpha\cup\beta}C(X)=F_\alpha C(X)+F_\beta C(X)$ need not hold.
Indeed, for $(\alpha,\beta)\le(\alpha',\beta')$ with $\beta\setminus\alpha = \beta'\setminus\alpha'$ we have that
\[
\omega(F_\beta U) -\omega(F_\alpha U) = \omega(F_{\beta'} U) -\omega(F_{\alpha'} U),
\]
v.\ \cite{lsa3}, and therefore, by the invariance of the Conley index, the homomorphism $\ell\colon$
$H(F_\beta U,F_\alpha U) \to H(F_{\beta'} U,F_{\alpha'} U)$, induced by the inclusion of index pairs, is an isomorphism, v.\ \cite{salamon},
which proves the excisive property of $\bfE\bigl( \sT(\sU)\bigr)$. Since the $E$-terms are defined as homology, we say that $\bfE\bigl( \sT(\sU)\bigr)$ is \emph{chain-generated}. More generally, an $E$-system over $\sO(\sP)$ is called \emph{chain-generated} if it is isomorphic to the $E$-system defined for instance as above by a differential module $(C,\dff)$ with a filtering $\alpha\mapsto F_\alpha C$ by submodules satisfying $\dff F_\alpha C\subset F_\alpha C$ and $F_{\alpha\cap\beta}C=F_\alpha C\cap F_\beta C$, v.\ \cite{Vander2}.
\begin{theorem}
    \label{CEsysinvforM}
    Let $\sT(\sU)\cong \sT(\sU')$ be two Morse tessellations for the Morse representation $\sM(\sA)$. Then,
    $\bfE\bigl( \sT(\sU)\bigr) \cong \bfE\bigl( \sT(\sU')\bigr)$ as Cartan--Eilenberg systems.
\end{theorem}
\begin{proof}
Since $\omega(F_\alpha U)=F_\alpha A=\omega(F_\alpha U')$ and $\omega$ is a lattice homomorphism, $F_\alpha V:=F_\alpha U\cap F_\alpha U'$ is a trapping region with $\omega(F_\alpha V)=F_\alpha A$ for all $\alpha\in\sO(\sM)$. The map $\alpha\mapsto F_\alpha V$ is monotone, which is all that is used in \eqref{shefilt1} and \eqref{shefilt2}. Replacing $F_\alpha U$ by $F_\alpha V$ therefore defines a Cartan--Eilenberg system $\bfE(V)$ with $E$-terms $H(F_\beta V,F_\alpha V)$. For $\alpha\subset\beta$ the inclusion $(F_\beta V,F_\alpha V)\subset(F_\beta U,F_\alpha U)$ induces chain maps between the short exact sequences in \eqref{shefilt1}, which commute with the maps in \eqref{shefilt2}. By naturality of the connecting homomorphisms the induced homomorphisms $h^\beta_\alpha\colon H(F_\beta V,F_\alpha V)\to H(F_\beta U,F_\alpha U)$ define a morphism $h\colon\bfE(V)\to\bfE\bigl(\sT(\sU)\bigr)$. Both pairs are index pairs for the Morse set $F_\beta A-F_\alpha A$. Hence each $h^\beta_\alpha$ is an isomorphism, v.\ \cite{salamon}, and $h$ is an isomorphism. In the same way $\bfE(V)\cong\bfE\bigl(\sT(\sU')\bigr)$.
\qed
\end{proof}
By Theorem~\ref{CEsysinvforM} the isomorphism class of the Cartan--Eilenberg system $\bfE\bigl( \sT(\sU)\bigr)$ is an invariant for $(\sM,\le)$, which we denote by $\bfE(\sM)$. In this case we write
\begin{equation}
    \label{CIEsystem12}
E^\beta_\alpha = \HC(M),\quad \forall \alpha\subset\beta,~~\text{with}~~M=\Inv\bigl(G_{\beta\setminus\alpha}U \bigr),
\end{equation}
i.e.\ up to isomorphism the $E$-terms of $\bfE(\sM)$ are the Conley indices of the Morse sets $\Inv\bigl(G_{\beta\setminus\alpha}U\bigr)=\omega(F_\beta U)-\omega(F_\alpha U)$.
\begin{remark}
    The term chain-generated is reminiscent of Franzosa's approach of assigning invariants to Morse representations, v.\ \cite{fran}. The approach in \cite{fran} is more restrictive since the complexes $C^\beta_\alpha$ are required to be excisive, i.e.\ to depend only on $\beta\setminus\alpha$ (equivalently, $F_{\alpha\cup\beta}C=F_\alpha C+F_\beta C$). This property is not needed in the theory of Morse representations presented in this paper.
\end{remark}

\section{Spectral systems}
\label{specsys}

Cartan--Eilenberg systems over a poset $\sQ$, v.\ Sect.\ \ref{CAsystems}, 
may be regarded as page-1 of a spectral system as introduced in \cite{Matschke}.
In this section we give an axiomatic definition of spectral systems that are not a priori generated by a Cartan--Eilenberg system.
\begin{definition}
    \label{specsysdefn}
A \emph{spectral system} $\bfS$, or \emph{$S$-system}, over a poset $\sQ$ consists of:
\begin{enumerate}[label=\textup{(\roman*)},leftmargin=2.2em,itemsep=0.1em,topsep=0pt]
    \item [(i)] $S$-terms $S_{\alpha\beta}^{\gamma\delta}\in \sRmodR$, for all  quadruples $(\alpha,\beta,\gamma,\delta)\in \bm{\sQ_4}$;
    \item [(ii)] homomorphisms $\ell\colon S^{\gamma\delta}_{\alpha\beta} \to S^{\gamma'\delta'}_{\alpha'\beta'}$ for all $(\alpha,\beta,\gamma,\delta)\le(\alpha',\beta',\gamma',\delta')$ in $\bm{\sQ_4}$, such that $\ell'\circ\ell=\ell''$, where $\ell'\colon S^{\gamma'\delta'}_{\alpha'\beta'} \to S^{\gamma''\delta''}_{\alpha''\beta''}$ and $\ell''\colon S^{\gamma\delta}_{\alpha\beta} \to S^{\gamma''\delta''}_{\alpha''\beta''}$,
for all triples of quadruples $(\alpha,\beta,\gamma,\delta)\le (\alpha',\beta',\gamma',\delta')\le (\alpha'',\beta'',\gamma'',\delta'')$ in $\bm{\sQ_4}$. Identical quadruples  $(\alpha,\beta,\gamma,\delta) = (\alpha',\beta',\gamma',\delta')$ yield the identity map $\ell=\id$. Conditions (i) and (ii) define a functor $\bm{\sQ_4}\to\sRmodR$;
\item[(iii)] differentials
$\dff\colon S_{\gamma\delta}^{\epsilon\zeta} \to S_{\alpha\beta}^{\gamma\delta}$ for all sextuples 
$(\alpha,\beta,\gamma,\delta,\epsilon,\zeta)\in \bm{\sQ_6}$,
and such that 
\begin{equation}
    \label{homdef1}
S_{\beta\delta}^{\epsilon\eta} \cong \frac{\ker\colon S_{\gamma\delta}^{\epsilon\zeta}\xrightarrow[]{\dff} S_{\alpha\beta}^{\gamma\delta} }{\image\colon S_{\epsilon\zeta}^{\eta\theta}\xrightarrow[]{\dff} S_{\gamma\delta}^{\epsilon\zeta}} = H\bigl(S_{\gamma\delta}^{\epsilon\zeta}\bigr),\quad\text{(homology)},
\end{equation}
for all octuples  $(\alpha,\beta,\gamma,\delta,\epsilon,\zeta,\eta,\theta)\in \bm{\sQ_8}$, where $S_{\epsilon\zeta}^{\eta\theta}\xrightarrow[]{\dff} S_{\gamma\delta}^{\epsilon\zeta} \xrightarrow[]{\dff} S_{\alpha\beta}^{\gamma\delta}$ satisfies $\dff^2=0$. The isomorphism in \eqref{homdef1} is induced by the homomorphisms $\ell$ in (ii), i.e.\ $\ell\colon S^{\epsilon\zeta}_{\beta\delta}\to S^{\epsilon\zeta}_{\gamma\delta}$ is injective with image $\ker\dff$, and $\ell'\colon S^{\epsilon\zeta}_{\gamma\delta}\to S^{\epsilon\eta}_{\gamma\delta}$ is surjective with kernel $\image\dff$;
\item [(iv)]  commutative diagrams
\begin{equation}
   \label{commdwithl} 
\begin{tikzcd}
S^{\epsilon\zeta}_{\gamma\delta} \arrow[r, "\dff"] \arrow[d, "\ell"'] & S^{\gamma\delta}_{\alpha\beta} \arrow[d, "\ell"] \\
S^{\epsilon'\zeta'}_{\gamma'\delta'} \arrow[r, "\dff"]                 & S^{\gamma'\delta'}_{\alpha'\beta'}               
\end{tikzcd}
\end{equation}
for all sextuples $(\alpha,\beta,\gamma,\delta,\epsilon,\zeta)\le (\alpha',\beta',\gamma',\delta',\epsilon',\zeta')$ in $\bm{\sQ_6}$;
\item [(v)] short exact sequences
\begin{equation}
    \label{shortexactSt}
\begin{tikzcd}
0 \arrow[r] & S_{\alpha\beta}^{\gamma\epsilon} \arrow[r, "i"] & S_{\alpha\beta}^{\delta\epsilon} \arrow[r, "j"] & S_{\alpha\gamma}^{\delta\epsilon} \arrow[r] & 0,
\end{tikzcd}
\end{equation}
for all quintuples $(\alpha,\beta,\gamma,\delta,\epsilon)\in \bm{\sQ_5}$.\footnote{In the special cases $S_{\alpha\beta}^{\gamma\epsilon}\to S_{\alpha\beta}^{\delta\epsilon}$ and $S_{\alpha\beta}^{\delta\epsilon} \to S_{\alpha\gamma}^{\delta\epsilon}$ the homomorphisms $\ell$ are denoted by $i$ and $j$ respectively.}
\end{enumerate} 
\end{definition}

\begin{remark}
\label{homiso}
The isomorphism in \eqref{homdef1} is induced by the homomorphisms $\ell$ as follows. Let $(\alpha,\beta,\gamma,\delta,\epsilon,\zeta,\eta,\theta)\in\bm{\sQ_8}$. The two conditions in (iii) concern the diagram
\[
\begin{tikzcd}[row sep=large, column sep=6em]
 & S^{\epsilon\zeta}_{\beta\delta} \arrow[d, "{\ell}"', tail] & \\
S^{\eta\theta}_{\epsilon\zeta} \arrow[r, "\dff", "{\image\dff\,=\,\ker\ell'}"'] & S^{\epsilon\zeta}_{\gamma\delta} \arrow[r, "\dff", "{\ker\dff\,=\,\image\ell}"'] \arrow[d, "{\ell'}"', two heads] & S^{\gamma\delta}_{\alpha\beta}\\
 & S^{\epsilon\eta}_{\gamma\delta} &
\end{tikzcd}
\]
in which $\ell$ and $\ell'$ are homomorphisms as in (ii): $\ell$ is injective with image $\ker\dff$, the kernel of the right-hand differential, and $\ell'$ is surjective with kernel $\image\dff$, the image of the left-hand differential. Since $\image\dff\subset\ker\dff$, we have $\image\dff=\ker\dff\cap\ker\ell'$, and therefore $\ell'$ induces an isomorphism
\[
H\bigl(S^{\epsilon\zeta}_{\gamma\delta}\bigr)=\frac{\ker\dff}{\image~\dff}=\frac{\ell\bigl(S^{\epsilon\zeta}_{\beta\delta}\bigr)}{\ell\bigl(S^{\epsilon\zeta}_{\beta\delta}\bigr)\cap\ker\ell'}\xrightarrow[\ \ell'\ ]{\ \cong\ }\ell'\ell\bigl(S^{\epsilon\zeta}_{\beta\delta}\bigr)\subset S^{\epsilon\eta}_{\gamma\delta}.
\]
To identify the image $\ell'\ell\bigl(S^{\epsilon\zeta}_{\beta\delta}\bigr)$, consider the square
\[
\begin{tikzcd}[row sep=large, column sep=large]
S^{\epsilon\zeta}_{\beta\delta} \arrow[r, "{\ell''}", two heads] \arrow[d, "{\ell}"', tail] & S^{\epsilon\eta}_{\beta\delta} \arrow[d, "{\ell'''}", tail]\\
S^{\epsilon\zeta}_{\gamma\delta} \arrow[r, "{\ell'}"', two heads] & S^{\epsilon\eta}_{\gamma\delta}
\end{tikzcd}
\]
of homomorphisms as in (ii). It commutes, since both compositions are the homomorphism $S^{\epsilon\zeta}_{\beta\delta}\to S^{\epsilon\eta}_{\gamma\delta}$ of (ii). The map $\ell''$ is surjective by the second condition in (iii) for the octuple $(\beta,\beta,\beta,\delta,\epsilon,\zeta,\eta,\theta)$, and $\ell'''$ is injective by the first condition in (iii) for the octuple $(\alpha,\beta,\gamma,\delta,\epsilon,\eta,\theta,\theta)$. Hence $\ell'\ell\bigl(S^{\epsilon\zeta}_{\beta\delta}\bigr)=\ell'''\ell''\bigl(S^{\epsilon\zeta}_{\beta\delta}\bigr)=\ell'''\bigl(S^{\epsilon\eta}_{\beta\delta}\bigr)$, and $\ell'''$ maps $S^{\epsilon\eta}_{\beta\delta}$ isomorphically onto this image.

Combining the two isomorphisms, the isomorphism in \eqref{homdef1} is
\[
S^{\epsilon\eta}_{\beta\delta}\xrightarrow{\ \cong\ }H\bigl(S^{\epsilon\zeta}_{\gamma\delta}\bigr),\qquad x\mapsto\bigl[\ell(y)\bigr],\quad\text{where $y\in S^{\epsilon\zeta}_{\beta\delta}$ with $\ell''(y)=x$},
\]
with inverse $[\ell(y)]\mapsto\ell''(y)$. The class $[\ell(y)]$ does not depend on the choice of $y$: if $\ell''(y)=\ell''(y')$, then $\ell'\ell(y-y')=\ell'''\ell''(y-y')=0$, so that $\ell(y-y')\in\ker\ell'=\image\dff$.

Conversely, the two conditions in (iii) are equivalent to the exactness of the sequences \eqref{specialsheseq2111} in Lemma~\ref{exactinSpS}. For the spectral system $\bfS(\bfE)$ of a Cartan--Eilenberg system $\bfE$, v.\ Sect.~4.2, both conditions hold: with $S^{\gamma\delta}_{\alpha\beta}\cong\image\bigl(E^\gamma_\alpha\to E^\delta_\beta\bigr)$, $\ell$ is an inclusion of submodules of $E^\zeta_\delta$ and $\ell'$ is surjective, v.\ \cite[Lem.\ 4.11]{Matschke}.
\end{remark}
 
\begin{remark} Properties (ii) and (v) yield
    the short exact sequence 
    $0 \to S^{\alpha\alpha}_{\alpha\alpha}\xrightarrow[]{\id}S^{\alpha\alpha}_{\alpha\alpha}\xrightarrow[]{\id}S^{\alpha\alpha}_{\alpha\alpha}\to 0$, and therefore $S^{\alpha\alpha}_{\alpha\alpha}= 0$ for all $\alpha$. By \eqref{homdef1} this implies that $S_{\alpha\beta}^{\beta\gamma}= 0$ for all $(\alpha,\beta,\gamma)\in \bm{\sQ_3}$.
\end{remark}

\subsection{The $E$-terms for spectral systems}
From Definition~\ref{specsysdefn}(iii) it follows that
\begin{equation}
    \label{kerandcoker}
\ker\colon S_{\gamma\delta}^{\epsilon\zeta}\xrightarrow[]{\dff} S_{\alpha\beta}^{\gamma\delta}  \cong S^{\epsilon\zeta}_{\beta\delta},\quad\text{and}\quad \coker \colon
S_{\epsilon\zeta}^{\eta\theta}\xrightarrow[]{\dff} S_{\gamma\delta}^{\epsilon\zeta} \cong S_{\gamma\delta}^{\epsilon \eta},
\end{equation}
 for all octuples $(\alpha,\beta,\gamma,\delta,\epsilon,\zeta,\eta,\theta)\in \bm{\sQ_8}$, v.\ \cite[Lem.\ 4.11]{Matschke}.
In particular, for a triple $(\alpha,\beta,\gamma)\in \bm{\sQ_3}$, the differential $\dff\colon S^{\gamma\gamma}_{\beta\beta}\to S^{\beta\beta}_{\alpha\alpha}$  is defined via the sextuple $(\alpha,\alpha,\beta,\beta,\gamma,\gamma)$. Then, 
the short exact sequences in \eqref{shortexactSt} yield the short exact sequence
\begin{equation}
    \label{specialsheseq1}
\begin{tikzcd}
0 \arrow[r] & S^{\beta\gamma}_{\alpha\alpha} \arrow[r,"i"] & S^{\gamma\gamma}_{\alpha\alpha} \arrow[r,"j"] & S^{\gamma\gamma}_{\alpha\beta} \arrow[r] & 0,
\end{tikzcd}
\end{equation}
with $S^{\beta\gamma}_{\alpha\alpha}\cong \coker~\dff$
and $S^{\gamma\gamma}_{\alpha\beta}\cong \ker\dff$, by \eqref{kerandcoker}.
Equation \eqref{kerandcoker} implies another family of exact sequences.
\begin{lemma}
    \label{exactinSpS}
    The sequence
    \begin{equation}
    \label{specialsheseq2111}
\begin{tikzcd}
0 \arrow[r] & S_{\beta\delta}^{\epsilon\zeta} \arrow[r,"\ell"] & S^{\epsilon\zeta}_{\gamma\delta} \arrow[r,"\dff"] &   S^{\gamma\delta}_{\alpha\beta} \arrow[r,"\ell"] &S^{\gamma\epsilon}_{\alpha \beta} \arrow[r] & 0
\end{tikzcd}
\end{equation}
     is exact for all  sextuples $(\alpha,\beta,\gamma,\delta,\epsilon,\zeta)\in \bm{\sQ_6}$.
\end{lemma}
\begin{proof}
For $\dff\colon S^{\epsilon\zeta}_{\gamma\delta}\to S^{\gamma\delta}_{\alpha\beta}$ we have the exact sequence
$0 \longrightarrow \ker \dff \xrightarrow{\ \subset\ } S^{\epsilon\zeta}_{\gamma\delta} \xrightarrow{\ \dff\ } S^{\gamma\delta}_{\alpha\beta} \xrightarrow{\ \pi\ } \coker~\dff \longrightarrow 0$.
By Definition~\ref{specsysdefn}(iii), $\ell\colon S_{\beta\delta}^{\epsilon\zeta}\to S^{\epsilon\zeta}_{\gamma\delta}$ is injective with image $\ker\dff$, and $\ell\colon S^{\gamma\delta}_{\alpha\beta}\to S^{\gamma\epsilon}_{\alpha \beta}$ is surjective with kernel $\image\dff$ (apply (iii) to the octuple $(\alpha,\alpha,\alpha,\beta,\gamma,\delta,\epsilon,\zeta)$), which proves the lemma.
    \qed
\end{proof}

As a consequence of Lemma \ref{exactinSpS} we can build   exact triangles from the short exact sequences in \eqref{specialsheseq1}.
\begin{lemma}
    \label{extraextrian}
Let $\bfS$ be a 
spectral system over $\sQ$. Then,   
the triangles
\begin{equation}
    \label{theexacttris}
  \begin{tikzcd}[column sep=0.5cm]
S^{\beta\beta}_{\alpha\alpha} \arrow[rr, "i"] && S^{\gamma\gamma}_{\alpha\alpha} \arrow[ld, "j"] \\
& S^{\gamma\gamma}_{\beta\beta} \arrow[lu, "\dff"]
\end{tikzcd}
\end{equation}
are exact, for all triples $(\alpha,\beta,\gamma)\in \bm{\sQ_3}$.
\end{lemma}
\begin{proof}
By Lemma \ref{exactinSpS} the sequence
    \begin{equation}
    \label{specialsheseq2}
\begin{tikzcd}
0 \arrow[r] & S^{\gamma\gamma}_{\alpha\beta} \arrow[r,"j''"] & S^{\gamma\gamma}_{\beta\beta} \arrow[r,"\dff"] &   S^{\beta\beta}_{\alpha\alpha} \arrow[r,"i''"] &S^{\beta\gamma}_{\alpha\alpha} \arrow[r] & 0
\end{tikzcd}
\end{equation}
     is exact for all  triples $(\alpha,\beta,\gamma)\in \bm{\sQ_3}$.
     By combining \eqref{specialsheseq1} with \eqref{specialsheseq2} we obtain the 
 diagram
\[
\begin{tikzcd}[row sep=large, column sep=large]
            & 0                                                            &                                                                & 0 \arrow[d]                                             &   \\
0 \arrow[r] & S^{\beta\gamma}_{\alpha\alpha} \arrow[r, "i'"] \arrow[u]      & S^{\gamma\gamma}_{\alpha\alpha} \arrow[r, "j'"] \arrow[rd, "j"] & S^{\gamma\gamma}_{\alpha\beta} \arrow[r] \arrow[d, "j''"] & 0 \\
            & S^{\beta\beta}_{\alpha\alpha} \arrow[u, "i''"] \arrow[ru, "i"] &                                                                & S^{\gamma\gamma}_{\beta\beta} \arrow[ll, "\dff"]        &  
\end{tikzcd}
\]
Using   the composition property in (ii) we obtain $i=i'i''$ and $j=j''j'$. 
Since the horizontal sequence and $U$-shaped sequence from Lemma \ref{exactinSpS} are exact
we have
\[
\begin{aligned}
&\ker i'= 0,\quad \ker j'=\image i',\quad\text{and}\quad \image j'=S^{\gamma\gamma}_{\alpha\beta};\\
&\ker j''=0,\quad \ker\dff=\image j'',\quad \ker i''=\image \dff,\quad\text{and}\quad \image i''=S^{\beta\gamma}_{\alpha\alpha}.
\end{aligned}
\]
Now consider $\ker i$, i.e.\ $i'i''(x) =0$ if and only if $i''(x)=0$ since $i'$ is injective.
Therefore, $\ker i = \ker i'i''=\ker i'' = \image\dff$. Similarly,  $\ker j = \ker j''j'$, i.e.\ $j''j'(x)=0$ if and only if $j'(x) =0$ since $j''$ is injective. Therefore,
$\ker j = \ker j''j' = \ker j' = \image i' = \image i'i'' = \image i$, since $i''$ is surjective. Finally, $\ker \dff = \image j'' = \image j''j' =\image j$, since $j'$ is surjective, which proves 
 that the triangle is exact. 
    \qed
\end{proof}

For a spectral system $\bfS$ over $\sQ$ we define the $E$-terms $E^\beta_\alpha = S^{\beta\beta}_{\alpha\alpha}$, which then define a Cartan--Eilenberg system $\bfE(\bfS)$ over $\sQ$, with connecting homomorphisms $k=\dff\colon S^{\gamma\gamma}_{\beta\beta}\to S^{\beta\beta}_{\alpha\alpha}$. From Definition \ref{specsysdefn}(iii) it follows  that
\[
S^{\gamma\delta}_{\alpha\beta} \cong \frac{\ker\colon S_{\beta\beta}^{\gamma\gamma}\xrightarrow[]{\dff} S_{\alpha\alpha}^{\beta\beta} }{\image\colon S_{\gamma\gamma}^{\delta\delta}\xrightarrow[]{\dff} S_{\beta\beta}^{\gamma\gamma}} = \displaystyle  \frac{\ker \colon E^{\gamma}_{\beta} \xrightarrow{k} E^{\beta}_{\alpha}}{\image \colon E^{\delta}_{\gamma} \xrightarrow{k} E^{\gamma}_{\beta}} = H\bigl(S_{\beta\beta}^{\gamma\gamma}\bigr) = H\bigl(E^{\gamma}_{\beta}\bigr),
\]
by using octuples $(\alpha,\alpha,\beta,\beta,\gamma,\gamma,\delta,\delta)\in \bm{\sQ_8}$. The conditions for  Cartan--Eilenberg systems follow directly from Definition \ref{specsysdefn} and Lemma~\ref{extraextrian}.

\begin{example}
\label{linearorderSterms}
    In a linear order $\sP$ given by $1 \covv 2 \covv\cdots {\scriptstyle <} p\covv \cdots\covv n$, 
    where $n=\#\sP$,
    the lattice of down-sets $\sQ=\sO(\sP)$ is given by
    $\bt \covv \alpha_1 \covv \alpha_2\covv \cdots \covv \alpha_{p}\covv\cdots\covv\alpha_{n-1}\covv \tp$,\footnote{We say that $p$ \emph{covers} $q$, notation $q\prec p$, if $q<p$ and $q\le q' \le p$ implies that $q'=q$, or $q'=p$.}
    where $\alpha_p=\{1,\cdots,p\}$, $\bt=\varnothing$, and $\tp = \{1,\cdots,n\}$.
    For technical reasons it is convenient to take $p\in \Z$ for indexing the $S$-terms and define
    \[
    \alpha_p = \begin{cases}
        \bt & \text{for~~} p\le 0,\\
        \tp & \text{for~~} p \ge n.
    \end{cases}
    \]
    The $S$-terms  $S^{\alpha_p\alpha_{p-1+r}}_{\alpha_{p-r}\alpha_{p-1}}$ are the standard terms in a spectral sequence, where $r$ is the page number.
\end{example}

\subsection{The $S$-terms for Cartan--Eilenberg systems}
For a Cartan--Eilenberg system $\bfE$ over $\sQ$ we follow \cite{Matschke} for the construction of an associated spectral system $\bfS(\bfE)$. Define the $S$-terms:
\begin{equation}
    \label{Sterms}
S_{\alpha\beta}^{\gamma\delta}=\displaystyle  \frac{\ker \colon E^{\gamma}_{\beta} \xrightarrow{k} E^{\beta}_{\alpha}}{\image \colon E^{\delta}_{\gamma} \xrightarrow{k} E^{\gamma}_{\beta}},\quad \forall (\alpha,\beta,\gamma,\delta)\in \bm{\sQ_4}.
\end{equation}
Here $k$ denotes the connecting homomorphisms of $\bfE$.
For $(\alpha,\beta,\gamma,\delta)\le(\alpha',\beta',\gamma',\delta')$ in $\bm{\sQ_4}$
the  homomorphisms $\ell\colon E^\gamma_\beta \to E^{\gamma'}_{\beta'}$
induce   homomorphisms:
    \[
    \ell\colon S_{\alpha\beta}^{\gamma\delta}\rightarrow S_{\alpha'\beta'}^{\gamma'\delta'},
    \]
v.\ \cite[Sect.\ 4]{Matschke}. The connecting homomorphisms $k$ on the $E$-terms define differentials $\dff\colon S_{\gamma\delta}^{\epsilon\zeta} \to S_{\alpha\beta}^{\gamma\delta}$ on the $S$-terms for
sextuples $(\alpha,\beta,\gamma,\delta,\epsilon,\zeta)\in \bm{\sQ_6}$.
The differentials commute with the $\ell$-homomorphisms: $\ell\circ \dff=\dff\circ \ell$. 
The $S$-terms together with the above defined differentials yield a spectral system denoted by $\bfS(\bfE)$, v.\ \cite[Sect.\ 4]{Matschke}.
A \emph{morphism} $h\colon\bfS\to\bfS'$ of spectral systems over $\sQ$ consists of homomorphisms $h\colon S^{\gamma\delta}_{\alpha\beta}\to S'^{\gamma\delta}_{\alpha\beta}$ that commute with the homomorphisms $\ell$ and with the differentials $\dff$; $h$ is an \emph{isomorphism} if all these homomorphisms are isomorphisms.
Combining this with Lemma \ref{extraextrian} we obtain the following duality theorem:
\begin{theorem}
    \label{duality}
    Let $\bfS$ be a     spectral system and let $\bfE$ be a Cartan--Eilenberg system, both over a poset $\sQ$. Then,
    \[
    \bfS\bigl( \bfE(\bfS)\bigr) \cong \bfS,\quad\text{and}\quad \bfE\bigl( \bfS(\bfE)\bigr) = \bfE.
    \]
\end{theorem}
\begin{proof}
Let $\bfE$ be a Cartan--Eilenberg system. By \eqref{Sterms}, 
\[
S^{\beta\beta}_{\alpha\alpha}(\bfE)=\ker\bigl(E^\beta_\alpha\xrightarrow{k}E^\alpha_\alpha\bigr)/\image\bigl(E^\beta_\beta\xrightarrow{k}E^\beta_\alpha\bigr)=E^\beta_\alpha,
\]
since $E^\alpha_\alpha=E^\beta_\beta=0$, and the induced homomorphisms $\ell$ and differentials $\dff$ on these terms are the homomorphisms $\ell$ and $k$ of $\bfE$. This implies that $\bfE\bigl(\bfS(\bfE)\bigr)=\bfE$.

Let $\bfS$ be a spectral system and $\bfE=\bfE(\bfS)$, i.e.\ $E^\beta_\alpha=S^{\beta\beta}_{\alpha\alpha}$. For $(\alpha,\beta,\gamma,\delta)\in\bm{\sQ_4}$, Remark~\ref{homiso} with the octuple $(\alpha,\alpha,\beta,\beta,\gamma,\gamma,\delta,\delta)$ gives isomorphisms
\[
h\colon S^{\gamma\delta}_{\alpha\beta}\xrightarrow{\ \cong\ }H\bigl(S^{\gamma\gamma}_{\beta\beta}\bigr)=S^{\gamma\delta}_{\alpha\beta}(\bfE),\qquad x\mapsto\bigl[\ell(y)\bigr],
\]
where $y\in S^{\gamma\gamma}_{\alpha\beta}$ with $\ell(y)=x$. On the $E$-terms, $h$ is the identity. Since $h$ is defined by homomorphisms $\ell$, it commutes with the homomorphisms $\ell$, by (ii). For the differentials, (iv) applied to $(\alpha,\alpha,\gamma,\gamma,\epsilon,\epsilon)\le(\alpha,\beta,\gamma,\delta,\epsilon,\zeta)$ gives $\dff\circ\ell=\ell\circ k$ on $E^\epsilon_\gamma=S^{\epsilon\epsilon}_{\gamma\gamma}$, where $\ell\colon E^\epsilon_\gamma\to S^{\epsilon\zeta}_{\gamma\delta}$ is surjective by (v) and (iii). The same relation holds in $\bfS(\bfE)$, which satisfies (iv) and has the same $E$-terms and connecting homomorphisms $k$. Since $h$ commutes with $\ell$ and is the identity on $E$-terms, it follows that $h\circ\dff=\dff\circ h$. Hence $h\colon\bfS\to\bfS\bigl(\bfE(\bfS)\bigr)$ is an isomorphism of spectral systems.
\qed
\end{proof}
If  $\sQ=\sL$ is a  modular lattice, then a spectral system is said to be
\emph{excisive} if $\bfE(\bfS)$ is an excisive Cartan--Eilenberg system over $\sL$.
In terms of the $S$-terms this is equivalent, for $\sL=\sO(\sP)$, to
\[
\ell\colon S^{\gamma\delta}_{\alpha\beta} \xrightarrow{\cong} S^{\gamma'\delta'}_{\alpha'\beta'},\quad\forall
(\alpha,\beta,\gamma,\delta)\le(\alpha',\beta',\gamma',\delta'),
\]
with, $\beta\setminus\alpha = \beta'\setminus\alpha'$, $\gamma\setminus\beta = \gamma'\setminus\beta'$, and $\delta\setminus\gamma = \delta'\setminus\gamma'$, 
v.\ \cite[Lem.\ 4.16]{Matschke}.
In other words, in an excisive $S$-system the $S$-terms $S^{\gamma\delta}_{\alpha\beta}$ only depend on the convex sets $\beta\setminus\alpha,~\gamma\setminus\beta, \text{~and~} \delta\setminus\gamma$.

\subsection{Convergence and extensions}
\label{convergence}
A spectral system $\bfS$ over $\sQ$ may have various convergence properties. A poset $\sQ$ is \emph{bounded} if it contains unique elements $\bt$ and $\tp$ such that $\bt\le q\le \tp$ for all $q\in \sQ$.\footnote{Every poset $\sQ$ can be extended to a bounded poset $\sQ^*$ by adding a bottom and top element $\bt$ and $\tp$ respectively. In the setting of spectral systems it is convenient to consider bounded posets which we again denote by $\sQ$.} This allows the construction of  graded versions of the page-$1$ $E$-terms, v.\   \cite{Matschke}. 
In this setting the page-$\infty$ terms are given by $S_{\smbt\alpha}^{\beta \smtp}$.
The poset $\sQ$ induces the following filtering of $E^\smtp_\smbt$:
\[
\alpha \mapsto F_\alpha E^\smtp_\smbt = \image\colon E^\alpha_\smbt \xrightarrow[]{i} E^\smtp_\smbt, \quad\text{with}\quad 
\alpha\le\beta \implies F_\alpha E^\smtp_\smbt\subset F_\beta E^\smtp_\smbt,
\]
which follows from the composition of the  homomorphisms $i$, v.\ Defn.~\ref{CEsys}(ii). 
The poset filtering of $E_\smbt^\smtp$ yields graded versions of $E_\smbt^\smtp$.
Following \cite{Matschke} the spectral system yields the short exact sequences, for $\alpha\le\beta$:
\[
\begin{tikzcd}[row sep=large, column sep=large]
0 \arrow[r] & S^{\alpha\smtp}_{\smbt\smbt} \arrow[r] \arrow[d, "\cong"'] & S^{\beta\smtp}_{\smbt\smbt} \arrow[r] \arrow[d, "\cong"] & S^{\beta\smtp}_{\smbt\alpha} \arrow[r] \arrow[d, "\cong"] & 0 \\
0 \arrow[r] & F_\alpha E^\smtp_\smbt \arrow[r,"\subset"]                           & F_\beta E^\smtp_\smbt \arrow[r]                          & F_\beta E^\smtp_\smbt/F_\alpha E^\smtp_\smbt \arrow[r]    & 0
\end{tikzcd}
\]
A short exact sequence describes an extension $F_\beta E^\smtp_\smbt$ of $F_\beta E^\smtp_\smbt/F_\alpha E^\smtp_\smbt$ by $F_\alpha E^\smtp_\smbt$. Extensions are captured by the abelian group $\Ext^1_\bR\bigl(  S^{\beta\smtp}_{\smbt\alpha},S^{\alpha\smtp}_{\smbt\smbt} \bigr)$. For example if $S^{\beta\smtp}_{\smbt\alpha}$ is projective, then $\Ext^1_\bR\bigl(  S^{\beta\smtp}_{\smbt\alpha},S^{\alpha\smtp}_{\smbt\smbt} \bigr)=0$, and  therefore the exact sequence is split. This happens for example if we consider spectral systems in the category $\sKVec$.
A similar construction exists for the quotients $E^\smtp_\smbt/F_\alpha E^\smtp_\smbt$, for $\sQ=\sO(\sP)$.
Define
\[
\alpha \mapsto F_{\alpha^{\!\!\neg}} E^\smtp_\smbt = \coimage\colon E^\smtp_\smbt \xrightarrow[]{j} E^{\smtp}_{\alpha} = \coker \colon E^\alpha_\smbt \xrightarrow[]{i} E^\smtp_\smbt, \quad \alpha^{\!\!\neg}=\sP\setminus\alpha.
\]
Following \cite{Matschke} we obtain the short exact sequences, for $\alpha\le\beta$:
\[
\begin{tikzcd}[row sep=large, column sep=large]
0 \arrow[r] & S^{\beta\smtp}_{\smbt\alpha} \arrow[r] \arrow[d, "\cong"'] & S^{\smtp\smtp}_{\smbt\alpha} \arrow[r] \arrow[d, "\cong"] & S^{\smtp\smtp}_{\smbt\beta} \arrow[r] \arrow[d, "\cong"] & 0 \\
0 \arrow[r] & F_\beta E^\smtp_\smbt/F_\alpha E^\smtp_\smbt \arrow[r]                           & F_{\alpha^{\!\!\neg}} E^\smtp_\smbt \arrow[r]                          & F_{\beta^{\!\neg}} E^\smtp_\smbt \arrow[r]    & 0
\end{tikzcd}
\]
\begin{remark}
The same filtering can be applied to any $E$-term. This way we obtain grading information about a given $E$-term. Given a quadruple 
$(\alpha,\beta,\gamma,\delta) \in \bm{\sQ_4}$, then
\[
\beta \mapsto F_\beta E_\alpha^\delta = \image\colon E_\alpha^\beta \xrightarrow[]{i} E_\alpha^\delta\cong S^{\beta\delta}_{\alpha\alpha},\quad\text{and}\quad F_\gamma E_\alpha^\delta/F_\beta E_\alpha^\delta \cong S_{\alpha\beta}^{\gamma\delta}.
\]
\end{remark}

Consider the short exact sequences
\begin{equation*}
\begin{tikzcd}
0 \arrow[r] & S_{\alpha\beta}^{\gamma\epsilon} \arrow[r, "i"] & S_{\alpha\beta}^{\delta\epsilon} \arrow[r, "j"] & S_{\alpha\gamma}^{\delta\epsilon} \arrow[r] & 0,
\end{tikzcd}
\end{equation*}
from \eqref{shortexactSt}. The extensions of $S_{\alpha\gamma}^{\delta\epsilon}$ by $S_{\alpha\beta}^{\gamma\epsilon}$ are described by the extension groups
 $\Ext^1_\bR\bigl( S_{\alpha\gamma}^{\delta\epsilon}, S_{\alpha\beta}^{\gamma\epsilon}\bigr)$, which are related by associated distinguished triangles in the derived category $\Dh(\sRmodR)$ (of complexes):
\begin{equation}
   \label{distshortexactSt}
\begin{tikzcd}
 S_{\alpha\beta}^{\gamma\epsilon} \arrow[r, "i"] & S_{\alpha\beta}^{\delta\epsilon} \arrow[r, "j"] & S_{\alpha\gamma}^{\delta\epsilon} \arrow[r,"\kappa"] & S_{\alpha\beta}^{\gamma\epsilon}[1],
\end{tikzcd}
\end{equation} 
where the $S$-terms are regarded as complexes concentrated at $n=0$, with trivial differentials, and $\kappa \in \Hom_{\Dh(\sRmodR)}\bigl( S_{\alpha\gamma}^{\delta\epsilon},S_{\alpha\beta}^{\gamma\epsilon}[1]\bigr)\cong \Ext^1_\bR\bigl( S_{\alpha\gamma}^{\delta\epsilon}, S_{\alpha\beta}^{\gamma\epsilon}\bigr)$ which is not induced by an actual homomorphism in general. 
If $S_{\alpha\gamma}^{\delta\epsilon}$ is projective, then $\kappa=0$. Otherwise,
to determine $\kappa$ we consider a projective resolution for $S_{\alpha\gamma}^{\delta\epsilon}$, denoted by $P_\bullet\xrightarrow[\cong]{\pi} S_{\alpha\gamma}^{\delta\epsilon}$,
which is a quasi-isomorphism and $S_{\alpha\gamma}^{\delta\epsilon}\in \Ch(\sRmodR)$, concentrated at $n=0$.
Consider the commutative diagram:
\[
\begin{tikzcd}[row sep=large]
\cdots \arrow[r] & P_n \arrow[r, "\dff_n"] \arrow[d, "0"'] & \cdots \arrow[r] & P_2 \arrow[r, "\dff_2"] \arrow[d, "0"'] & P_1 \arrow[r, "\dff_1"] \arrow[d, "f_1", dotted] & P_0 \arrow[d, "f_0", dotted] \arrow[r, "\pi"] & S_{\alpha\gamma}^{\delta\epsilon} \arrow[r] \arrow[d, equal] & 0 \\
\cdots \arrow[r] & 0 \arrow[r]                             & \cdots \arrow[r] & 0 \arrow[r]                             & S_{\alpha\beta}^{\gamma\epsilon} \arrow[r,"i"]                                    & S_{\alpha\beta}^{\delta\epsilon} \arrow[r,"j"]                                             & S_{\alpha\gamma}^{\delta\epsilon} \arrow[r]                    & 0
\end{tikzcd}
\]
Since $P_0$ is projective  and $j$ is surjective there exists a lift $f_0\colon P_0\to S_{\alpha\beta}^{\delta\epsilon}$, such that $jf_0=\pi$. Observe that $jf_0 \dff_1 = \pi \dff_1 =0$ and therefore $\image (f_0\dff_1) \subset \ker j = \image i$. The map $i$ is injective and since $f_0\dff_1$ maps into $\image i$, the map $f_1:= i^{-1}f_0\dff_1$ is well-defined. The map $f_1\colon P_1\to S_{\alpha\beta}^{\gamma\epsilon}$ fits in the commutative square above. Since $if_1\dff_2=f_0\dff_1\dff_2=0$ and $i$ is injective, $f_1\dff_2=0$, so that $f_1$ provides an actual chain map $P_\bullet\to S_{\alpha\beta}^{\gamma\epsilon}[1]$. In a diagram this reads:
\[
\begin{tikzcd}[row sep=large]
\cdots \arrow[r] & P_n \arrow[r, "\dff_n"] \arrow[d, "0"'] & \cdots \arrow[r] & P_2 \arrow[r, "\dff_2"] \arrow[d, "0"'] & P_1 \arrow[r, "\dff_1"] \arrow[d, "f_1"]   & P_0 \arrow[d, "0", no head] \\
\cdots \arrow[r] & 0 \arrow[r]                             & \cdots \arrow[r] & 0 \arrow[r]                             & S_{\alpha\beta}^{\gamma\epsilon} \arrow[r] & 0                          
\end{tikzcd}
\]
The morphism $\kappa\colon S_{\alpha\gamma}^{\delta\epsilon}\to S_{\alpha\beta}^{\gamma\epsilon}[1]$ is represented by the homotopy class of the chain map $f_1\colon P_\bullet\to S_{\alpha\beta}^{\gamma\epsilon}[1]$, composed with the inverse of the quasi-isomorphism $\pi$. As we pointed out, $\kappa$ represents a class in $\Ext^1_\bR(S_{\alpha\gamma}^{\delta\epsilon},S_{\alpha\beta}^{\gamma\epsilon})$,
and $\kappa=0$ if and only if $f_1\sim 0$, which is equivalent to $f_1-0=\sigma_0 \dff_1 + 0\cdot\sigma_0$, i.e.\ $f_1 = \sigma_0\dff_1$ for some $\sigma_0\colon P_0\to S_{\alpha\beta}^{\gamma\epsilon}$.
Consequently, $\kappa\neq 0$ if and only if the short exact sequence \eqref{shortexactSt} does not split.
The homomorphism (in $\Dh(\sRmodR)$) $\kappa\colon S_{\alpha\gamma}^{\delta\epsilon}\to S_{\alpha\beta}^{\gamma\epsilon}[1]$ will be referred to as a \emph{(higher order) derived differential} between $S$-terms.

More explicitly, for an arbitrary ring $\bR$ the following statements are equivalent:
\begin{enumerate}[label=\textup{(\roman*)},leftmargin=2.2em,itemsep=0.1em,topsep=0pt]
\item[(i)] $\kappa=0$;
\item[(ii)] $f_1=\sigma_0\dff_1$ for some $\sigma_0\colon P_0\to S_{\alpha\beta}^{\gamma\epsilon}$;
\item[(iii)] the short exact sequence in \eqref{shortexactSt} splits, i.e.\ there exists a homomorphism $s\colon S_{\alpha\gamma}^{\delta\epsilon}\to S_{\alpha\beta}^{\delta\epsilon}$ with $js=\id$, or, equivalently, a homomorphism $r\colon S_{\alpha\beta}^{\delta\epsilon}\to S_{\alpha\beta}^{\gamma\epsilon}$ with $ri=\id$.
\end{enumerate}
The equivalence of (i) and (ii) was shown above. The equivalence of (i) and (iii) holds since $\kappa$ represents the class of the extension \eqref{shortexactSt} in $\Ext^1_\bR\bigl(S_{\alpha\gamma}^{\delta\epsilon},S_{\alpha\beta}^{\gamma\epsilon}\bigr)$, and this class vanishes if and only if the extension splits, v.\ \cite[Thm.~3.4.3]{weibel}.

\begin{remark}
\label{projkappa}
    If $S^{\delta\epsilon}_{\alpha\gamma}$ is projective, then the projective resolution can be chosen trivial and $\kappa=0$ follows immediately. Indeed, $P_0\cong S^{\delta\epsilon}_{\alpha\gamma}$ and $P_n=0$, $\dff_n=0$ for $n\ge 1$. This implies that $f_1=0$.
\end{remark}

\section{Spectral nets}
\label{spectralnets}
In this paper we restrict ourselves to Cartan--Eilenberg systems over a finite distributive lattice $\sL$, which we represent, using Birkhoff's representation theorem, as the lattice $\sO(\sP)$ of down-sets of a finite poset $\sP$.
Throughout this section, $\bfE$ is an excisive Cartan--Eilenberg system over $\sO(\sP)$, v.\ Sect.~\ref{CAsystems}, and $\bfS=\bfS(\bfE)$ is the associated spectral system, v.\ Thm.~\ref{duality}. Since $\bfE$ is excisive also  $\bfS$ is excisive.
A way to obtain the \say{limiting} $S$-terms in an excisive spectral system is to
define a generalization of spectral sequences to lattices, which we will refer to as a \emph{spectral net}.
A \emph{maximal chain} $\sK$ in $\sO(\sP)$ is a linearly ordered subset of $\sO(\sP)$ (regarded as a poset) that is not contained in a larger chain in $\sO(\sP)$.
The length of a maximal chain $\sK$ is equal to the number of elements $|\sP|$ of $\sP$, since $\alpha\subsetneq\beta$ is a cover in $\sO(\sP)$ if and only if $\beta\setminus\alpha$ is a singleton.
Therefore all maximal chains in $\sO(\sP)$ have the same length $|\sP|$ (Jordan--Dedekind chain condition).
As shown below, the number of maximal chains in $\sO(\sP)$ is equal to the number of linear extensions of $\sP$.\footnote{The number of linear extensions can be determined via the order polynomial $\rmOmega_\sP(\omega)$ of $\sP$. The coefficient of the highest order term $\omega^k$, $k=|\sP|$, multiplied by $k!$ gives the number of maximal chains in $\sO(\sP)$, v.\ \cite{Stanley1}.}
Let $\sK\subset \sO(\sP)$ be a maximal chain. If $\sK$ is represented as $\sO(\sP')$, then the inclusion $\sO(\sP')\rightarrowtail\sO(\sP)$ is an injective lattice homomorphism. Since $\sK$ is a sublattice of $\sO(\sP)$ which contains $\bt$ and $\tp$,   the dual order-preserving map $\sP\twoheadrightarrow\sP'$ is surjective. Since $|\sP'|=|\sK|-1=|\sP|$, the map $\sP\twoheadrightarrow\sP'$ is bijective (in general not an order-isomorphism) and $\sP'$ is a linear extension of the order on $\sP$. Conversely, a linear extension $p_1,\ldots,p_n$ of $\sP$ yields the maximal chain of down-sets $\{p_1,\ldots,p_k\}$, $k=0,\ldots,n$. This gives a one-to-one correspondence between maximal chains in $\sO(\sP)$ and linear extensions of $\sP$.
Define the following exact couple:
\[
\begin{tikzcd}[row sep=large]
D^1 \arrow[rr, "i"] &                     & D^1 \arrow[ld, "j"] \\
                    & E^1 \arrow[lu, "k"] &
\end{tikzcd}
\quad
D^1=\bigoplus_{\sK\subset\sO(\sP)}\bigoplus_{\alpha\in \sK}E^\alpha_\smbt,~~~~
E^1=\bigoplus_{\sK\subset\sO(\sP)}\bigoplus_{\alpha\in \sK\setminus\{\smbt\}}E^\alpha_{\alpha^\pred}
\]
where $\sK$ ranges over the maximal chains in $\sO(\sP)$, and $\alpha^\pred$ is the unique immediate predecessor of $\alpha$ in $\sK$. The homomorphisms $i$, $j$ and $k$ are those of the exact triangles of $\bfE$ for the pairs $\alpha^\pred\subset\alpha$, v.\ Sect.~\ref{CAsystems}, i.e.\ $i\colon E^{\alpha^\pred}_\smbt\to E^\alpha_\smbt$, $j\colon E^\alpha_\smbt\to E^\alpha_{\alpha^\pred}$ and $k\colon E^\alpha_{\alpha^\pred}\to E^{\alpha^\pred}_\smbt$.\footnote{As in Example~\ref{linearorderSterms}, $\sK$ is indexed by $m\in\Z$, with $\alpha_m=\bt$ for $m\le0$ and $\alpha_m=\tp$ for $m\ge|\sP|$, the sums are taken over $m\in\Z$, and $i\colon E^{\alpha_{m-1}}_\smbt\to E^{\alpha_m}_\smbt$ is the identity for $m>|\sP|$. This makes the couple exact at $E^\smtp_\smbt$, and does not change $E^1$, since $E^\smtp_\smtp=0$.}
The composition $\dff^1= jk$ defines a (first order) differential $\dff^1\colon E^1\to E^1$, and by the theory of exact couples, consecutive evaluation of homology with respect to the induced differentials yields the higher pages and differentials, i.e.\ $E^{r+1} = H(E^r,\dff^r)$,  with $\dff^r=j\,i^{1-r}k$, where $i^{1-r}$ denotes taking a preimage under $i^{r-1}$, which is well defined on $E^r$, $r\ge 1$, cf.\ Sect.\ \ref{theN5} and \cite[Sect.~5.9]{weibel}.

Since $\bfE$ is excisive, there exist $\bR$-modules $G_pC\cong E^\beta_\alpha$, for all $\beta\setminus\alpha=\{p\}$ and for all $p\in \sP$, with isomorphisms that are compatible with the homomorphisms $\ell$. For example, 
$G_pC=E^{\beta}_{\alpha}$, where $\beta = \bigl\downarrow p$ is the principal down-set generated by $p$ and 
$\alpha=\beta^\pred = \bigl\downarrow p\setminus\{p\}$,
with the isomorphisms $\ell\colon E^\beta_\alpha\to E^{\beta'}_{\alpha'}$, v.\ Sect.~\ref{CAsystems}.
Define
\[
C^1 := \bigoplus_{\sP\twoheadrightarrow\sP'} \bigoplus_{p\in \sP'} G_pC,
\]
where the outer sum ranges over all linear extensions $\sP'$ of $\sP$ determined by the choice of maximal chain $\sK$.
The isomorphisms from the excisive property can be used to define a differential $\dff^1$ on $C^1$, as expressed in the diagram below:
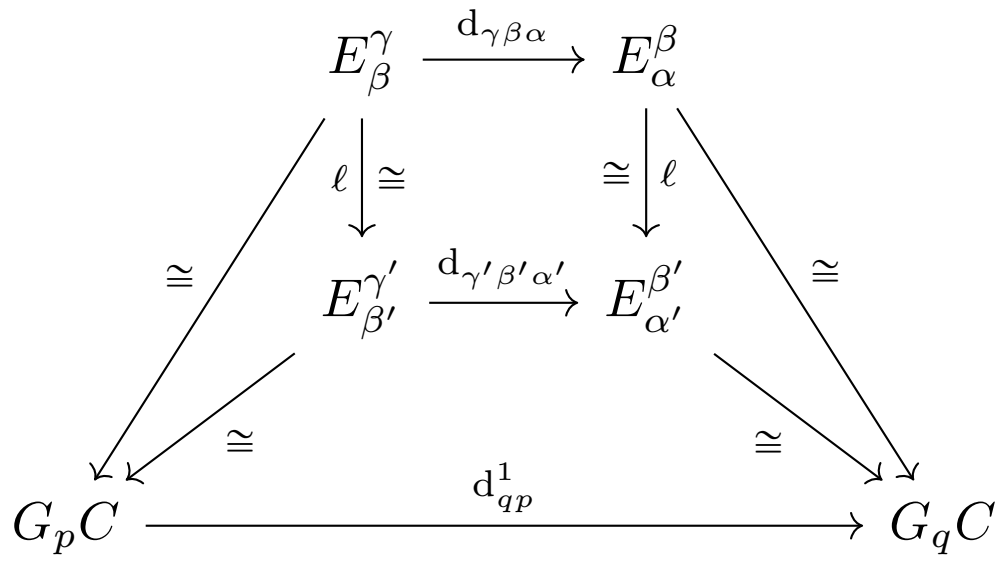
\begin{equation}
    \label{equivofdiff12}
\begin{tikzcd}[row sep=0.8cm, column sep=1cm]
                                & E^\gamma_\beta \arrow[r, "\dff_{\gamma\beta\alpha}"] \arrow[d, "\ell"',"\cong"] \arrow[ldd, "\cong"'] & E^\beta_\alpha \arrow[d, "\ell","\cong"'] \arrow[rdd, "\cong"] &      \\
                                & E^{\gamma'}_{\beta'} \arrow[r, "\dff_{\gamma'\beta'\alpha'}"] \arrow[ld, "\cong"]                          & E^{\beta'}_{\alpha'} \arrow[rd, "\cong"']                  &      \\
G_pC \arrow[rrr, "\dff^1_{qp}"] &                                                                                                   &                                                       & G_qC
\end{tikzcd}
\end{equation}
Here $\gamma\setminus\beta=\gamma'\setminus\beta'=\{p\}$, $\beta\setminus\alpha=\beta'\setminus\alpha'=\{q\}$ and $(\alpha,\beta,\gamma)\le(\alpha',\beta',\gamma')$. The square commutes by the naturality of the differentials, and the triangles commute since the isomorphisms with $G_pC$ and $G_qC$ are compatible with $\ell$. Two such triples are both below their join, which is again such a triple. Hence $\dff^1_{qp}$ does not depend on the choice of $(\alpha,\beta,\gamma)$.
The isomorphisms $G_pC\cong E^\alpha_{\alpha^\pred}$, with $\alpha\setminus\alpha^\pred=\{p\}$, yield an isomorphism $C^1\cong E^1$, under which $\dff^1$ on $C^1$ corresponds to $\dff^1=jk$ on $E^1$. Indeed, for consecutive elements $\alpha\subset\beta\subset\gamma$ of $\sK$, the composition $jk\colon E^\gamma_\beta\to E^\beta_\smbt\to E^\beta_\alpha$ is the differential $\dff_{\gamma\beta\alpha}$, by the naturality of the exact triangles. Via this isomorphism the higher pages and differentials of the exact couple define $C^r$ and $\dff^r$.
The result is the \emph{(fundamental) spectral net} $C^{r+1} = H(C^r,\dff^r)$ for $\bfS$, consisting of sub-quotients isomorphic to the $S$-terms and associated differentials. If $\sK$ is written as $\alpha_m$, $m\in\Z$, as in Example~\ref{linearorderSterms}, then the summand of $C^r$ that corresponds to $\alpha_m\setminus\alpha_{m-1}$, $1\le m\le|\sP|$, is isomorphic to $S^{\alpha_m\alpha_{m-1+r}}_{\alpha_{m-r}\alpha_{m-1}}$, cf.\ Sect.~\ref{theN5}.
 Since $\sO(\sP)$ is a finite
lattice, the spectral net converges to $C^\infty$, i.e.\ $C^r=C^\infty$ for $r\ge|\sP|$, which determines the $S$-terms $S^{\beta\smtp}_{\smbt\alpha}$ for all $\beta\setminus\alpha=\{p\}$ and for all $p
\in \sP$, and  computes the graded object of $E^\smtp_\smbt$ with respect to the filtrations $\alpha\mapsto F_\alpha E^\smtp_\smbt$, $\alpha\in\sK$, v.\ Sect.~\ref{convergence}. For the fundamental spectral net no requirements  on the ring $\bR$ are needed. By considering sublattices of $\sO(\sP)$ we obtain associated spectral nets that reveal other higher differentials and $E$-terms, v.\ Defn.~\ref{subnetdefn} below. Since a spectral net is defined over a partial order, there exist various additional relations between differentials, which are explained in detail in Section \ref{firstordinfo}.

\begin{definition}
\label{subnetdefn}
Let $\sO(\sQ)\subset\sO(\sP)$ be a sublattice, where, by Birkhoff's representation theorem, $\sQ$ is the poset of join-irreducible elements of the sublattice. The restriction of $\bfE$ to $\sO(\sQ)$ is an excisive Cartan--Eilenberg system over $\sO(\sQ)$, since meets and joins in $\sO(\sQ)$ are those in $\sO(\sP)$. The \emph{spectral net over $\sO(\sQ)$}, or \emph{subnet}, is the fundamental spectral net of this restriction, i.e.\ the above construction with $\sQ$ in place of $\sP$, where the elements of $\sO(\sQ)$ are regarded as elements of $\sO(\sP)$. For an interval $\sO(\sQ)=[\alpha,\beta]=\{\gamma\in\sO(\sP)\mid\alpha\subset\gamma\subset\beta\}$ it is called the \emph{spectral net of the $E$-term $E^\beta_\alpha$}.
\end{definition}
Let $\bt_{\sQ}$ and $\tp_\sQ$ be the bottom and top elements of $\sO(\sQ)$. The spectral net over $\sO(\sQ)$ converges to the graded object of $E^{\tp_\sQ}_{\bt_\sQ}$ with respect to the filtrations along the maximal chains in $\sO(\sQ)$. Its page~$1$ consists of the $E$-terms $E^\beta_\alpha$ for covers $\alpha\subsetneq\beta$ in $\sO(\sQ)$, and the differentials between them. An interval $[\alpha,\beta]$ is isomorphic to $\sO(\beta\setminus\alpha)$ via $\gamma\mapsto\gamma\setminus\alpha$, i.e.\ $\sQ\cong\beta\setminus\alpha$, which is a convex subset of $\sP$. In this case page~$1$ consists of the $G_pC$, $p\in\beta\setminus\alpha$, with the differentials $\dff^1_{qp}$ of the fundamental spectral net. For other sublattices, $\beta\setminus\alpha$ need not be a singleton for a cover $\alpha\subsetneq\beta$ in $\sO(\sQ)$, and page~$1$ contains $E$-terms and differentials which do not occur on page~$1$ of the fundamental spectral net.

\section{Properties of first and  higher order differentials}
\label{properties}
In this section we look at various properties of higher differentials and spectral matrices, and their interpretation for Morse representations. In this section the results are stated in the category of $\bR$-modules.
\subsection{Perspective pairs}
\label{firstordinfo}
In \cite{fran} entries in a connection matrix can sometimes be computed using the fact that a connection matrix, or spectral matrix, is a differential. In principle a  spectral matrix recovers all information about the spectral system.
The interpretation of the entries of connection matrices  in \cite{fran} only applies to the flow-order of 
 adjacent Morse sets.
In this section we demonstrate that such information can also be obtained from the spectral system as properties of a spectral net. 
\begin{figure}[h] 
  \centering
\begin{tikzpicture}[
  scale=0.7,
  every node/.style={draw, rounded corners, fill=white, inner sep=2pt, font=\small},
  tick/.style={postaction={decorate}, decoration={
      markings, mark=at position 0.5 with {\draw (0,-2pt)--(0,2pt);}}},
  ticktick/.style={postaction={decorate}, decoration={
      markings,
      mark=at position 0.5 with {\draw (-1.5pt,-2pt)--(-1.5pt,2pt);
                                  \draw ( 1.5pt,-2pt)--( 1.5pt,2pt);}}}
]
  \node[circle] (La) at (0,0)    {$~\alpha~$};
  \node[circle] (Lb) at (-1.6,1.6) {$~\beta~$};
  \node[circle] (Lg) at ( 1.6,1.6) {$~\gamma~$};
  \node[circle] (Ld) at (0,3.2)  {$~\delta~$};

  \draw[tick] (La) -- (Lg);
  \draw[tick] (Lb) -- (Ld);
  \draw[ticktick] (La) -- (Lb);
  \draw[ticktick] (Ld) -- (Lg);

  \begin{scope}[xshift=7cm]
    \node (Ra)  at (0,0)     {$~\alpha~$};
    \node (Rb)  at (-1.6,1.6){$~\beta~$};
    \node (Rg)  at ( 1.6,1.6){$~\gamma~$};
    \node (Rd)  at (0,3.2)   {$~\delta~$};
    \node (Rdp) at (0,4.6)   {$\delta_{+}$};
    \node (Ram) at (0,-1.4)  {$\alpha_{-}$};

    \draw (Ra) -- (Rb);
    \draw (Ra) -- (Rg);
    \draw (Rb) -- (Rd);
    \draw (Rd) -- (Rg);
    \draw (Rd) -- (Rdp);
    \draw (Ra) -- (Ram);
  \end{scope}
\end{tikzpicture}
\caption{A perspective pair  $\beta\setminus\alpha=\delta\setminus\gamma$ and $\gamma\setminus\alpha=\delta\setminus\beta$ [left] and the embedding into $\sO(\sP)$ [right].
}\label{expair123}
\vspace{-2em}
\end{figure}

A \emph{perspective pair} in   $\sO(\sP)$ is a pair of triples
$(\alpha,\beta,\delta), (\alpha,\gamma,\delta)\in \bm{\sO(\sP)_3}$ such that $\beta\setminus\alpha=\delta\setminus\gamma$ and $\gamma\setminus\alpha=\delta\setminus\beta$, which is equivalent to
$\alpha=\beta\cap\gamma$ and $\delta = \beta\cup \gamma$, v.\ Fig.\ \ref{expair123}. We assume in addition that there exist \emph{tail elements} $\alpha_-\subsetneq\alpha$ and $\delta\subsetneq\delta_+$ in $\sO(\sP)$.
We have the following exact triangles:
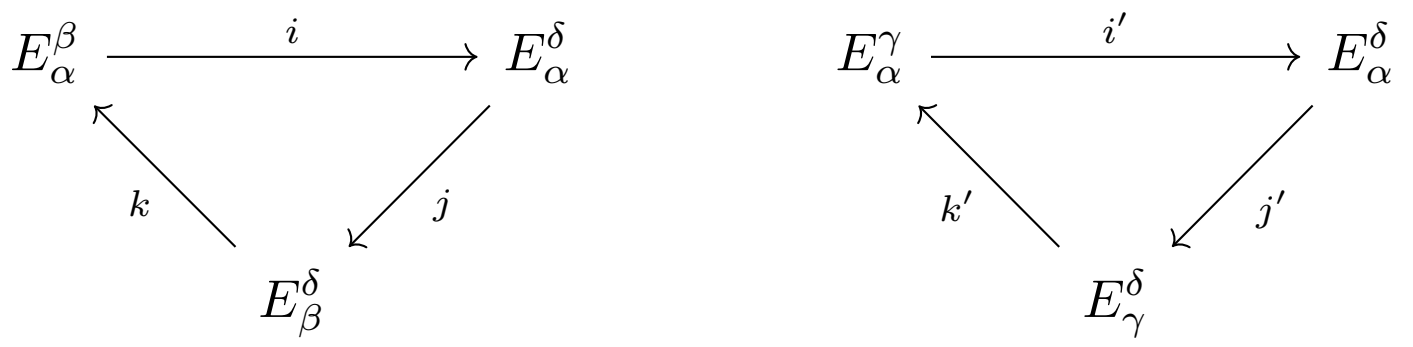
\[
\begin{tikzcd}[row sep=large]
E^{\beta}_{\alpha} \arrow[rr, "i"] &                                    & E^{\delta}_{\alpha} \arrow[ld, "j"] \\
                                   & E^{\delta}_{\beta} \arrow[lu, "k"] &                                    
\end{tikzcd}
\qquad\qquad
\begin{tikzcd}[row sep=large]
E^{\gamma}_{\alpha} \arrow[rr, "i'"] &                                      & E^{\delta}_{\alpha} \arrow[ld, "j'"] \\
                                     & E^{\delta}_{\gamma} \arrow[lu, "k'"] &                                     
\end{tikzcd}
\]
and isomorphisms (due to the excisive property):
 \[
\begin{aligned}
\ell  &: E^{\beta}_{\alpha} \xrightarrow{\;\cong\;} E^{\delta}_{\gamma}, &
\qquad j'i   &: E^{\beta}_{\alpha} \longrightarrow E^{\delta}_{\gamma};\\
\ell' &: E^{\gamma}_{\alpha} \xrightarrow{\;\cong\;} E^{\delta}_{\beta}, &
\qquad ji'   &: E^{\gamma}_{\alpha} \longrightarrow E^{\delta}_{\beta}.
\end{aligned}
\]
By the functoriality of the homomorphisms $\ell$, v.\ Defn.~\ref{CEsys}(ii), we have $j'i=\ell$ and $ji'=\ell'$, i.e.
\[
\begin{aligned}
\ell^{-1}\,j'i &= \id \text{~ on~ } E^{\beta}_{\alpha}, & \text{~and~~~~~}
\ell'^{-1}\,ji' &= \id \text{~ on~ } E^{\gamma}_{\alpha},\\
j'i\,\ell^{-1} &= \id \text{ ~on~ } E^{\delta}_{\gamma}, & \text{~and~~~~~}
ji'\,\ell'^{-1} &= \id \text{~ on~ } E^{\delta}_{\beta},
\end{aligned}
\]
v.\ \cite{SpV4}. Consequently, $i,i'$ are injective and $j,j'$ are surjective. Hence we obtain split short exact sequences
\[
\begin{tikzcd}[column sep=large]
0 \arrow[r] & E^{\beta}_{\alpha} \arrow[r, "i"] & E^{\delta}_{\alpha} \arrow[r, "j", two heads] & E^{\delta}_{\beta} \arrow[l, "i'\ell'^{-1}", tail, dashed, bend left, shift left] \arrow[r] & 0
\end{tikzcd}
\]
\[
\begin{tikzcd}[column sep=large]
0 \arrow[r] & E^{\gamma}_{\alpha} \arrow[r, "i'"] & E^{\delta}_{\alpha} \arrow[r, "j'", two heads] & E^{\delta}_{\gamma} \arrow[l, "i\ell^{-1}", tail, dashed, bend left, shift left] \arrow[r] & 0
\end{tikzcd}
\]
For the connecting homomorphisms $k$ and $k'$ this implies: $k=0$ and $k'=0$.
For the differentials we have the following anti-commutativity relation:
\begin{theorem}
    \label{firstantocomm}
    Consider a perspective pair  given by $(\alpha,\beta,\delta), (\alpha,\gamma,\delta)\in \bm{\sO(\sP)_3}$ with tail elements $\alpha_-$ and $\delta_+$.
    Then, the differentials $\dff_{\delta_+\delta\beta}\colon E^{\delta_+}_\delta\to E^\delta_\beta$ and $\dff_{\delta_+\delta\gamma}\colon E^{\delta_+}_\delta\to E^\delta_\gamma$, and $\dff_{\beta\alpha\alpha_-}\colon E^\beta_\alpha\to E^\alpha_{\alpha_-}$ and $\dff_{\gamma\alpha\alpha_-}\colon E^\gamma_\alpha\to E^\alpha_{\alpha_-}$ satisfy
    \begin{equation}
    \label{firstantocomm12}
        \dff_{\beta\alpha\alpha_{-}}\,\ell^{-1}\,\dff_{\delta_{+}\delta\gamma}
\;+\;
\dff_{\gamma\alpha\alpha_{-}}\,\ell'^{-1}\,\dff_{\delta_{+}\delta\beta}
\;=\;0
    \end{equation}
\end{theorem}
\begin{proof}
Consider the differentials
\[
\begin{tikzcd}[column sep=large]
E^{\delta_{+}}_{\delta} \arrow[r,"\dff_{\delta_{+}\delta\alpha}"] &
E^{\delta}_{\alpha} \arrow[r,"\dff_{\delta\alpha\alpha_{-}}"] &
E^{\alpha}_{\alpha_{-}},
\end{tikzcd}
\quad \text{with}\quad\dff_{\delta\alpha\alpha_{-}}\, \dff_{\delta_{+}\delta\alpha}=0
\]
The two factorizations of $E^{\delta}_{\alpha}$ give the diagrams
\[
\begin{tikzcd}[row sep=large, column sep=huge]
E^{\delta_{+}}_{\delta} \arrow[r,"\dff_{\delta_{+}\delta\alpha}"] \arrow[d,equal] &
  E^{\delta}_{\alpha} = i E^{\beta}_{\alpha}\oplus i'\ell'^{-1}E^{\delta}_{\beta}
  \arrow[d,"j", two heads] \\
E^{\delta_{+}}_{\delta} \arrow[r,"\dff_{\delta_+\delta\beta}"'] & E^{\delta}_{\beta}
\end{tikzcd}
\quad
\begin{tikzcd}[row sep=large, column sep=huge]
E^{\delta_{+}}_{\delta} \arrow[r,"\dff_{\delta_{+}\delta\alpha}"] \arrow[d,equal] &
  E^{\delta}_{\alpha} = i' E^{\gamma}_{\alpha}\oplus i\ell^{-1}E^{\delta}_{\gamma}
  \arrow[d,"j'", two heads] \\
E^{\delta_{+}}_{\delta} \arrow[r,"\dff_{\delta_+\delta\gamma}"'] & E^{\delta}_{\gamma}
\end{tikzcd}
\]
so that $j\,\dff_{\delta_{+}\delta\alpha}=\dff_{\delta_+\delta\beta}$ and
$j'\,\dff_{\delta_{+}\delta\alpha}=\dff_{\delta_+\delta\gamma}$.
Note that $i'\ell'^{-1}E^\delta_\beta=i'E^\gamma_\alpha$ and $i\ell^{-1}E^\delta_\gamma=iE^\beta_\alpha$, so that both splittings are the decomposition $E^\delta_\alpha=iE^\beta_\alpha\oplus i'E^\gamma_\alpha$.

\medskip
\noindent Write
$\displaystyle \dff_{\delta_{+}\delta\alpha}=\begin{pmatrix}\dff_{1}\\ \dff_{2}\end{pmatrix}$
with respect to the splitting $E^{\delta}_{\alpha}=iE^{\beta}_{\alpha}\oplus i'\ell'^{-1}E^{\delta}_{\beta}$.
Then
\[
\dff_{\delta_+\delta\beta}x \;=\; j\,\dff_{\delta_{+}\delta\alpha}x
            \;=\; j\bigl(\dff_{1}x+\dff_{2}x\bigr) \;=\; j\,\dff_{2}x,
\]
since the splitting
\[
\begin{tikzcd}[column sep=large]
0\arrow[r] & iE^{\beta}_{\alpha} \arrow[r,tail] &
E^{\delta}_{\alpha} \arrow[r,two heads] &
i'\ell'^{-1}E^{\delta}_{\beta} \arrow[r] & 0
\end{tikzcd}
\]
identifies the projection.

\medskip
\noindent The map $j\colon i'\ell'^{-1}E^{\delta}_{\beta}\xrightarrow{\;\cong\;}E^{\delta}_{\beta}$
has inverse $j^{-1}=i'\ell'^{-1}$, hence
\[
\dff_{2} \;=\; j^{-1}\dff_{\delta_+\delta\beta} \;=\; i'\ell'^{-1}\,\dff_{\delta_+\delta\beta}.
\]
Symmetrically (using the other splitting and $j'$),
$\dff_{1}=i\ell^{-1}\,\dff_{\delta_+\delta\gamma}$. Therefore
\[
\dff_{\delta_{+}\delta\alpha}
   \;=\;\begin{pmatrix} i\ell^{-1}\,\dff_{\delta_+\delta\gamma}\\[2pt] i'\ell'^{-1}\,\dff_{\delta_+\delta\beta}\end{pmatrix}.
\]
As for the second differential we argue as follows. 
Write
$\dff_{\delta\alpha\alpha_-}=\begin{pmatrix}\dff'_{1} ~~\dff'_{2}\end{pmatrix}$
with respect to the splitting $E^{\delta}_{\alpha}=iE^{\beta}_{\alpha}\oplus i'\ell'^{-1}E^{\delta}_{\beta}$.
As before we have the 
commutative diagrams:
\[
\begin{tikzcd}[row sep=large, column sep=huge]
i E^{\beta}_{\alpha}\oplus i'\ell'^{-1}E^{\delta}_{\beta}=E^{\delta}_{\alpha} \arrow[r,"\dff_{\delta\alpha\alpha_{-}}"] & E^{\alpha}_{\alpha_{-}} \\
E^{\beta}_{\alpha} \arrow[u,"i",tail] \arrow[r,"\dff_{\beta\alpha\alpha_-}"'] & E^{\alpha}_{\alpha_{-}} \arrow[u,equal]
\end{tikzcd}
\quad
\begin{tikzcd}[row sep=large, column sep=huge]
i' E^{\gamma}_{\alpha}\oplus i\ell^{-1}E^{\delta}_{\gamma}=E^{\delta}_{\alpha} \arrow[r,"\dff_{\delta\alpha\alpha_{-}}"] & E^{\alpha}_{\alpha_{-}} \\
E^{\gamma}_{\alpha} \arrow[u,"{i'}",tail] \arrow[r,"\dff_{\gamma\alpha\alpha_-}"'] & E^{\alpha}_{\alpha_{-}} \arrow[u,equal]
\end{tikzcd}
\]
we get $\dff_{\delta\alpha\alpha_{-}}\,i=\dff_{\beta\alpha\alpha_-}$, i.e.
$\dff_{\beta\alpha\alpha_-}x=\dff_{\delta\alpha\alpha_{-}}i(x)=\dff'_{1}i(x)$.
Since $i\colon E^{\beta}_{\alpha}\xrightarrow{\cong}iE^{\beta}_{\alpha}$,
\[
\dff_{\beta\alpha\alpha_-}\,i^{-1}=\dff'_{1},\qquad
\dff_{\gamma\alpha\alpha_-}\,i'^{-1}=\dff'_{2},
\]
so
\[
\dff_{\delta\alpha\alpha_{-}}\;=\;\bigl(\dff_{\beta\alpha\alpha_-}\,i^{-1}~~~\dff_{\gamma\alpha\alpha_-}\,i'^{-1}\bigr).
\]
The relation $\dff_{\delta\alpha\alpha_{-}}\,\dff_{\delta_{+}\delta\alpha}=0$ yields
\[
\bigl(\dff_{\beta\alpha\alpha_-}\,i^{-1}~~~\;\dff_{\gamma\alpha\alpha_-}\,i'^{-1}\bigr)
\begin{pmatrix} i\ell^{-1}\,\dff_{\delta_+\delta\gamma}\\[2pt] i'\ell'^{-1}\,\dff_{\delta_+\delta\beta}\end{pmatrix}
\;=0,
\]
and therefore
\[
\dff_{\beta\alpha\alpha_{-}}\,\ell^{-1}\,\dff_{\delta_{+}\delta\gamma}
\;+\;
\dff_{\gamma\alpha\alpha_{-}}\,\ell'^{-1}\,\dff_{\delta_{+}\delta\beta}
\;=\;0,
\]
which proves the theorem.
    \qed
\end{proof}

\begin{figure}
\centering
\begin{tikzcd}[row sep=1.25cm, column sep=0.9cm, rvpath/.style={densely dotted, semithick, no head},
  execute at end picture={
    \begin{scope}[densely dotted, semithick]
    \draw[red] (r1a.center) -- (r1b.center); \draw[red] (r1c.center) -- (r1d.center);
    \draw[red] (r2a.center) -- (r2b.center); \draw[red] (r2c.center) -- (r2d.center);
    \draw[blue] (b1a.center) -- (b1b.center); \draw[blue] (b1c.center) -- (b1d.center);
    \draw[blue] (b2a.center) -- (b2b.center); \draw[blue] (b2c.center) -- (b2d.center);
    \end{scope}}]
 & & G_3C \arrow[rd, equal] \arrow[r, "0"]
     \arrow[rd, red, rvpath, shift right=5pt, ""{description, name=r1b, pos=0, inner sep=0pt, outer sep=0pt}, ""{description, name=r1c, pos=1, inner sep=0pt, outer sep=0pt}]
   & G_4C \arrow[r, "\dff^1_{24}"]
     \arrow[r, red, rvpath, -{Stealth[length=4pt]}, shift right=4pt, ""{description, name=r2d, pos=0, inner sep=0pt, outer sep=0pt}]
   & G_2C \arrow[rd, "\dff^1_{12}"]
     \arrow[rd, blue, rvpath, -{Stealth[length=4pt]}, shift right=4pt, ""{description, name=b2d, pos=0, inner sep=0pt, outer sep=0pt}]
   & & \\
0 \arrow[r] & G_5C \arrow[ru, "\dff^1_{35}"] \arrow[r, "\dff^1_{45}"]
     \arrow[ru, red, rvpath, shift right=4pt, ""{description, name=r1a, pos=1, inner sep=0pt, outer sep=0pt}]
     \arrow[r, red, rvpath, shift right=4pt, ""{description, name=r2a, pos=1, inner sep=0pt, outer sep=0pt}]
   & G_4C \arrow[ru, equal] \arrow[r, "0"] \arrow[rd, "\dff^1_{24}"]
     \arrow[ru, red, rvpath, shift right=5pt, ""{description, name=r2b, pos=0, inner sep=0pt, outer sep=0pt}, ""{description, name=r2c, pos=1, inner sep=0pt, outer sep=0pt}]
     \arrow[r, blue, rvpath, shift right=4pt, ""{description, name=b1a, pos=1, inner sep=0pt, outer sep=0pt}]
     \arrow[rd, blue, rvpath, shift right=4pt, ""{description, name=b2a, pos=1, inner sep=0pt, outer sep=0pt}]
   & G_3C \arrow[rd, equal] \arrow[ru, "0"]
     \arrow[ru, red, rvpath, -{Stealth[length=4pt]}, shift left=4pt, ""{description, name=r1d, pos=0, inner sep=0pt, outer sep=0pt}]
     \arrow[rd, blue, rvpath, shift right=5pt, ""{description, name=b1b, pos=0, inner sep=0pt, outer sep=0pt}, ""{description, name=b1c, pos=1, inner sep=0pt, outer sep=0pt}]
   & & G_1C \arrow[r] & 0 \\
 & & & G_2C \arrow[ruu, equal] \arrow[r, "0"]
     \arrow[ruu, blue, rvpath, shift right=5pt, ""{description, name=b2b, pos=0, inner sep=0pt, outer sep=0pt}, ""{description, name=b2c, pos=1, inner sep=0pt, outer sep=0pt}]
   & G_3C \arrow[ru, "\dff^1_{13}"]
     \arrow[ru, blue, rvpath, -{Stealth[length=4pt]}, shift left=4pt, ""{description, name=b1d, pos=0, inner sep=0pt, outer sep=0pt}]
   & &
\end{tikzcd}
        \caption{The paths in red and blue correspond to the perspective pairs (i) and (ii) respectively. It indicates how the differentials are composed to give zero. 
        }\label{perspec1}
        \vspace{-1.5em}
    \end{figure}
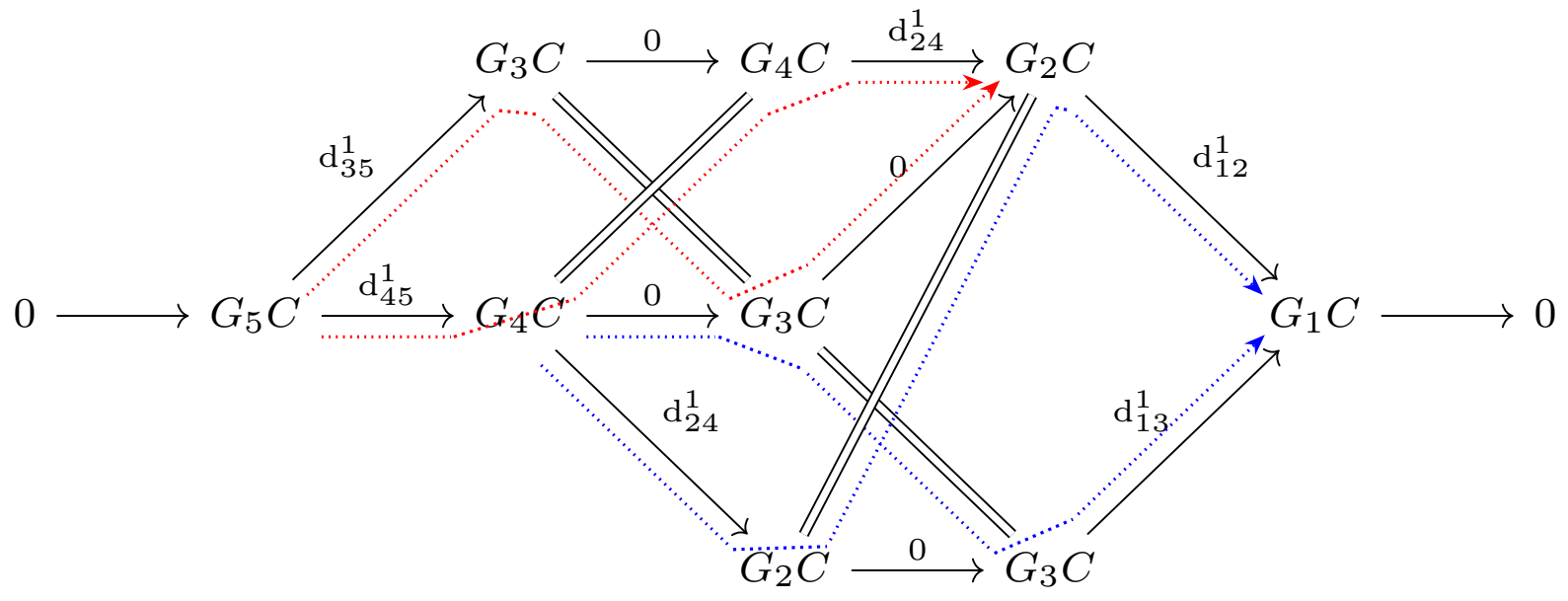

If we consider the fundamental spectral net then perspective pairs yield relations on the differentials. This is best illustrated with the example in Section \ref{theN5}.
In the lattice $\sO(\sN_5)$ we have three perspective pairs. We start with two immediate pairs: (i) $\beta,\delta,\epsilon$ and $\zeta$, with tail elements $\alpha$ and $\tp$, and (ii) $\alpha,\beta,\gamma$ and $\delta$, with tail elements $\bt$ and $\zeta$.
For (i) and (ii) we obtain the following relation on the differentials in Figure \ref{perspec1}:
\[
\begin{aligned}
&{\rm (i)\colon} \quad 0\cdot \dff^1_{35} + \dff^1_{24} \dff^1_{45} = \dff^1_{24} \dff^1_{45}=0;\\
&{\rm (ii)\colon} \quad \dff^1_{13} \cdot 0 + \dff^1_{12} \dff^1_{24} = \dff^1_{12} \dff^1_{24}=0.
\end{aligned}
\]
Here $\dff^1_{23}=0$ and $\dff^1_{34}=0$, since the connecting homomorphisms $k$ and $k'$ of a perspective pair vanish: $\dff^1_{23}=\dff_{\delta\beta\alpha}$ for (ii) and $\dff^1_{34}=\dff_{\zeta\delta\beta}$ for (i).
Figure \ref{perspec1} shows the paths of the differentials induced by the perspective pairs. For (i) in red and for (ii) in blue. Observe that in Figure \ref{fig:downset-lattice}[left] the associated orders do not necessarily form a diamond. The diamond in a perspective pair is found in the lattice of down-sets. 

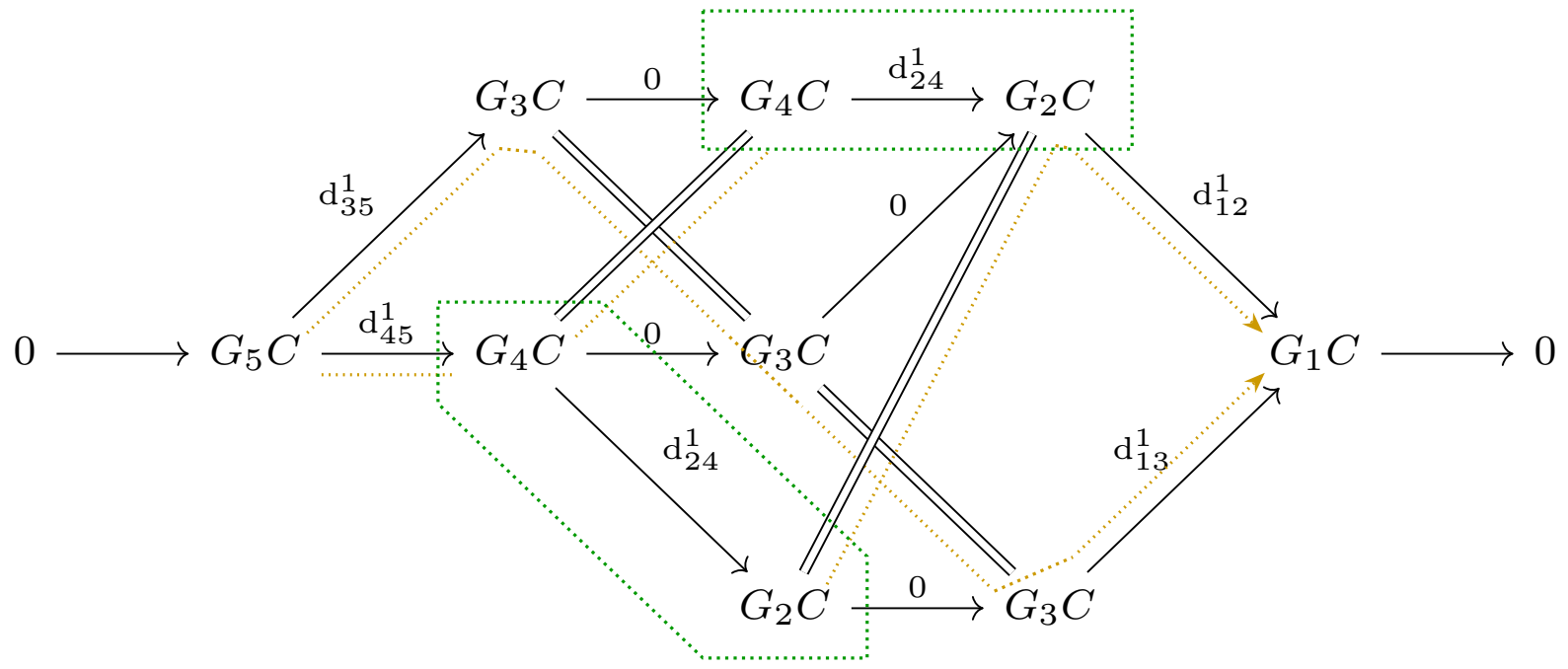
\begin{figure}[h]
\centering
\begin{tikzcd}[row sep=1.25cm, column sep=0.9cm, rvpath/.style={densely dotted, semithick, no head},
  execute at end picture={
    \begin{scope}[densely dotted, semithick]
    \draw[draw={rgb,255:red,204;green,153;blue,0}] (y1a.center) -- (y1b.center); \draw[draw={rgb,255:red,204;green,153;blue,0}] (y1c.center) -- (y1d.center); \draw[draw={rgb,255:red,204;green,153;blue,0}] (y1e.center) -- (y1f.center);
    \draw[draw={rgb,255:red,204;green,153;blue,0}] (y2e.center) -- (y2f.center);
    \draw[green!60!black, sharp corners] ($(G4t.north west)+(-3pt,10pt)$) rectangle ($(G2t.south east)+(3pt,-3pt)$);
    \draw[green!60!black, sharp corners] ($(G4m.north west)+(-3pt,3pt)$) -- ($(G4m.north east)+(3pt,3pt)$) -- ($(G2b.north east)+(3pt,3pt)$) -- ($(G2b.south east)+(3pt,-3pt)$) -- ($(G2b.south west)+(-3pt,-3pt)$) -- ($(G4m.south west)+(-3pt,-3pt)$) -- cycle;
    \end{scope}}]
 & & G_3C \arrow[rd, equal] \arrow[r, "0"]
     \arrow[rd, draw={rgb,255:red,204;green,153;blue,0}, rvpath, shift right=5pt, ""{description, name=y1b, pos=0, inner sep=0pt, outer sep=0pt}, ""{description, name=y1c, pos=1, inner sep=0pt, outer sep=0pt}]
   & |[alias=G4t]| G_4C \arrow[r, "\dff^1_{24}"]
   & |[alias=G2t]| G_2C \arrow[rd, "\dff^1_{12}"]
     \arrow[rd, draw={rgb,255:red,204;green,153;blue,0}, rvpath, -{Stealth[length=4pt]}, shift right=4pt, ""{description, name=y2f, pos=0, inner sep=0pt, outer sep=0pt}]
   & & \\
0 \arrow[r] & G_5C \arrow[ru, "\dff^1_{35}"] \arrow[r, "\dff^1_{45}"]
     \arrow[ru, draw={rgb,255:red,204;green,153;blue,0}, rvpath, shift right=4pt, ""{description, name=y1a, pos=1, inner sep=0pt, outer sep=0pt}]
     \arrow[r, draw={rgb,255:red,204;green,153;blue,0}, rvpath, shift right=4pt]
   & |[alias=G4m]| G_4C \arrow[ru, equal] \arrow[r, "0"] \arrow[rd, "\dff^1_{24}"]
     \arrow[ru, draw={rgb,255:red,204;green,153;blue,0}, rvpath, shift right=5pt]
   & G_3C \arrow[rd, equal] \arrow[ru, "0"]
     \arrow[rd, draw={rgb,255:red,204;green,153;blue,0}, rvpath, shift right=5pt, ""{description, name=y1d, pos=0, inner sep=0pt, outer sep=0pt}, ""{description, name=y1e, pos=1, inner sep=0pt, outer sep=0pt}]
   & & G_1C \arrow[r] & 0 \\
 & & & |[alias=G2b]| G_2C \arrow[ruu, equal] \arrow[r, "0"]
     \arrow[ruu, draw={rgb,255:red,204;green,153;blue,0}, rvpath, shift right=5pt, ""{description, name=y2e, pos=1, inner sep=0pt, outer sep=0pt}]
   & G_3C \arrow[ru, "\dff^1_{13}"]
     \arrow[ru, draw={rgb,255:red,204;green,153;blue,0}, rvpath, -{Stealth[length=4pt]}, shift left=4pt, ""{description, name=y1f, pos=0, inner sep=0pt, outer sep=0pt}]
   & &
\end{tikzcd}
        \caption{The paths in beige correspond to the perspective pair in (iii). 
        The green boxes indicate the Morse set corresponding to the convex set $\{2,4\}$. 
        }\label{perspec2}
        \vspace{-2.5em}
    \end{figure}

A third diamond can be found in the down-set lattice $\sO(\sN_5)$, v.\ Fig.\ \ref{fig:downset-lattice}[right], given by (iii) $\alpha,\gamma,\epsilon$ and $\zeta$, with tail elements $\bt$ and $\tp$.  In Figure \ref{perspec2} the path is indicated in beige. This pair involves, besides differentials on page~$1$ of the fundamental spectral net, differentials from and to the convex set $\{2,4\}$. These are differentials on page~$1$ of the spectral net over the sublattice $\{\bt,\alpha,\gamma,\epsilon,\zeta,\tp\}$, v.\ Defn.~\ref{subnetdefn}:
\[
{\rm (iii)\colon} \quad \dff^1_{13} \dff^1_{35} + \dff^1_{1,24} \dff^1_{24,5} = 0.
\]
The example shows that in the fundamental spectral net there are additional relations between differentials due to perspective pairs.
In Example \ref{franexample} we will apply these additional relations to an example discussed in \cite{fran}.

\begin{remark}
    \label{perspec3}
    The perspective pairs in (i) and (ii) in the above example are elementary because they give an immediate relation between differentials on page-1 of the fundamental spectral net. In \cite{fran} such perspective pairs are related to \emph{adjacent elements} in the order $\sP$.
    In particular the argument about additional connecting orbits between non-adjacent elements in \cite[Fig.\ 5.1]{fran}, cf.\ \cite{MischRein}, can be repeated using perspective pairs. We emphasize that the existence of a spectral matrix, or connection matrix, is not required. In \cite{CastVander2} we will study perspective pairs in general for first and higher order differentials.
\end{remark}

\subsection{Derived differentials}
\label{deriveddiff}
If we consider spectral systems of vector spaces, i.e.\ $\bR=\K$ is a field, then the short exact sequences in \eqref{shortexactSt} are split exact and therefore the differentials carry all information. However, for ring coefficients the short exact sequences in \eqref{shortexactSt} are not necessarily split exact. The extension problem can be solved via the information in the spectral system as discussed in Section \ref{convergence}. The differentials in a spectral net do not suffice in general. We define derived differentials in this setting which carry the information about non-split extensions. In Section \ref{connorbits} the latter is applied to the problem of connecting orbits in Morse representations.

Consider a subsystem of one exact triangle:
\begin{equation}
    \label{simpleextr}
\begin{tikzcd}
0 \arrow[r] & E^\gamma_\beta \arrow[rr, "k"] &                                 & E^\beta_\alpha \arrow[r] \arrow[ld, "i"] & 0 \\
            &                                & E^\gamma_\alpha \arrow[lu, "j"] &                                          &  
\end{tikzcd}
\end{equation}
On the second  (and final) page we obtain:
\[
\begin{tikzcd}[column sep=large]
0 \arrow[rr, "0", bend left] & S^{\gamma\gamma}_{\alpha\beta} \arrow[rr, "0"', bend right] & S^{\beta\gamma}_{\alpha\alpha} & 0
\end{tikzcd}
\]
The associated short exact sequence of $S$-terms is given by \eqref{shortexactSt}:
\begin{equation*}
\begin{tikzcd}
0 \arrow[r] & S_{\alpha\alpha}^{\beta\gamma} \arrow[r, "i"] & S^{\gamma\gamma}_{\alpha\alpha} \arrow[r, "j"] & S^{\gamma\gamma}_{\alpha\beta} \arrow[r] & 0.
\end{tikzcd}
\end{equation*}
This yields a distinguished triangle
in the derived category $\Dh(\sRmodR)$: 
\begin{equation*}
\begin{tikzcd}
 S_{\alpha\alpha}^{\beta\gamma} \arrow[r, "i"] & E^{\gamma}_{\alpha} \arrow[r, "j"] & S^{\gamma\gamma}_{\alpha\beta} \arrow[r,"\kappa"] & S_{\alpha\alpha}^{\beta\gamma}[1],
\end{tikzcd}
\end{equation*}
where $\kappa\colon S^{\gamma\gamma}_{\alpha\beta}\to S_{\alpha\alpha}^{\beta\gamma}[1]$ is the \emph{derived differential} of the triangle \eqref{simpleextr}. If $\beta\setminus\alpha=\{p\}$ and $\gamma\setminus\beta=\{q\}$ are singletons, then $k$ is the first order differential $\dff^1_{pq}$, and $\kappa$ is called a \emph{first order derived differential}.
The derived differential allows us to characterize when the exact triangle in \eqref{simpleextr} is split, i.e.\ $k=0$ and the short exact sequence $0\to E^\beta_\alpha\to E^\gamma_\alpha\to E^\gamma_\beta\to0$ splits. In particular, then $E^\gamma_\alpha\cong E^\beta_\alpha \oplus E^\gamma_\beta$, but the converse does not hold in general, v.\ Remark~\ref{splitnoniso}.
We consider the pair $(k,\kappa)\in \Hom(E^\gamma_\beta,E^\beta_\alpha)\times \Hom_{\Dh(\sRmodR)}\bigl(S^{\gamma\gamma}_{\alpha\beta},S^{\beta\gamma}_{\alpha\alpha}[1]\bigr)$. 
\begin{theorem}
\label{exactfortr}
Let $\iota(a)=(a,0)$ and $\pi(a,c)=c$. If $(k,\kappa)=0$, then there exists an isomorphism $\psi\colon E^\beta_\alpha \oplus E^\gamma_\beta\to E^\gamma_\alpha$ with $\psi\iota=i$ and $j\psi=\pi$. In particular, $E^\gamma_\alpha\cong E^\beta_\alpha \oplus E^\gamma_\beta$ in the   exact triangle in \eqref{simpleextr}.
Conversely, if $\psi\colon E^\beta_\alpha \oplus E^\gamma_\beta\xrightarrow[]{\cong} E^\gamma_\alpha$ satisfies $\psi\iota=i$ and $j\psi=\pi$, then, $(k,\kappa)=0$.
\end{theorem}
\begin{proof}
The condition $k=0$ implies that $S^{\beta\gamma}_{\alpha\alpha}\cong E^\beta_\alpha$ and $S^{\gamma\gamma}_{\alpha\beta}\cong E^\gamma_\beta$, which yields the short exact sequence
\[
\begin{tikzcd}
0 \arrow[r] & E_{\alpha}^{\beta} \arrow[r, "i"] & E^{\gamma}_{\alpha} \arrow[r, "j"] & E^{\gamma}_{\beta} \arrow[r] & 0.
\end{tikzcd}
\]
Since the derived differential $\kappa\colon E^\gamma_\beta\to E^\beta_\alpha[1]$ is trivial, this sequence splits, v.\ Sect.~\ref{convergence}, i.e.\ there exists $s\colon E^\gamma_\beta\to E^\gamma_\alpha$ with $js=\id$. Then $\psi(a,c)=i(a)+s(c)$ is an isomorphism with $\psi\iota=i$ and $j\psi=\pi$, and we conclude that $E^\gamma_\alpha\cong E^\beta_\alpha \oplus E^\gamma_\beta$.

For the converse we consider the exact sequences:
\[
\begin{tikzcd}
\cdots \arrow[r] & E_{\alpha}^{\beta} \arrow[r, "i"]                      & E^{\gamma}_{\alpha} \arrow[r, "j"]                                             & E^{\gamma}_{\beta} \arrow[r, "k"]                   & E^\beta_\alpha \arrow[r]                   & \cdots \\
\cdots \arrow[r] & E_{\alpha}^{\beta} \arrow[r, "\iota"] \arrow[u, "\id"] & E_{\alpha}^{\beta}\oplus E^{\gamma}_{\beta} \arrow[r, "\pi"] \arrow[u, "\psi","\cong"'] & E^{\gamma}_{\beta} \arrow[r, "k"] \arrow[u, "\id"'] & E^\beta_\alpha \arrow[r] \arrow[u, "\id"'] & \cdots
\end{tikzcd}
\]
From the bottom exact sequence we derive that $\ker k=E^\gamma_\beta$ and $\image~ k=\ker \iota=0$ and thus $k=0$. This implies that $\kappa=0$ since the short exact sequence
$0\to E^\beta_\alpha \to E_{\alpha}^{\beta}\oplus E^{\gamma}_{\beta} \to E^\gamma_\beta\to 0$ corresponds to the trivial class in $\Ext_\bR^1(E^\gamma_\beta,E^\beta_\alpha)$.
    \qed
\end{proof}

The negation of the above result can be formulated as follows:
\begin{corollary}
    \label{exactfortrneg}
    The conditions $k\neq 0$, or $\kappa\neq 0$, imply that
    there is no isomorphism $\psi\colon E^\beta_\alpha\oplus E^\gamma_\beta\to E^\gamma_\alpha$ with $\psi\iota=i$ and $j\psi=\pi$. 
\end{corollary}
\begin{proof}
The  statement is the contrapositive of the second part of Theorem~\ref{exactfortr}. 
    \qed
\end{proof}

\begin{remark}
\label{splitnoniso}
If $\bR$ is commutative and Noetherian and the $E$-terms are finitely generated, then $k\neq0$ or $\kappa\neq0$ implies $E^\gamma_\alpha\not\cong E^\beta_\alpha\oplus E^\gamma_\beta$. For $k=0$ this follows from \cite[Thm.~1]{Miyata}. Over such a ring, a short exact sequence of finitely generated modules whose middle term is isomorphic to the direct sum of its outer terms splits. For $k\neq0$ it follows by comparing lengths after localizing at a minimal prime of the support of $\image~k$. Without these assumptions the statement fails. For $\bR=\Z$, the sequence $0\to\Z\to\Z\oplus(\Z/2)^{(\N)}\to\Z/2\oplus(\Z/2)^{(\N)}\to0$, with $n\mapsto(2n,0)$ and $(n,x)\mapsto(\bar n,x)$, does not split and occurs in \eqref{simpleextr} with $k=0$ and $\kappa\neq0$, but its middle term is isomorphic to the direct sum of its outer terms.
\end{remark}

\begin{remark}
    Let $\bR$ be hereditary, for example a principal ideal domain, and let all $E$-terms $E^\beta_\alpha$ with $\beta\setminus\alpha$ a singleton be projective. Then $\kappa=0$ for every triangle \eqref{simpleextr} with $\beta\setminus\alpha$ and $\gamma\setminus\beta$ singletons. Indeed, $S^{\gamma\gamma}_{\alpha\beta}=\ker k$ is then a submodule of the projective module $E^\gamma_\beta$, hence projective, v.\ Rmk.~\ref{projkappa}. For a chain-generated $E$-system whose differential module is a chain complex $C$ vanishing in negative degrees, with subquotients $C^\beta_\alpha=F_\beta C/F_\alpha C$ consisting of projective modules, such as $\bfE\bigl(\sT(\sU)\bigr)$ in Sect.~\ref{invforsect}, the same holds over an arbitrary ring $\bR$. 
\end{remark}

\subsection{Higher differentials}
\label{highdiffspec}
The term higher differential is used to indicate differentials between $S$-terms, other than those between $E$-terms.  

\begin{definition}
\label{zigzag}
    Let $\sP$ be a poset. Two elements $p,q\in \sP$ are \emph{ connected inside} a convex set $I\in \sCo(\sP)$ if they belong to the same connected component of the Hasse diagram of $I$, regarded as an undirected graph. 
    Equivalently,   $p$ and $q$ are connected via a path $\{q_k\}\subset I$ of covering relations
    \begin{equation}
    \label{covrel123}
q=q_0 ~\overset{\prec}{_\succ}~ q_1 ~\overset{\prec}{_\succ}~ \cdots  ~\overset{\prec}{_\succ}~ q_k ~\overset{\prec}{_\succ}~\cdots ~\overset{\prec}{_\succ}~ q_{m-1} ~\overset{\prec}{_\succ}~ q_m=p,\quad m\ge 1,
\end{equation}
   where each step  $q_k ~\overset{\prec}{_\succ}~q_{k+1}$ indicates either $q_k\prec q_{k+1}$, or $q_k\succ q_{k+1}$.
\end{definition}

For example in Figure \ref{fig:downset-lattice} the elements $2$ and $3$ are connected in $\sP=\sN_5$ via the covering path $2\succ 1 \prec 3$, and $1$ and $5$ are connected in $\sN_5$ via the covering path $1\prec 3\prec  5$.
The following theorem relates disconnectedness of elements to a property of  differentials.
\begin{theorem}\label{propo0}
    Let $\bfS$ be an excisive spectral system over $\sO(\sP)$. 
    Let $p,q\in \sP$.
    Then, the differentials  
    \[
    \dff\colon S^{\epsilon\zeta}_{\gamma\delta}\rightarrow S^{\gamma\delta}_{\alpha\beta}, \quad \forall (\alpha,\beta,\gamma,\delta,\epsilon,\zeta) \in \bm{\sO(\sP)_6}, 
    \]
    with $\epsilon\setminus\delta = \lbrace p\rbrace$ and $\gamma\setminus\beta=\lbrace q\rbrace$, and such that $p$ and $q$ are not connected inside $\epsilon\setminus\beta$, are trivial, i.e.\ $\dff=0$.
\end{theorem}
\begin{proof}
Using the commutative squares for $S$-terms as given by Definition \ref{specsysdefn}(iv) we obtain:
\begin{equation}
    \label{anydiag}
\begin{tikzcd}[row sep=large, column sep=huge]
S^{\epsilon\epsilon}_{\gamma\delta} \arrow[r, "\dff^\dagger"] \arrow[d, "i"'] & S^{\gamma\delta}_{\beta\beta} \arrow[d, "\id"] \\
S^{\epsilon\zeta}_{\gamma\delta} \arrow[r, "\dff^\ddagger"]                & S^{\gamma\delta}_{\beta\beta}                 \\
S^{\epsilon\zeta}_{\gamma\delta} \arrow[u, "\id"] \arrow[r, "\dff"]        & S^{\gamma\delta}_{\alpha\beta} \arrow[u, "j"'] 
\end{tikzcd} \quad\implies \quad 
\begin{tikzcd}[row sep=large, column sep=huge]
S^{\epsilon\epsilon}_{\gamma\delta} \arrow[r, "\dff^\dagger"] \arrow[d, "i",two heads]                    & S^{\gamma\delta}_{\beta\beta}  \\
S^{\epsilon\zeta}_{\gamma\delta}  \arrow[r, "\dff"] & S^{\gamma\delta}_{\alpha\beta} \arrow[u, "j", tail]                                                                                        
\end{tikzcd}
\end{equation}
Therefore, $\dff^\dagger = j\dff i$ and, since $i$ is surjective and $j$ is injective, v.\ Defn.~\ref{specsysdefn}(iii), $\dff=0$ if and only if $\dff^\dagger=0$. It is enough to prove the theorem for $\dff^\dagger$.

Define $I$ as the path connected component of $q\in\epsilon\setminus\beta$, that is, the set of all points $q'$ that are connected to $q$ inside $\epsilon\setminus\beta$. Now, define $J=(\epsilon\setminus\beta)\setminus I$.
Define $\omega^- = \!\bigl\downarrow J$, and $\omega^+ = \!\bigl\downarrow I$. This gives the chain:
\[
\omega^-\cap \omega^+ \subset \beta\subsetneq \gamma \subset \delta\subsetneq \epsilon,\qquad \epsilon\setminus\beta \subset \omega^-\cup \omega^+.
\]
Here $\omega^-\cap \omega^+ \subset \beta$ holds  because $I$ and $J$ are unions of connected components of the convex set $\epsilon\setminus\beta$: an element of $\omega^-\cap\omega^+$ which is not in $\beta$ lies in $\epsilon\setminus\beta$ and below elements of $I$ and of $J$, and would therefore connect $I$ and $J$. Moreover, $\epsilon\setminus\beta=I\cup J\subset\omega^-\cup\omega^+$ by construction. In particular, $p\in J$, $q\in I$, and no element of $I$ is comparable to an element of $J$.

Define $\beta^\pm=\beta\cap \omega^\pm$, $\gamma^\pm = \gamma\cap \omega^\pm$, $\delta^\pm=\delta\cap \omega^\pm$, and $\epsilon^\pm = \epsilon\cap \omega^\pm$.
From the excisive property for $\bfE$ we conclude that 
\[
S^{\epsilon\epsilon}_{\gamma\delta}\cong S^{\epsilon^-\epsilon^-}_{\gamma^-\delta^-} \bigoplus
S^{\epsilon^+\epsilon^+}_{\gamma^+\delta^+},\quad
S^{\gamma\delta}_{\beta\beta} \cong S^{\gamma^-\delta^-}_{\beta^-\beta^-}\bigoplus  S^{\gamma^+\delta^+}_{\beta^+\beta^+},
\]
through the natural isomorphisms. Indeed, let $\chi'\subset\chi$ be elements of the chain $\beta\subset\gamma\subset\delta\subset\epsilon$. Then $\chi\setminus\chi'\subset I\cup J$, and the down-sets $\chi'\cup(\chi\cap I)$ and $\chi'\cup(\chi\cap J)$ form a perspective pair with $\chi'$ and $\chi$, so that $E^\chi_{\chi'}\cong E^{\chi'\cup(\chi\cap I)}_{\chi'}\oplus E^{\chi'\cup(\chi\cap J)}_{\chi'}$ by Sect.~\ref{firstordinfo}. By excision, $E^{\chi'\cup(\chi\cap I)}_{\chi'}\cong E^{\chi^+}_{\chi'^+}$ and $E^{\chi'\cup(\chi\cap J)}_{\chi'}\cong E^{\chi^-}_{\chi'^-}$. These isomorphisms are compatible with the homomorphisms $\ell$, which yields the isomorphisms for the $S$-terms.
The following squares are commutative
\[
\begin{tikzcd}
S^{\epsilon^-\epsilon^-}_{\gamma^-\delta^-} \arrow[rr, "\ell^-_\epsilon"] \arrow[dd, "\dff^-"] &  & S^{\epsilon\epsilon}_{\gamma\delta} \arrow[dd, "\dff^\dagger"] &  & S^{\epsilon^+\epsilon^+}_{\gamma^+\delta^+} \arrow[rr, "\ell^+_\epsilon"] \arrow[dd, "\dff^+"] &  & S^{\epsilon\epsilon}_{\gamma\delta} \arrow[dd, "\dff^\dagger"] \\
                                                                                                &  &                                                                             &  &                                                                                                 &  &                                                                             \\
S^{\gamma^-\delta^-}_{\beta^-\beta^-} \arrow[rr, "\ell^-_\beta"]                                &  & S^{\gamma\delta}_{\beta\beta}                                               &  & S^{\gamma^+\delta^+}_{\beta^+\beta^+} \arrow[rr, "\ell^+_\beta"]                                &  & S^{\gamma\delta}_{\beta\beta}                                              
\end{tikzcd}
\]
Under direct sum we obtain:
\[
\begin{tikzcd}
S^{\epsilon^-\epsilon^-}_{\gamma^-\delta^-}\oplus S^{\epsilon^+\epsilon^+}_{\gamma^+\delta^+} \arrow[rr, "\ell^-_\epsilon+ \ell^+_\epsilon"] \arrow[dd, "\dff^-\oplus \dff^+"'] &  & S^{\epsilon\epsilon}_{\gamma\delta} \arrow[dd, "\dff^\dagger"] \\
                                                                                                                                                                            &  &                                                                             \\
S^{\gamma^-\delta^-}_{\beta^-\beta^-}\oplus S^{\gamma^+\delta^+}_{\beta^+\beta^+} \arrow[rr, "\ell^-_\beta+\ell^+_\beta "]                                                        &  & S^{\gamma\delta}_{\beta\beta}                                              
\end{tikzcd}
\]
Due to linearity and commutativity we have that 
\[
\begin{aligned}
    \dff^\dagger \bigl(\ell^-_\epsilon+ \ell^+_\epsilon \bigr)(a,b) &= \dff^\dagger\ell^-_\epsilon(a)+ \dff^\dagger \ell^+_\epsilon(b)= \ell^-_\beta \dff^-a+ \ell^+_\beta \dff^+ b \\
    &= \bigl(\ell^-_\beta+\ell^+_\beta\bigr) (\dff^-\oplus \dff^+)(a,b)
\end{aligned}
\]
and since both $\ell^-_\epsilon+ \ell^+_\epsilon$ and $\ell^-_\beta+\ell^+_\beta $ are the natural isomorphisms, the differential $\dff^\dagger$ transforms to $\widetilde\dff^\dagger:=(\ell^-_\beta+\ell^+_\beta)^{-1}\,\dff^\dagger\,(\ell^-_\epsilon+\ell^+_\epsilon)$, given by:
\[
\widetilde\dff^\dagger = \begin{bmatrix}
    ~\dff^-~&~0~\\
    ~0~&~\dff^+~
\end{bmatrix}\colon
S^{\epsilon^-\epsilon^-}_{\gamma^-\delta^-} \bigoplus
S^{\epsilon^+\epsilon^+}_{\gamma^+\delta^+}\to
S^{\gamma^-\delta^-}_{\beta^-\beta^-}\bigoplus  S^{\gamma^+\delta^+}_{\beta^+\beta^+},
\]
v.\ \cite[Lem.\ 4.18]{Matschke}.
Since $\epsilon=\delta\cup\lbrace p\rbrace$ and $p\notin \omega^+$, it follows that $\epsilon^+=\delta^+$.
Similarly, since $\gamma=\beta \cup\lbrace q\rbrace$, with $q\notin \omega^-$, it follows that  $\gamma^-=\beta^-$. Consequently,
\[
S^{\epsilon^+\epsilon^+}_{\gamma^+\delta^+}=0,\quad\text{and}\quad 
S^{\gamma^-\delta^-}_{\beta^-\beta^-}=0.
\]
and thus $\widetilde\dff^\dagger=0$, i.e.\ $\dff^\dagger=0$, which proves the theorem.
\qed
\end{proof}

\begin{remark}\label{correctionremark}
    For any pair of parallel elements $p\Vert  q$ there  always exist down-sets $\epsilon$ and  $\beta$ such that 
    $p$ and $q$ are \emph{not} connected inside $\epsilon\setminus\beta$, by taking
  $\epsilon = \bigl\downarrow \lbrace q,p\rbrace$ and $\beta=\bigl\downarrow \lbrace q,p\rbrace\setminus\{q,p\}$. This implies $\epsilon\setminus\beta=\lbrace p,q\rbrace$, which shows that $p$ and $q$ are not connected in $\epsilon\setminus\beta$. 
  Hence, by Theorem~\ref{propo0}, $\dff=0$ for this choice of $\beta$ and $\epsilon$.
\end{remark}
\begin{corollary}\label{propo0eterms}
    For any pair of parallel elements $p\Vert  q$, the differential
    \[\dff\colon E^\epsilon_\delta\rightarrow E^\delta_\gamma \]
    with $\epsilon\setminus\delta=\lbrace p \rbrace$, and $\delta\setminus\gamma=\lbrace q \rbrace$, is trivial.
    \end{corollary}
    \begin{proof}
      Recall that $E^\epsilon_\delta=S^{\epsilon\epsilon}_{\delta\delta}$ and $E^\delta_\gamma=S^{\delta\delta}_{\gamma\gamma}$. By assumption, $\delta=\gamma\cup \lbrace q\rbrace$, and $\epsilon=\delta\cup \lbrace p\rbrace$. By construction $\epsilon\setminus\gamma = \{p,q\}$, and since $p\Vert q$, then $p$ and $q$
      are disconnected in $\epsilon\setminus\gamma$, and thus $\dff=0$ by Theorem \ref{propo0}.
    \end{proof}
\begin{remark}
\label{anyhigher}
    By \eqref{anydiag}, $\dff=0$ if and only if $\dff^\dagger=0$, and $\dff^\dagger$ does not depend on $\alpha$ and $\zeta$. In particular, the differential 
\[
    \dff^\ddagger\colon S^{\epsilon\smtp}_{\gamma\delta}\rightarrow S^{\gamma\delta}_{\smbt\beta},
\]
is trivial if and only if $\dff\colon S^{\epsilon\zeta}_{\gamma\delta}\rightarrow S^{\gamma\delta}_{\alpha\beta}$ is trivial. Equivalently, $\dff^\ddagger\neq 0$ if and only if $\dff\neq 0$.
\end{remark}

\begin{corollary}
    \label{propo1}
    Let $\bfS$ be an excisive spectral system over $\sO(\sP)$. 
    Suppose that a differential
    \[
    \dff\colon S^{\epsilon\zeta}_{\gamma\delta}\rightarrow S^{\gamma\delta}_{\alpha\beta}, \quad \text{for some } (\alpha,\beta,\gamma,\delta,\epsilon,\zeta) \in \bm{\sO(\sP)_6}, 
    \]
    with $\epsilon\setminus\delta = \lbrace p\rbrace$, $p\neq q$ and $\gamma\setminus\beta=\lbrace q\rbrace$, is non-trivial, i.e.\ $\dff\neq 0$.
Then, there is a path of covering relations in $\epsilon\setminus\beta$ such that $q\prec q_1 ~\overset{\prec}{_\succ}~ \cdots  ~\overset{\prec}{_\succ}~ q_k ~\overset{\prec}{_\succ}~\cdots ~\overset{\prec}{_\succ}~ q_{m-1}\prec p$, with $m\geq1$.  
\end{corollary}
\begin{proof}
The hypothesis that $\dff\neq 0$ implies that $p,q$ are connected, which is equivalent to the existence of a path $q ~\overset{\prec}{_\succ}~ q_1 ~\overset{\prec}{_\succ}~\dots ~\overset{\prec}{_\succ}~ q_{m-1} ~\overset{\prec}{_\succ}~ p$ inside $\epsilon\setminus\beta$. 
Note that $\gamma = \beta \cup \{q\}$ and $\epsilon=\delta\cup \{p\}$. Suppose there exists $p'\in \epsilon$ such that $p<p'$. Then, $p'\in \delta$ and thus $p\in \delta$, a contradiction. Consequently, $p$ is maximal in $\epsilon$ and thus maximal in $\epsilon\setminus\beta$.
Suppose there exists $q'\in \epsilon\setminus\beta$ such that $q'<q$. Then, $q'\in \gamma$ and thus $q'\in \beta$, a contradiction. Consequently, $q$ is minimal in $\epsilon\setminus\beta$.
Hence the first step of the path is $q\prec q_1$, and the last step is $q_{m-1}\prec p$. 
Hence, $q \prec q_1 ~\overset{\prec}{_\succ}~\dots ~\overset{\prec}{_\succ}~ q_{m-1} \prec p$. If $m=1$, this implies $q\prec p$.
    \qed
\end{proof}

A path in Corollary \ref{propo1} is called a \emph{directed path} from $p$ to $q$, and $p$ is said to be  \emph{directly path connected} to $q$ inside $\epsilon\setminus\beta$ (not a symmetric relation!).
 \begin{figure}[h] 
  \centering
\begin{tikzpicture}[scale=0.6, every node/.style={circle, draw, minimum size=2mm, inner sep=2pt, font=\small}]
  \node (1) at (0,0) {{\small $1$}};         
  \node (2) at (3,0) {{\small $2$}};
  \node (2') at (0,3) {{\small $2'$}};
  \node (3) at (3,3) {{\small $3$}};

  \draw (1) -- (2');
  \draw (2) -- (2');
  \draw (2) -- (3);
\end{tikzpicture}
\quad\quad\quad
\begin{tikzpicture}[
  scale=0.6,
  every node/.style={draw, rounded corners, fill=white, inner sep=2pt, font=\small}
]

  \node (D1) at (0,0) {$\bt=\emptyset$};

  \node (D2) at (-2,1.2) {$\alpha=\{1\}$};
  \node (D3) at (2,1.2) {$\beta=\{2\}$};

  \node (D4) at (0,2.4) {$\gamma=\{1,2\}$};
  \node (D5) at (4,2.4) {$\delta=\{2,3\}$};

  \node (D6) at (-2.0,3.6) {$\epsilon=\{1,2,2'\}$};
  \node (D7) at (2.0,3.6) {$\zeta=\{1,2,3\}$};
  
  \node (D8) at (0,4.8) {$\tp=\{1,2,2',3\}$};

  \draw (D1) -- (D2);
  \draw (D1) -- (D3);
  \draw (D2) -- (D4);
  \draw (D5) -- (D7);
  \draw (D3) -- (D5);
  \draw (D3) -- (D4);
  \draw (D4) -- (D6);
  \draw (D4) -- (D7);
  \draw (D6) -- (D8);
  \draw (D7) -- (D8);

\end{tikzpicture}
\caption{The $N$-poset $\sP$ [left] and its lattice of down-sets $\sO(\sP)$ [right]. }\label{Nexample}
\vspace{-1.5em}
\end{figure}
\begin{example}\label{examzigzag}
Consider $\sP$ and  $\sO(\sP)$ as given in Figure \ref{Nexample}. The elements $1$ and $3$ are connected or not, depending on the choices of $\epsilon$ and $\beta$.
First we take $\tp\setminus\bt=\sP$, in which case $1$ and $3$ are connected, where $1$ is minimal and $3$ is maximal in $\sP$. The directed path is given by
$1\prec 2'\succ 2\prec 3$.
If we   take $\tp\setminus\beta$, $\zeta\setminus\bt$ or $\zeta\setminus\beta$, then $1$ and $3$ cannot be connected.
\end{example}
  \begin{figure}[h] 
  \centering
\begin{tikzpicture}[scale=0.6, every node/.style={circle, draw, minimum size=2mm, inner sep=2pt, font=\small}]
  \fill[gray!20,semitransparent] (-0.5,-0.5) -- (-0.5,3.5) -- (3.5,3.5) -- (3.5,-0.5) -- cycle;
  \fill[red,nearly transparent] (-0.5,-0.5) -- (-0.5,-0.4) -- (3.5,-0.4)-- (3.5,-0.5) -- cycle; 
  \node (1) at (0,0) {{\small $1$}};         
  \node (2) at (3,0) {{\small $2$}};
  \node (2') at (0,3) {{\small $2'$}};
  \node (3) at (3,3) {{\small $3$}};

  \draw (1) -- (2');
  \draw (2) -- (2');
  \draw (2) -- (3);

\end{tikzpicture}
\quad\quad\quad
\begin{tikzpicture}[
  scale=0.6,
  every node/.style={draw, circle, fill=white, inner sep=2pt, font=\small}
]
 \fill[gray!20,semitransparent] (-0.5,-0.5) -- (-0.5,0.5) --(2.5,3.5) -- (3.5,3.5) -- (3.5,-0.5) -- cycle;
 \fill[red,nearly transparent] (-0.5,-0.5) -- (-0.5,-0.4) -- (3.5,-0.4)-- (3.5,-0.5) -- cycle; 
 
  \node (1) at (0,0) {{\small $1$}};         
  \node (2) at (3,0) {{\small $2$}};
  \node[draw = red] (2') at (0,3) {{\small \textcolor{red}{$2'$}}};
  \node (3) at (3,3) {{\small $3$}};

  \draw[red] (1) -- (2');
  \draw[red] (2) -- (2');
  \draw (2) -- (3);
 
\end{tikzpicture}
\quad\quad\quad
\begin{tikzpicture}[
  scale=0.6,
  every node/.style={draw, circle, fill=white, inner sep=2pt, font=\small}
]
\fill[gray!20,semitransparent] (-0.5,-0.5) -- (-0.5,3.5) -- (3.5,3.5) -- (3.5,-0.5) -- cycle;
\fill[red,nearly transparent] (0.5,-0.5) -- (3.5,2.5) -- (3.5,-0.5) -- cycle; 

  \node (1) at (0,0) {{\small $1$}};         
   \node[draw = red] (2) at (3,0) {{\small \textcolor{red}{$2$}}};
  \node (2') at (0,3) {{\small $2'$}};
  \node (3) at (3,3) {{\small $3$}};

  \draw (1) -- (2');
  \draw[red] (2) -- (2');
  \draw[red] (2) -- (3);

\end{tikzpicture}
\quad\quad\quad
\begin{tikzpicture}[
  scale=0.6,
  every node/.style={draw, circle, fill=white, inner sep=2pt, font=\small}
]
 \fill[gray!20,semitransparent] (-0.5,-0.5) -- (-0.5,0.5) --(2.5,3.5) -- (3.5,3.5) -- (3.5,-0.5) -- cycle;
  \fill[red,nearly transparent] (0.5,-0.5) -- (3.5,2.5) -- (3.5,-0.5) -- cycle;
  \node (1) at (0,0) {{\small $1$}};         
  \node[draw = red] (2) at (3,0) {{\small \textcolor{red}{$2$}}};
  \node[draw = red] (2') at (0,3) {{\small \textcolor{red}{$2'$}}};
  \node (3) at (3,3) {{\small $3$}};

  \draw[red] (1) -- (2');
  \draw[red] (2) -- (2');
  \draw[red] (2) -- (3);

\end{tikzpicture}
\caption{In the notation of Theorem \ref{propo0}, the gray areas represent $\epsilon$  and the red areas represent $\beta$. The difference $\epsilon\setminus\beta$ is also shaded gray. From left to right we consider the cases 
$\tp\setminus\bt$, $\zeta\setminus\bt$, $\tp\setminus\beta$ and $\zeta\setminus\beta$ used in Example \ref{examzigzag}. 
The red nodes and red edges do not belong to $\epsilon\setminus\beta$. 
Black nodes and black edges do belong to $\epsilon\setminus\beta$. Only in the leftmost figure can we connect nodes $1$ and $3$ through black lines.}\label{figexamzigzag}
\vspace{-1.5em}
\end{figure}

\begin{remark}
    All results described above are also valid for convex sets $\epsilon\setminus\delta$ and $\gamma\setminus\beta$ and not just single elements in $\sP$, where \say{$p$ and $q$ are not connected} is replaced by \say{no element of $\epsilon\setminus\delta$ is connected to an element of $\gamma\setminus\beta$ inside $\epsilon\setminus\beta$}.
\end{remark}

\section{Connecting orbits for Morse representations}
\label{connorbits}
For attractor-repeller pairs the connecting homomorphism in the long exact sequence of Conley indices  contains information about connecting orbits, i.e.\ if the connecting homomorphism is non-trivial, then the set of connecting orbits between attractor and repeller is non-void, v.\ \cite{conley:cbms,Vander1} and Lemma~\ref{arpair} below. This idea is translated to connection matrices by Franzosa for flow-orders, v.\ \cite[Prop.\ 5.3]{fran}.
In Franzosa's approach
non-zero entries of a connection matrix associated with adjacent Morse sets of the flow-order imply
connections between the adjacent Morse sets.
In \cite{Bartl1} Bart{\l}omiejczyk uses attractor-repeller pairs to relate first order differentials in spectral sequences to connecting orbits between adjacent Morse sets in any linear Morse representation. We show that the differentials in a spectral system also carry information about connecting orbits between Morse sets.
Throughout this section, $(\sM,\le)$ is a Morse representation, and the $S$-terms and differentials are those of the spectral system $\bfS=\bfS\bigl(\bfE(\sM)\bigr)$, v.\ Sect.~\ref{invforsect} and Thm.~\ref{duality}. Since $\bfE(\sM)$ is excisive, so is $\bfS$.

\begin{lemma}
\label{arpair}
Let $\alpha\subset\beta\subset\gamma$ in $\sO(\sM)$, and let $(A,A^*)$ be the attractor-repeller pair in $\Inv(G_{\gamma\setminus\alpha}U)$ given by $A=\Inv(G_{\beta\setminus\alpha}U)$ and $A^*=\Inv(G_{\gamma\setminus\beta}U)$. If $C(A^*,A)=\emptyset$, then there exists an isomorphism $\psi\colon E^\beta_\alpha\oplus E^\gamma_\beta\to E^\gamma_\alpha$ with $\psi\iota=i$ and $j\psi=\pi$, so that $(k,\kappa)=0$, v.\ Thm.~\ref{exactfortr}. In particular, if the connecting homomorphism $k\colon E^\gamma_\beta\to E^\beta_\alpha$ is non-trivial, then $C(A^*,A)\neq\emptyset$.
\end{lemma}
\begin{proof}
The chain $\bt\subset\alpha\subset\beta\subset\gamma\subset\tp$ is a sublattice of $\sO(\sM)$. It defines a coarser Morse representation $(\sM^\dagger,\le^\dagger)$, v.\ Sect.~\ref{tessMD}, whose Morse sets are the non-empty sets among $\Inv(F_\alpha U)$, $A$, $A^*$ and $\Inv(G_{\tp\setminus\gamma}U)$, linearly ordered in this order. Its lattice of down-sets is identified with the chain, so that $\beta\setminus\alpha=\{A\}$ and $\gamma\setminus\beta=\{A^*\}$, and by Remark~\ref{changeorder} and Birkhoff duality, $\bfE(\sM^\dagger,\le^\dagger)$ is the restriction of $\bfE(\sM)$ to the chain. Let $\le^\ddagger$ be the order obtained from $\le^\dagger$ by removing the relation $A<^\dagger A^*$, so that $A\Vert A^*$. The condition in Theorem~\ref{thm:mainMDthm} for $\le^\ddagger$ differs from the one for $\le^\dagger$ only for orbits from $A^*$ to $A$, and $C(A^*,A)=\emptyset$. Hence $(\sM^\dagger,\le^\ddagger)$ is again a Morse representation for $\phi$, and by Remark~\ref{changeorder} and Birkhoff duality we obtain an injection $\bfE(\sM^\dagger,\le^\dagger)\rightarrowtail\bfE(\sM^\dagger,\le^\ddagger)$. Besides the chain, $\sO(\sM^\dagger,\le^\ddagger)$ contains the down-set $\beta'$ with $\beta'\setminus\alpha=\{A^*\}$ and $\gamma\setminus\beta'=\{A\}$.

Since $\beta\cap\beta'=\alpha$ and $\beta\cup\beta'=\gamma$, the triples $(\alpha,\beta,\gamma)$ and $(\alpha,\beta',\gamma)$ form a perspective pair in $\sO(\sM^\dagger,\le^\ddagger)$; in the notation of Sect.~\ref{firstordinfo}, $\gamma$ and $\delta$ are replaced by $\beta'$ and $\gamma$. Since $\beta'\setminus\alpha=\gamma\setminus\beta$ and $\gamma\setminus\beta'=\beta\setminus\alpha$, excision gives isomorphisms $\ell\colon E^\beta_\alpha\to E^\gamma_{\beta'}$ and $\ell'\colon E^{\beta'}_\alpha\to E^\gamma_\beta$ with $j'i=\ell$ and $ji'=\ell'$, where $i'\colon E^{\beta'}_\alpha\to E^\gamma_\alpha$ and $j'\colon E^\gamma_\alpha\to E^\gamma_{\beta'}$. Hence $\ell^{-1}j'i=\id$, so that $i$ is injective and $k=0$, and $ji'\ell'^{-1}=\id$, so that $s=i'\ell'^{-1}$ is a section of $j$. Then $\psi(a,c)=i(a)+s(c)$ is the required isomorphism, v.\ the proof of Thm.~\ref{exactfortr}. The tail elements of Sect.~\ref{firstordinfo} are not needed. Via the injection $\bfE(\sM^\dagger,\le^\dagger)\rightarrowtail\bfE(\sM^\dagger,\le^\ddagger)$ the same holds in $\bfE(\sM^\dagger,\le^\dagger)$, that is, in $\bfE(\sM)$.
    \qed
\end{proof}

\subsection{Connecting orbits via higher order differentials}
\label{connorbhigherdiff}
Corollary \ref{propo1} yields an immediate result for a flow-order.

\begin{theorem}
    \label{conn2aa}
    Let $(\sM,\le)$ be a flow-order. Suppose the differential
    \[
    \dff\colon S^{\delta\smtp}_{\beta\gamma}\rightarrow S^{\beta\gamma}_{\smbt\alpha},
    \]
    is non-trivial, with $\beta\setminus\alpha=\{M\}$ and $\delta\setminus\gamma = \{M'\}$. Then, there exists a directed path $\{M_k\}{}_{k=0}^{n}\subset \delta\setminus\alpha$
    from $M'=M_n$ to $M=M_0$, such that for all $1\le k\le n$ either
    $C(M_k,M_{k-1})\neq \emptyset$, or $C(M_{k-1},M_{k})\neq \emptyset$, depending on whether $M_k\succ M_{k-1}$, or $M_k\prec M_{k-1}$ respectively.
\end{theorem}
\begin{proof}
Since $\dff\neq 0$, Corollary \ref{propo1} implies   that there exists a directed path $\{M_k\}$  from $M'$ to $M$.
Since $(\sM,\le)$ is a flow-order, each such covering relation yields a connecting orbit between $M_{k-1}$ and $M_k$, where the direction is determined by the order. Indeed, if $M_{k-1}\prec M_k$, then in the characterization of flow-orders after Defn.~\ref{flowdefn} the chain of connecting orbits from $M_k$ to $M_{k-1}$ has no intermediate elements, so that $C(M_k,M_{k-1})\neq\emptyset$.
    \qed
\end{proof}
Theorem \ref{conn2aa} may be interpreted as the spectral systems analogue of \cite[Prop.\ 5.3]{fran}.
For \emph{arbitrary} Morse representations the higher differentials in spectral systems still contain information about connecting orbits.

\begin{theorem}
    \label{conn2}
Let $(\sM,\le)$ be a Morse representation and suppose the differential
    \[
    \dff\colon S^{\delta\smtp}_{\beta\gamma}\rightarrow S^{\beta\gamma}_{\smbt\alpha},
    \]
    is non-trivial, with $\beta\setminus\alpha=\{M\}$ and $\delta\setminus\gamma = \{M'\}$. Then, there exists a sequence $\{M_k\}_{k=0}^{n} \\ \subset\delta\setminus\alpha$ from $M'=M_n$ to $M=M_0$, with $M_0<M_1$ and $M_{n-1}<M_n$, such that for all $1\le k\le n$ either $C(M_k,M_{k-1})\neq\emptyset$, or $C(M_{k-1},M_k)\neq\emptyset$, depending on whether $M_k>M_{k-1}$, or $M_k<M_{k-1}$ respectively.
\end{theorem}
\begin{proof}
Let $\le'$ be the flow-order obtained from $\le$ by removing order relations, v.\ the characterization after Defn.~\ref{flowdefn}. By Theorem~\ref{thm:mainMDthm}, $(\sM,\le')$ is again a Morse representation for $\phi$, the identity map is an order-surjection $(\sM,\le')\twoheadrightarrow(\sM,\le)$, and $\sO(\sM,\le)$ is a sublattice of $\sO(\sM,\le')$. By Remark~\ref{changeorder} and Birkhoff duality we obtain an injection $\bfE(\sM,\le)\rightarrowtail\bfE(\sM,\le')$, so that $\dff$ is also a non-trivial differential for $(\sM,\le')$. Theorem~\ref{conn2aa} for the flow-order $(\sM,\le')$ yields a directed path $\{M_k\}_{k=0}^{n}\subset\delta\setminus\alpha$ of $(\sM,\le')$ from $M'=M_n$ to $M=M_0$, with connecting orbits as stated. By Corollary~\ref{propo1}, the first and the last relation in this path are $M_0\prec'M_1$ and $M_{n-1}\prec'M_n$, since $M$ is minimal and $M'$ is maximal in $\delta\setminus\alpha$. Hence $M_0<M_1$ and $M_{n-1}<M_n$. Since $\le'\subset\le$, this proves the theorem.
    \qed
\end{proof}
If two Morse sets are adjacent in a Morse representation Theorem \ref{conn2} yields a connecting orbit whenever $\dff\colon S^{\gamma\smtp}_{\beta\beta}\to S^{\beta\beta}_{\smbt\alpha}$ is non-trivial, which is an extension of
 Bart{\l}omiejczyk's approach via spectral sequences. By Remark~\ref{anyhigher} and Theorem~\ref{conn2} we obtain:

\begin{corollary}
    \label{conn1}
    Let $\alpha \subsetneq\beta\subsetneq\gamma$ and $\beta\setminus\alpha=\{M\}$ and $\gamma\setminus\beta=\{M'\}$.
    Suppose the connecting homomorphism
    \[
    k\colon E_\beta^\gamma \to E_\alpha^\beta,
    \]
    is non-trivial. Then, the set of connecting orbits $C(M',M)$  from
    $M'$ to $M$ is non-empty.
\end{corollary}
\begin{proof}
The connecting homomorphism $k$ is the differential $\dff\colon S^{\gamma\gamma}_{\beta\beta}\to S^{\beta\beta}_{\alpha\alpha}$. By Remark~\ref{anyhigher}, $k\neq0$ if and only if $\dff\colon S^{\gamma\smtp}_{\beta\beta}\to S^{\beta\beta}_{\smbt\alpha}$ is non-trivial. Theorem~\ref{conn2}, with $(\alpha,\beta,\beta,\gamma)$ in place of $(\alpha,\beta,\gamma,\delta)$, applies, and since $\gamma\setminus\alpha=\{M,M'\}$, it yields $M<M'$ and $C(M',M)\neq\emptyset$. This also follows from Lemma~\ref{arpair}.
    \qed
\end{proof}

The next result reveals how   non-trivial higher differentials in a spectral system $\bfS$ contain information about  connecting orbits to and from specific Morse sets.

\begin{theorem}
\label{connviahigherdiff}
Let $(\sM,\le)$ be a Morse representation and
suppose the differential
    \[
    \dff\colon S^{\epsilon\zeta}_{\gamma\delta}\rightarrow S^{\gamma\delta}_{\alpha\beta}
    \]
    is non-trivial, with $\gamma\setminus\beta=\{M\}$ and $\epsilon\setminus\delta = \{M'\}$. Then,
    \begin{enumerate}[label=\textup{(\roman*)},leftmargin=2.2em,itemsep=0.1em,topsep=0pt]
        \item [{\rm (i)}] there exists a Morse set $M''\in \delta\setminus\beta$, $M''<M'$, such that $C(M',M'')\neq \emptyset$;
        \item[{\rm (ii)}] there exists a Morse set $M'''\in \epsilon\setminus\gamma$, $M<M'''$, such that $C(M''',M)\neq \emptyset$.
    \end{enumerate}
\end{theorem}

\begin{proof}
By Remark~\ref{anyhigher}, the differential $\dff\colon S^{\epsilon\smtp}_{\gamma\delta}\to S^{\gamma\delta}_{\smbt\beta}$ is non-trivial. Theorem~\ref{conn2}, with $(\beta,\gamma,\delta,\epsilon)$ in place of $(\alpha,\beta,\gamma,\delta)$, yields a sequence $\{M_k\}_{k=0}^{n}\subset\epsilon\setminus\beta$ from $M'=M_n$ to $M=M_0$, with $M_0<M_1$ and $M_{n-1}<M_n$, and with connecting orbits between consecutive elements. Hence (i) holds with $M''=M_{n-1}$: $C(M',M'')\neq\emptyset$, $M''<M'$, and $M''\in(\epsilon\setminus\beta)\setminus\{M'\}=\delta\setminus\beta$. Similarly, (ii) holds with $M'''=M_1$: $C(M''',M)\neq\emptyset$, $M<M'''$, and $M'''\in(\epsilon\setminus\beta)\setminus\{M\}=\epsilon\setminus\gamma$.
    \qed
\end{proof}

\begin{remark}\label{remarkhigher}
If $M\Vert M'$ in Theorem \ref{connviahigherdiff}, then the Morse sets $M''$ and $M'''$ in (i) and (ii) are different, since $M''=M'''$ would give $M<M'''=M''<M'$. Hence there exist
    two different Morse sets $M'''>M$ and $M''<M'$, and connections  $C(M''',M)\neq \emptyset$ and  $C(M',M'')\neq \emptyset$.
    Example \ref{Ssysnr2a} will illustrate that considering spectral sequences as opposed to a spectral net will not reveal most of the orbit structure.
\end{remark}

\begin{remark}\label{monotonehigher}
The following property of spectral systems gives a direct proof of Theorem~\ref{connviahigherdiff}. Let $\alpha_p$, $p\in\Z$, be a maximal chain in $\sO(\sM)$ as in Example~\ref{linearorderSterms}, with $\cE^k_p=S^{\alpha_p\alpha_{p-1+k}}_{\alpha_{p-k}\alpha_{p-1}}$ and $\dff^k_p\colon\cE^k_p\to\cE^k_{p-k}$, and let $r\ge 1$. For $1\le k\le r$ consider the differentials
\[
\dff^k\colon S^{\alpha_p\alpha_{p+k-1}}_{\alpha_{p-k}\alpha_{p-1}}\to S^{\alpha_{p-k}\alpha_{p-1}}_{\alpha_{p-r-k}\alpha_{p-r-1}}.
\]
For $k=r$ this is the differential $\dff^r_p\colon\cE^r_p\to\cE^r_{p-r}$ of the spectral sequence of the chain, and for $k=1$ it is the first order differential $\dff_{\alpha_p\alpha_{p-1}\alpha_{p-r-1}}\colon E^{\alpha_p}_{\alpha_{p-1}}\to E^{\alpha_{p-1}}_{\alpha_{p-r-1}}$. For $1<k<r$, $\dff^k$ is a differential in the spectral net which is not part of the spectral sequence. These differentials satisfy
\[
\dff^{k+1}\neq 0\quad\Longrightarrow\quad\dff^k\neq 0,\qquad 1\le k\le r-1.
\]
This follows from the commutative diagram
\[
\begin{tikzcd}[row sep=large, column sep=large, ampersand replacement=\&]
\cE^k_p \arrow[r, equals] \& S^{\alpha_{p}\alpha_{p+k-1}}_{\alpha_{p-k}\alpha_{p-1}} \arrow[r, "\dff^k"] \arrow[d, "j"', two heads] \& S^{\alpha_{p-k}\alpha_{p-1}}_{\alpha_{p-r-k}\alpha_{p-r-1}} \\
\& S^{\alpha_{p}\alpha_{p+k}}_{\alpha_{p-k}\alpha_{p-1}} \arrow[r, "\widetilde\dff"] \& S^{\alpha_{p-k}\alpha_{p-1}}_{\alpha_{p-r-k-1}\alpha_{p-r-1}} \arrow[u, "i"', tail] \\
\cE^{k+1}_p \arrow[r, equals] \& S^{\alpha_{p}\alpha_{p+k}}_{\alpha_{p-k-1}\alpha_{p-1}} \arrow[r, "\dff^{k+1}"] \arrow[u, "i", tail] \& S^{\alpha_{p-k-1}\alpha_{p-1}}_{\alpha_{p-r-k-1}\alpha_{p-r-1}} \arrow[u, "i"', tail]
\end{tikzcd}
\]
The maps $i$ are injective and $j$ is surjective, v.\ Defn.~\ref{specsysdefn}(iii) and (v), and the squares commute by Defn.~\ref{specsysdefn}(iv). The upper square gives $\dff^k=i\,\widetilde\dff\,j$, so that $\dff^k\neq 0$ if and only if $\widetilde\dff\neq 0$. This is also Remark~\ref{anyhigher}, since $\dff^k$ and $\widetilde\dff$ differ only in the labels $\zeta$ and $\alpha$. The lower square gives $\widetilde\dff\,i=i\,\dff^{k+1}$, that is, $\dff^{k+1}$ is the restriction of $\widetilde\dff$ to subterms. If $\dff^{k+1}a\neq 0$, then $\widetilde\dff(ia)=i\,\dff^{k+1}a\neq 0$, since $i$ is injective, so that $\widetilde\dff\neq 0$, and hence $\dff^k\neq 0$. Starting with $\dff^r=\dff^r_p$ and descending through $k=r-1,\dots,1$, a non-trivial differential $\dff^r_p$ forces a staircase of non-trivial differentials $\dff^{r-1},\dots,\dff^1$ in the spectral net. Along the staircase the sources run through the terms $\cE^k_p$ of the spectral sequence, from $\cE^r_p$ down to $\cE^1_p=E^{\alpha_p}_{\alpha_{p-1}}=\HC(M_p)$, while the targets grow from $\cE^r_{p-r}$ to $E^{\alpha_{p-1}}_{\alpha_{p-r-1}}$, the Conley index of the Morse set of $\{M_{p-r},\dots,M_{p-1}\}$. Similarly, the differentials
\[
\dff^k_*\colon S^{\alpha_p\alpha_{p+k-1}}_{\alpha_{p-r}\alpha_{p-r+k-1}}\to S^{\alpha_{p-r}\alpha_{p-r+k-1}}_{\alpha_{p-r-k}\alpha_{p-r-1}},\qquad 1\le k\le r,
\]
connect $\dff^r_*=\dff^r_p$ with the first order differential $\dff^1_*=\dff_{\alpha_p\alpha_{p-r}\alpha_{p-r-1}}\colon E^{\alpha_p}_{\alpha_{p-r}}\to E^{\alpha_{p-r}}_{\alpha_{p-r-1}}$, and $\dff^{k+1}_*\neq 0$ implies $\dff^k_*\neq 0$; here $\dff^{k+1}_*$ is induced by the corresponding differential $\widetilde\dff$ on quotients. By Lemma~\ref{arpair}, the two first order differentials yield connecting orbits from $M_p$ to the Morse set of $\{M_{p-r},\dots,M_{p-1}\}$ and from the Morse set of $\{M_{p-r+1},\dots,M_p\}$ to $M_{p-r}$. Together with Theorem~\ref{thm:mainMDthm} this gives Theorem~\ref{connviahigherdiff}. 
\end{remark}

\subsection{Connecting orbits via derived differentials}
\label{connorbderdiff}
Corollary \ref{conn1} gives a statement about first order differentials and the existence of connecting orbits. It may happen that $k=0$ even though connecting orbits may be forced by the spectral system, v.\ Ex.\ \ref{franexample2}.
In this case the first order derived differential may detect such connecting orbits.

Let $\alpha \subsetneq\beta\subsetneq\gamma$ and $\beta\setminus\alpha=\{M\}$ and $\gamma\setminus\beta=\{M'\}$. In this case the Morse set associated to $\gamma\setminus\alpha$ is denoted by $M_{\gamma\setminus\alpha}=\Inv(G_{\gamma\setminus\alpha}U)$ and
\[
M_{\gamma\setminus\alpha} = M\cup M'\cup C(M',M).
\]
The Conley indices are given by $\HC(M)=E^\beta_\alpha$, $\HC(M') = E^\gamma_\beta$ and
$\HC(M_{\gamma\setminus\alpha}) = E^\gamma_\alpha$.
\begin{lemma}
    \label{critconniff}
    $C(M',M)=\emptyset$ implies that there exists an isomorphism $\psi\colon E^\beta_\alpha\oplus E^\gamma_\beta\to E^\gamma_\alpha$ with $\psi\iota=i$ and $j\psi=\pi$, v.\ Thm.~\ref{exactfortr}. In particular, $E^\gamma_\alpha \cong E^\beta_\alpha\oplus E^\gamma_\beta$.
\end{lemma}
\begin{proof}
This is Lemma~\ref{arpair} with $\beta\setminus\alpha=\{M\}$ and $\gamma\setminus\beta=\{M'\}$, for which, in the notation of Lemma~\ref{arpair}, $A=M$ and $A^*=M'$.
    \qed
\end{proof}
\begin{remark}
    To prove Lemma \ref{critconniff} we can also use the property that the Conley index of a disjoint union is isomorphic to the direct sum of the Conley indices. In general
$E^\gamma_\alpha \cong E^\beta_\alpha\oplus E^\gamma_\beta$ does not imply that $C(M',M)=\emptyset$ as is shown in Example \ref{franexample2}.

\end{remark}
Using derived differentials in combination with Lemma \ref{critconniff} allows us to prove the following extension of Corollary \ref{conn1}:

\begin{theorem}
    \label{conn1forder}
    Let $\alpha \subsetneq\beta\subsetneq\gamma$ and $\beta\setminus\alpha=\{M\}$ and $\gamma\setminus\beta=\{M'\}$.
    Consider the connecting homomorphism (or first order differential) and the first order derived differential
    \[
    k\colon E_\beta^\gamma \to E_\alpha^\beta,\quad \kappa\colon S^{\gamma\gamma}_{\alpha\beta} \to S^{\beta\gamma}_{\alpha\alpha}[1].
    \]
    If $(k,\kappa) \neq 0$, then the set of connecting orbits $C(M',M)$  from
    $M'$ to $M$ is non-empty.
\end{theorem}
\begin{proof}
By Corollary \ref{exactfortrneg} the condition $(k,\kappa)\neq 0$ implies that there is no isomorphism $\psi\colon E^\beta_\alpha\oplus E^\gamma_\beta\to E^\gamma_\alpha$ with $\psi\iota=i$ and $j\psi=\pi$. By the contrapositive of Lemma \ref{critconniff} this implies that $C(M',M)\neq \emptyset$.
\qed
\end{proof}

In the case that $k=0$ the derived differential is given by $\kappa\colon E^\gamma_\beta \to E^\beta_\alpha[1]$ in the derived category, which justifies the terminology derived differential and reveals the higher order information hidden in the exact triangles. The occurrence of derived differentials for differentials between $S$-terms in general is a subject of future study, cf.\ \cite{CastVander2}.

\section{Examples}
\label{examplesofconn}
In this section we give some examples on how spectral nets are utilized to compute homology and find connecting orbits. The apparatus of spectral nets serves both purposes.

\begin{example}
\label{Ssysnr1}
Consider a flow $\phi$ on the 2-disk $D^2$ as depicted in Figure \ref{theflowonD2}.
The flow exhibits  a stable periodic orbit $M_1$,  a degenerate periodic orbit $M_2$ and  an unstable spiral $M_3$.
Figure \ref{filtering123} displays a filtering of trapping regions, which is a $4$-element lattice isomorphic to $\sO(\sP)$, where $\sP$ is the linear poset $1\prec 2\prec 3$.
\begin{figure}[ht]
    \begin{multicols}{2}
        \centering
        \begin{tikzpicture}[x=0.83cm,y=0.83cm,
          tip/.style={-{Latex[length=2mm,width=1.5mm]}},
          rtip/.style={red,-{Latex[length=1.6mm,width=1.2mm]}}]
          \definecolor{sviiiblue}{RGB}{173,216,230}
          \definecolor{sviiigold}{RGB}{255,215,0}
          \definecolor{sviiipink}{RGB}{255,192,203}
          \definecolor{sviiigreen}{RGB}{144,238,144}
          \definecolor{sviiinavy}{RGB}{0,0,145}
          \fill[sviiiblue] (0,0) circle (3);
          \draw[very thin,fill=sviiigold] (0,0) circle (2.23);
          \draw[very thin,fill=sviiipink] (0,0) circle (0.73);
          \draw[thick] (0,0) circle (3);
          \draw[thick] (0,0) circle (1.48);
          \draw[tip] (0,3) -- (-0.02,3);
          \draw[tip] (0,-3) -- (0.02,-3);
          \draw[tip] (0,1.48) -- (-0.02,1.48);
          \draw[tip] (0,-1.48) -- (0.02,-1.48);
          \fill (0,0) circle (0.07);
          \node at (0.1,2.7) {$M_1$};
          \node at (0.1,1.2) {$M_2$};
          \node at (0.1,-0.43) {$M_3$};
          \draw[rtip] (1.49,-0.53) .. controls (1.95,0.3) and (1.95,1.5) .. (1.52,2.23);
          \draw[rtip] (-1.49,0.58) .. controls (-1.95,-0.25) and (-1.95,-1.5) .. (-1.53,-2.24);
          \draw[rtip] (0.12,-0.08) .. controls (0.6,-0.06) and (0.86,0.3) .. (0.88,0.84);
          \draw[rtip] (-0.12,0.08) .. controls (-0.6,0.06) and (-0.86,-0.3) .. (-0.8,-0.93);
        \end{tikzpicture}
       \columnbreak
\[
\begin{tikzpicture}[scale=0.8, every node/.style={circle, draw, minimum size=2mm, inner sep=1pt, font=\small}]
  \node (1) at (0,0) {{\small $M_1$}};         
  \node (3) at (0,2) {{\small $M_2$}};
  \node (5) at (0,4) {{\small $M_3$}};       

  \draw (1) -- (3);
  \draw (3) -- (5);
\end{tikzpicture}
\quad\quad\quad\quad
\begin{tikzpicture}[
  scale=0.8,
  every node/.style={draw, rounded corners, fill=white, inner sep=2pt, font=\small}
]

  \node (D0) at (0,0) {$\bt=\emptyset$};

  \node (D1) at (0,1.6) {$\alpha=\{M_1\}$};

  \node (D2) at (0,3.2) {$\beta=\{M_1,M_2\}$};

  \node (D3) at (0,4.8) {$\tp=\{M_1,M_2,M_3\}$};

  \draw (D0) -- (D1);
  \draw (D1) -- (D2);
  \draw (D2) -- (D3);

\end{tikzpicture}
\]
         \end{multicols}
        \caption{The flow $\phi$ with Morse tessellation[left], the Morse representation $\sM$[middle], and the down-sets $\sO(\sM)$[right].}
     \label{theflowonD2}
    \end{figure}
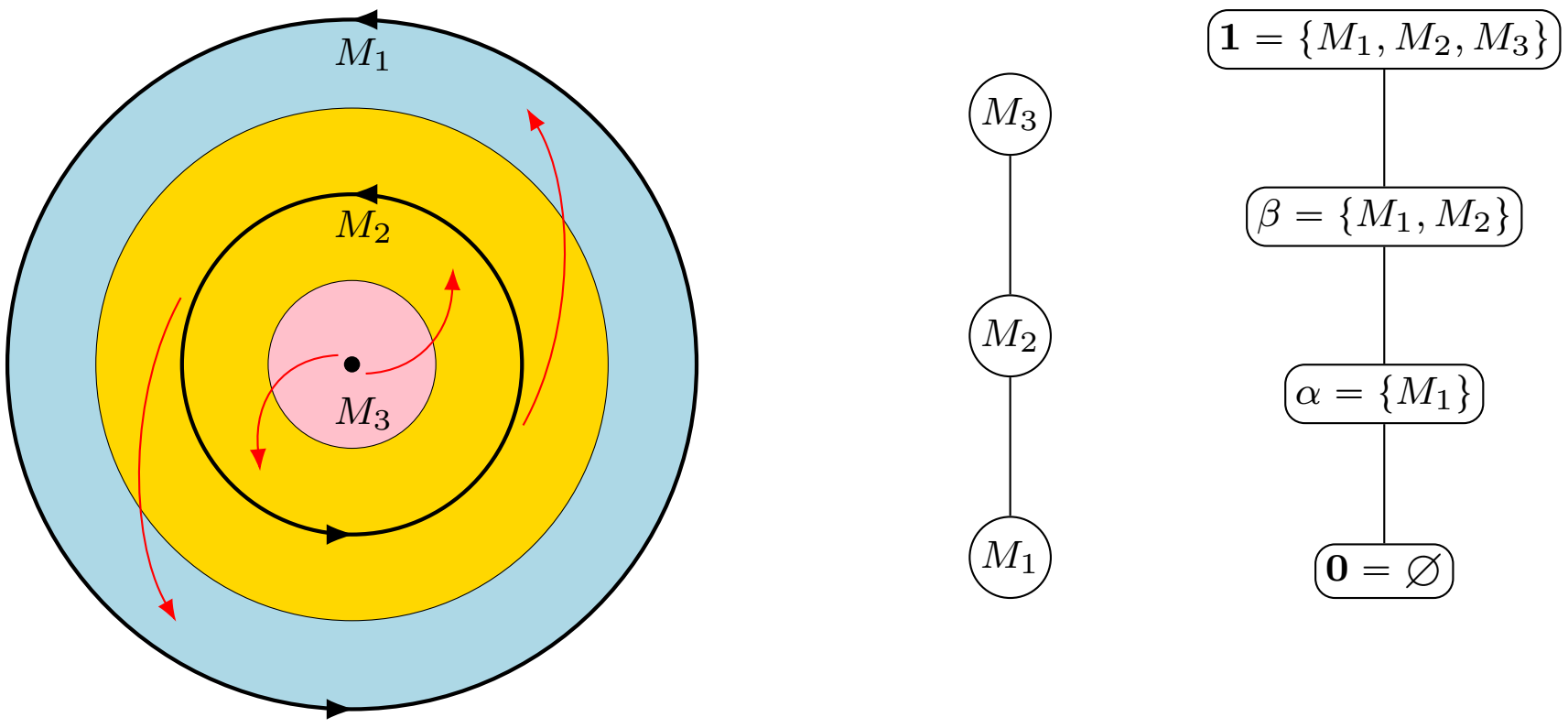
Specifically,
$F_\bt D^2=\emptyset$,
$F_\alpha D^2$ and  $F_\beta D^2$ are annuli and $F_\tp D^2=D^2$ is the whole disk.
The associated Morse tiles $G_pD^2$, $p=1,2,3$, give isolating neighborhoods for the Morse sets $M_1$, $M_2$ and $M_3$ respectively.
The filtering may be regarded as a  map $\sO(\sM)\to \sSet(D^2)$, with $\sO(\sM)
=\{\bt,\alpha,\beta,\tp\}$, $\alpha=\{M_1\}$, $\beta = \{M_1,M_2\}$.
This implies that the flow allows a Morse representation $\sM$ given by $M_1\prec M_2 \prec M_3$. We emphasize that at this point we do not know the latter is a flow-order.
In this example we can describe the Cartan-Eilenberg system $\bfE(\sM)$ with ring coefficients $R=\Z$. For the filtering we have:
\[
E^\alpha_\smbt = H(F_\alpha D^2)\cong \Z\oplus\Z,\quad E^\beta_\smbt = H(F_\beta D^2)\cong \Z\oplus \Z,\quad E^\smtp_\smbt = H(D^2)\cong \Z.
\]
Other $E$-terms that can be evaluated are
\[
\begin{aligned}
    \HC(M_1) &= E^\alpha_\smbt \cong \Z\oplus\Z;\\
    \HC(M_2) &= E^\beta_\alpha = H\bigl(F_\beta D^2,F_\alpha D^2\bigr)\cong 0;\\
    \HC(M_3) &= E^\smtp_\beta = H(F_\smtp D^2,F_\beta D^2)\cong\Z.
\end{aligned}
\]
If we invoke the spectral net approach most of the information about connecting orbits can be obtained from this data without further knowledge of the Cartan--Eilenberg system.
\begin{figure}[ht]
    \centering
    \begin{tikzpicture}[x=1cm,y=1cm,
      tip/.style={-{Latex[length=2mm,width=1.5mm]}}]
      \definecolor{sviiiblue}{RGB}{173,216,230}
      \definecolor{sviiigold}{RGB}{255,215,0}
      \definecolor{sviiipink}{RGB}{255,192,203}
      \definecolor{sviiigreen}{RGB}{144,238,144}
      \definecolor{sviiinavy}{RGB}{0,0,145}
      \node at (0.4,0) {\Large$\emptyset$};
      \begin{scope}[shift={(2.67,0)}]
        \fill[sviiiblue,even odd rule] (0,0) circle (1.43) (0,0) circle (1.08);
        \draw[very thin] (0,0) circle (1.08);
        \draw[thick] (0,0) circle (1.43);
        \draw[tip] (0,1.43) -- (-0.02,1.43);
        \draw[tip] (0,-1.43) -- (0.02,-1.43);
      \end{scope}
      \begin{scope}[shift={(6.27,0)}]
        \fill[sviiigold,even odd rule] (0,0) circle (1.43) (0,0) circle (0.36);
        \draw[very thin] (0,0) circle (0.36);
        \draw[thick] (0,0) circle (1.43);
        \draw[thick] (0,0) circle (0.72);
        \draw[tip] (0,1.43) -- (-0.02,1.43);
        \draw[tip] (0,-1.43) -- (0.02,-1.43);
        \draw[tip] (0,0.72) -- (-0.02,0.72);
        \draw[tip] (0,-0.72) -- (0.02,-0.72);
      \end{scope}
      \begin{scope}[shift={(9.83,0)}]
        \fill[sviiipink] (0,0) circle (1.43);
        \draw[thick] (0,0) circle (1.43);
        \draw[thick] (0,0) circle (0.72);
        \fill (0,0) circle (0.06);
        \draw[tip] (0,1.43) -- (-0.02,1.43);
        \draw[tip] (0,-1.43) -- (0.02,-1.43);
        \draw[tip] (0,0.72) -- (-0.02,0.72);
        \draw[tip] (0,-0.72) -- (0.02,-0.72);
      \end{scope}
      \node at (0.45,-1.86) {$F_\bt D^2$};
      \node at (2.67,-1.86) {$F_\alpha D^2$};
      \node at (6.27,-1.86) {$F_\beta D^2$};
      \node at (9.83,-1.86) {$F_\tp D^2$};
    \end{tikzpicture}
    \caption{The filtering $\sO(\sM)\to \sSet(D^2)$, with the periodic orbits highlighted.}
\label{filtering123}
\end{figure}

Consider the   page-$1$   spectral sequence based on the Conley index data:
\[
\begin{tikzcd}[column sep=large]
0 \arrow[r] & \HC(M_3) \arrow[r, "\dff^1_{23}"] & \HC(M_2) \arrow[r, "\dff^1_{12}"] & \HC(M_1) \arrow[r] & 0
\end{tikzcd}
\]
If we incorporate the information about the Conley indices and differentials we obtain:
\[
\begin{tikzcd}[column sep=large]
0 \arrow[r] & \Z \arrow[r, "0"] & 0 \arrow[r, "0"] & \Z\oplus \Z \arrow[r] & 0
\end{tikzcd}
\]
Computing homology we obtain the    page-$2$, where the middle differential is still to be determined:
\[
\begin{tikzcd}[column sep=large]\label{diagramsecondpageexample}
\cdots \arrow[rr, bend right] & 0 \arrow[rr, "0", bend left] & \Z \arrow[rr, "\dff_{13}^2", bend right] & 0 \arrow[rr, "0", bend left] & \Z\oplus \Z \arrow[rr, bend right] & 0 & \cdots
\end{tikzcd}
\]
By a classical result of spectral sequence theory, v.\ Sect.~\ref{convergence}, the infinity page $\bigoplus_p \cE^\infty_p$ is isomorphic to ${\rm Gr} E^\smtp_\smbt= {\rm Gr} H(D^2)$, and ${\rm rank}({\rm Gr} E^\smtp_\smbt) ={\rm rank}(H(D^2))=1$.
If $\dff_{13}^2=0$, then ${\rm Gr}E^\smtp_\smbt\cong \Z\oplus\Z\oplus\Z$, which is a contradiction. We conclude that $\dff^2_{13}\neq 0$ and $\dff_{13}^2$ is an injective map.
Consequently,  the system converges at page-$3$ and $\dff^r=0$ for $r\geq 3$.
If we keep track
  of the grading of homology, then the disk  has non-trivial homology only in dimension $0$. Using the grading on homology, and ordering the summands by decreasing degree, $\HC(M_1)=\HC_1(M_1)\oplus\HC_0(M_1)\cong\Z\oplus\Z$ and $\HC(M_3)=\HC_2(M_3)\cong\Z$, we obtain, up to the choice of generators:
  \begin{equation}
    \label{formfor13}
\dff_{13}^2=\begin{bmatrix}1\\0\end{bmatrix}\colon \Z\rightarrow \Z\oplus\Z,
\end{equation}
 because the generators in dimensions $1$ and $2$ have to be annihilated.
Since in this example all $E$-terms are given we can compute all differentials explicitly.

With this example we want to illustrate that with spectral nets/sequences information about connecting orbits may be obtained without knowing the $E$-system ($S$-system, or spectral matrix).
Note that,
\[
\dff_{13}^2\colon S^{\smtp\smtp}_{\alpha\beta}\rightarrow S^{\alpha\beta}_{\smbt\smbt} \neq 0
\]
which is a second order differential, and we emphasize that we are not assuming that $(\sM,\le)$ is a flow-order. By
Theorem \ref{connviahigherdiff} there exists a Morse set $M'''\in \lbrace M_2,M_3\rbrace$ such that $C(M''',M_1)\neq \emptyset$ and a Morse set $M''\in \lbrace M_1,M_2\rbrace$ such that $C(M_3,M'')\neq \emptyset$.
The periodic orbit $M_2$ is an invariant circle which separates $M_3$ from $M_1$ in $D^2$, so that no orbit goes from $M_3$ to $M_1$, i.e.\ $C(M_3,M_1)=\emptyset$. Hence, the only possible Morse set is $M''=M'''=M_2$ and therefore $C(M_2,M_1)\neq \emptyset$ and $C(M_3,M_2)\neq \emptyset$.
The latter also proves that $(\sM,\le)$ is
 flow-order, v.\ Defn.\ \ref{flowdefn}.

   This example shows that minimal algebraic information provides insight into Morse representations and even flow-orders. In particular,  higher order differentials play a role in finding connecting orbits. Observe  that algebraically the Morse set $M_2$ is \say{trivial}. Under small perturbations the connections found by Theorem \ref{connviahigherdiff} may become a connection from $M_3$ to $M_1$. In the next example we adjust Example \ref{Ssysnr1} slightly such that  higher order differentials are used for finding connections between topologically non-trivial Morse sets.

\end{example}

\begin{example}
    \label{Ssysnr2a}
    Consider   a similar flow $\phi$ on the $2$-disk $D^2$ as shown in Figure \ref{theflowonD3}. We consider a Morse representation
    consisting of $\{M_1,M_2,M_{2'},M_3\}$ where
    $M_1$ and $M_2$ are stable periodic orbits, $M_{2'}$ is an unstable periodic orbit,  and $M_3$ is an unstable spiral.
    We consider a filtering of $D^2$ based on the down-sets of the Morse sets, v.\ Fig.\ \ref{theflowonD3}[right]. The filtering determines a Morse representation $(\sM,\le)$, v.\ Fig.\ \ref{theflowonD3}[middle].
    As before we do not assume that $(\sM,\le)$ is a flow-order.
    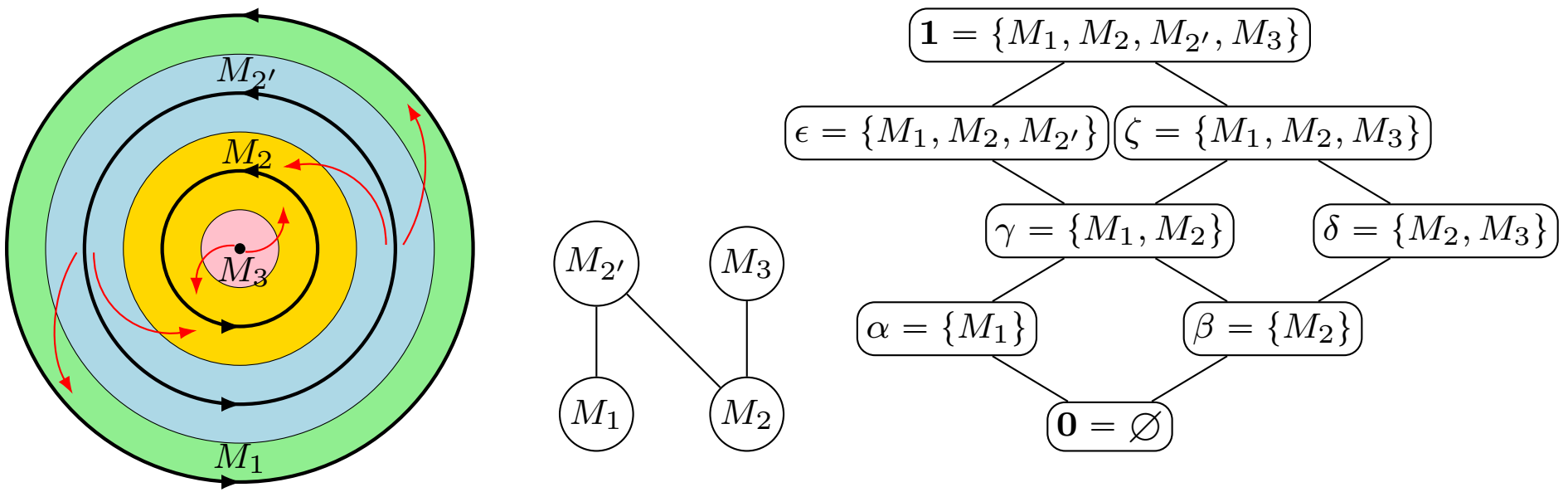
\begin{figure}[h]
    \begin{multicols}{3}
        \centering
        \vspace{5mm}
        \begin{tikzpicture}[x=0.62cm,y=0.62cm,font=\footnotesize,
          tip/.style={-{Latex[length=1.6mm,width=1.2mm]}},
          rtip/.style={red,-{Latex[length=1.4mm,width=1mm]}}]
          \definecolor{sviiiblue}{RGB}{173,216,230}
          \definecolor{sviiigold}{RGB}{255,215,0}
          \definecolor{sviiipink}{RGB}{255,192,203}
          \definecolor{sviiigreen}{RGB}{144,238,144}
          \definecolor{sviiinavy}{RGB}{0,0,145}
          \fill[sviiigreen] (0,0) circle (3);
          \draw[very thin,fill=sviiiblue] (0,0) circle (2.5);
          \draw[very thin,fill=sviiigold] (0,0) circle (1.5);
          \draw[very thin,fill=sviiipink] (0,0) circle (0.5);
          \draw[thick] (0,0) circle (3);
          \draw[thick] (0,0) circle (2);
          \draw[thick] (0,0) circle (1);
          \draw[tip] (0,3) -- (-0.02,3);
          \draw[tip] (0,-3) -- (0.02,-3);
          \draw[tip] (0,2) -- (-0.02,2);
          \draw[tip] (0,-2) -- (0.02,-2);
          \draw[tip] (0,1) -- (-0.02,1);
          \draw[tip] (0,-1) -- (0.02,-1);
          \fill (0,0) circle (0.07);
          \node at (0.1,2.26) {$M_{2'}$};
          \node at (0.1,1.22) {$M_2$};
          \node at (0.05,-0.3) {$M_3$};
          \node at (0,-2.7) {$M_1$};
          \draw[rtip] (2.1,0.05) .. controls (2.45,0.6) and (2.45,1.35) .. (2.15,1.87);
          \draw[rtip] (-2.1,-0.05) .. controls (-2.45,-0.6) and (-2.45,-1.35) .. (-2.15,-1.87);
          \draw[rtip] (1.88,0.05) .. controls (1.88,0.8) and (1.25,1.15) .. (0.55,1.05);
          \draw[rtip] (-1.88,-0.05) .. controls (-1.88,-0.8) and (-1.25,-1.15) .. (-0.55,-1.05);
          \draw[rtip] (0.07,-0.04) .. controls (0.4,-0.06) and (0.56,0.2) .. (0.56,0.55);
          \draw[rtip] (-0.07,0.04) .. controls (-0.4,0.06) and (-0.56,-0.2) .. (-0.56,-0.6);
        \end{tikzpicture}
       \columnbreak
\[
\begin{tikzpicture}[scale=0.4, every node/.style={circle, draw, minimum size=2mm, inner sep=1pt, font=\small}]
  \node (1) at (0,0) {{\small $M_1$}};         
  \node (2) at (3,0) {{\small $M_2$}};
  \node (2') at (0,3) {{\small $M_{2'}$}};
  \node (3) at (3,3) {{\small $M_3$}};

  \draw (1) -- (2');
  \draw (2) -- (2');
  \draw (2) -- (3);
\end{tikzpicture}
\vspace{10mm}
 \columnbreak
\begin{tikzpicture}[
  scale=0.65,
  every node/.style={draw, rounded corners, fill=white, inner sep=2pt, font=\small}
]

  \node (D1) at (0,0) {$\bt=\emptyset$};

  \node (D2) at (-2,1.2) {$\alpha=\{M_1\}$};
  \node (D3) at (2,1.2) {$\beta=\{M_2\}$};

  \node (D4) at (0,2.4) {$\gamma=\{M_1,M_2\}$};
  \node (D5) at (4,2.4) {$\delta=\{M_2,M_3\}$};

  \node (D6) at (-2.0,3.6) {$\epsilon=\{M_1,M_2,M_{2'}\}$};
  \node (D7) at (2.0,3.6) {$\zeta=\{M_1,M_2,M_3\}$};

  \node (D8) at (0,4.8) {$\tp=\{M_1,M_2,M_{2'},M_3\}$};

  \draw (D1) -- (D2);
  \draw (D1) -- (D3);
  \draw (D2) -- (D4);
  \draw (D5) -- (D7);
  \draw (D3) -- (D5);
  \draw (D3) -- (D4);
  \draw (D4) -- (D6);
  \draw (D4) -- (D7);
  \draw (D6) -- (D8);
  \draw (D7) -- (D8);

\end{tikzpicture}
\]
\end{multicols}
        \caption{The flow $\phi$ with Morse tessellation[left], the Morse representation $\sM$[middle], and the down-sets $\sO(\sM)$[right].}
     \label{theflowonD3}
    \end{figure}
For the Cartan--Eilenberg system $\bfE(\sM)$ (with ring coefficients $R=\Z$) we consider the $E$-terms associated with the lattice of down-sets $\sO(\sM)$:
\[
E^\alpha_\smbt = H(F_\alpha D^2)\cong \Z\oplus\Z,~~ E^\beta_\smbt = H(F_\beta D^2)\cong \Z\oplus \Z,~~  E^\gamma_\smbt = H(F_\gamma D^2)\cong \Z\oplus\Z\oplus\Z\oplus\Z,
\]
\[E^\delta_\smbt = H(F_\delta D^2)\cong \Z,~~
E^\epsilon_\smbt = H(F_\epsilon D^2)\cong \Z\oplus\Z,~~ E^\zeta_\smbt = H(F_\zeta D^2)\cong \Z\oplus\Z\oplus\Z.
\]
and the homology of the whole space $E^\smtp_\smbt = H(D^2)\cong \Z.$ The Conley indices of the Morse sets are given by
\[
    \HC(M_1) = E^\alpha_\smbt \cong E^\zeta_\delta \cong E^\gamma_\beta\cong\Z\oplus\Z, \quad  \HC(M_{2'}) = E^\epsilon_\gamma \cong E^\smtp_\zeta\cong \Z\oplus\Z;\]
    \[
    \HC(M_2) = E^\beta_\smbt \cong E^\gamma_\alpha\cong\Z\oplus\Z, \quad  \HC(M_3) =  E^\delta_\beta \cong E^\zeta_\gamma \cong E^\smtp_\epsilon\cong \Z.\\
\]
In total there are $23$ (strict) $E$-terms which are not all computed.
The lattice $\sO(\sM)$ has $5$ maximal chains and there are
 $3$ different perspective pairs consisting of covering relations, v.\ Sect.\ \ref{firstordinfo}.
 The objective of this example is to illustrate that higher order differentials may contain information about connecting orbits in case first order differentials do not.

 As   in Section~\ref{theN5}, page-$1$ of the spectral net is given by:
 
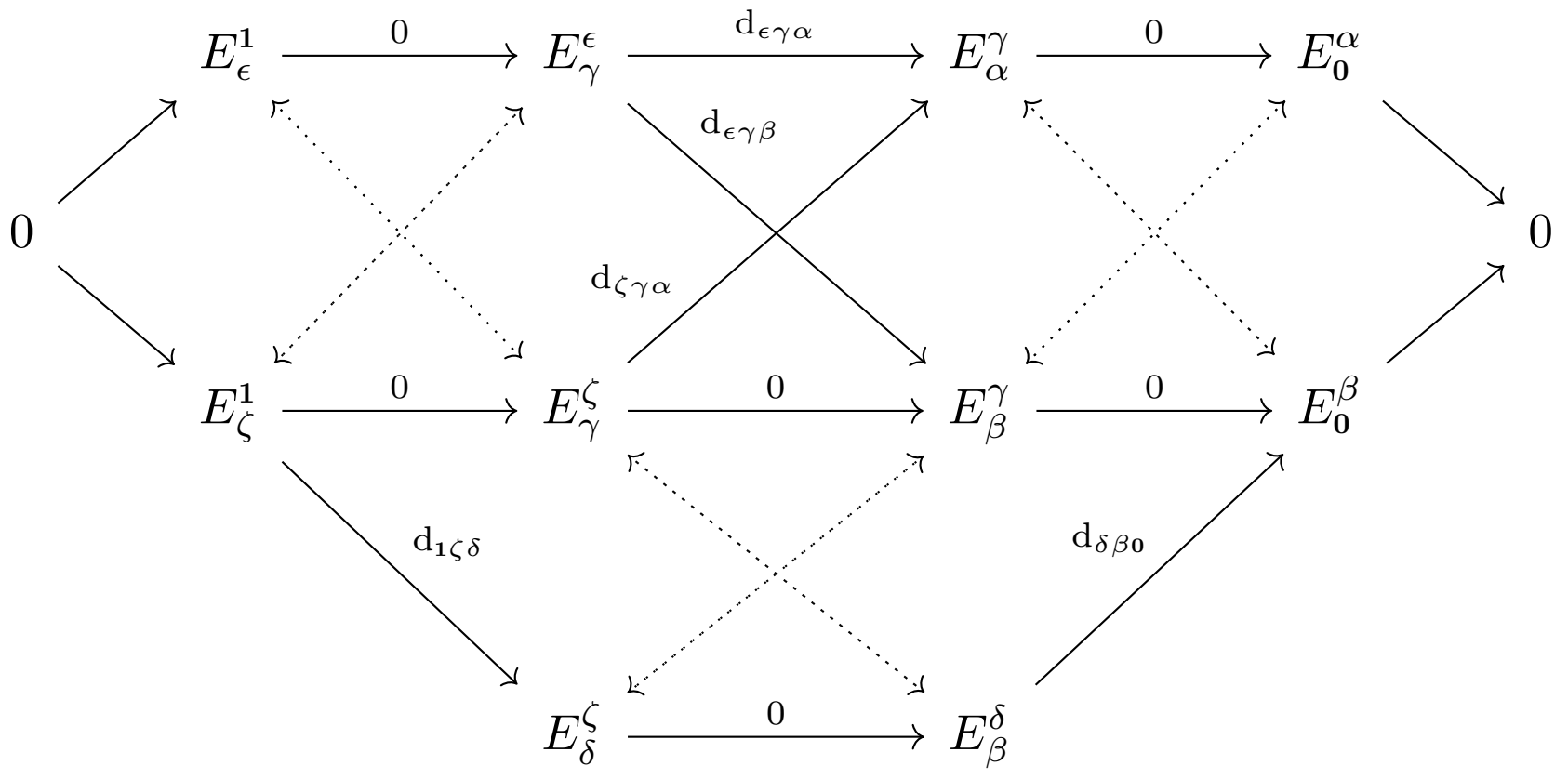
\[
 \begin{tikzcd}
 \label{diamant}
 & E^\smtp_\epsilon \arrow[rr, "0"] \arrow[rrdd, dotted]                                      &  & E^\epsilon_\gamma \arrow[rr, "\dff_{\epsilon\gamma\alpha}"] \arrow[rrdd, "\dff_{\epsilon\gamma\beta}"{pos=0.20}] \arrow[lldd, dotted] &              & E^\gamma_\alpha \arrow[rr, "0"] \arrow[rrdd, dotted]                        &  & E^\alpha_\smbt \arrow[rd] \arrow[lldd, dotted] &   \\
0 \arrow[ru] \arrow[rd] &                                                                                            &  &                                                                                                                             &              &                                                                             &  &                                                & 0 \\
                        & E^\smtp_\zeta \arrow[rr, "0"] \arrow[rrdd, "\dff_{\smtp\zeta\delta}"] \arrow[rruu, dotted] &  & E^\zeta_\gamma \arrow[rr, "0"] \arrow[rrdd, dotted] \arrow[lluu, dotted] \arrow[rruu, "\dff_{\zeta\gamma\alpha}"{pos=0.20}]          &              & E^\gamma_\beta \arrow[rr, "0"] \arrow[rruu, dotted] \arrow[lldd, dotted]    &  & E^\beta_\smbt \arrow[ru] \arrow[lluu, dotted]  &   \\
                        &                                                                                            &  &                                                                                                                             & \hspace{1mm} &                                                                             &  &                                                &   \\
                        &                                                                                            &  & E^\zeta_\delta \arrow[rr, "0"] \arrow[rruu, dotted]                                                                         &              & E^\delta_\beta \arrow[rruu, "\dff_{\delta\beta\smbt}"] \arrow[lluu, dotted] &  &                                                &
\end{tikzcd}
\]
The dotted arrows are the isomorphisms from excision and by
Corollary \ref{propo0eterms} all differentials in the squares with isomorphisms are trivial, v.\
Fig.\ \ref{theflowonD3}.
Moreover, the  differentials $\dff_{\smtp\zeta\delta}$ and $\dff_{\epsilon\gamma\beta}$, and the differentials
$\dff_{\zeta\gamma\alpha}$ and $\dff_{\delta\beta\smbt}$
are isomorphic. To illustrate the role of higher order differentials we denote the differential from $\HC(M_{2'})$ to $\HC(M_2)$ by $\dff^1_{22'}$. The latter can be computed via $\dff_{\epsilon\gamma\alpha}$.
To do so we unroll the exact triangles for $M_2$ and $M_{2'}$ as long exact sequences using the standard grading of singular homology:
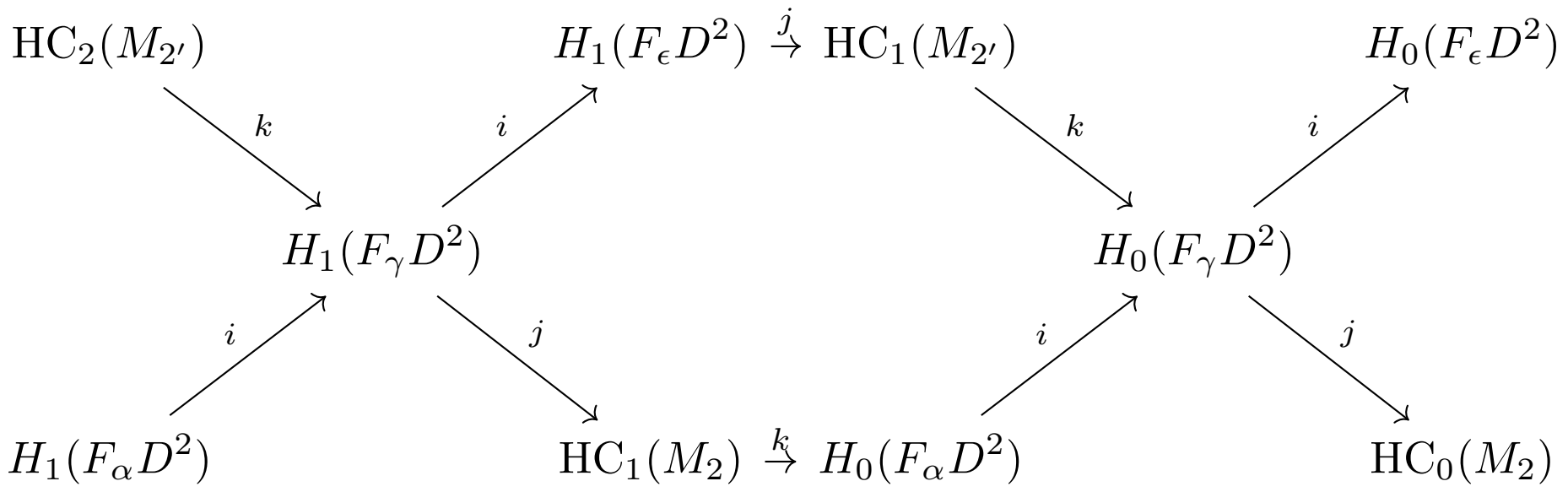
\[\begin{tikzcd}[row sep=large, column sep=tiny]
\HC_2(M_{2'}) \arrow[rd, "k"]      &                                                    &  H_1(F_\epsilon D^2) \arrow[r, "j"] & \HC_1(M_{2'}) \arrow[rd, "k"]      &                                                    &  H_0(F_\epsilon D^2) \\
                                   &  H_1(F_\gamma D^2) \arrow[rd, "j"] \arrow[ru, "i"] &                                     &                                    &  H_0(F_\gamma D^2) \arrow[rd, "j"] \arrow[ru, "i"] &                      \\
 H_1(F_{\alpha} D^2) \arrow[ru, "i"] &                                                    & \HC_1(M_{2}) \arrow[r, "k"]         &  H_0(F_{\alpha} D^2) \arrow[ru, "i"] &                                                    & \HC_0(M_{2})
\end{tikzcd}\]
and determining all homomorphisms gives:
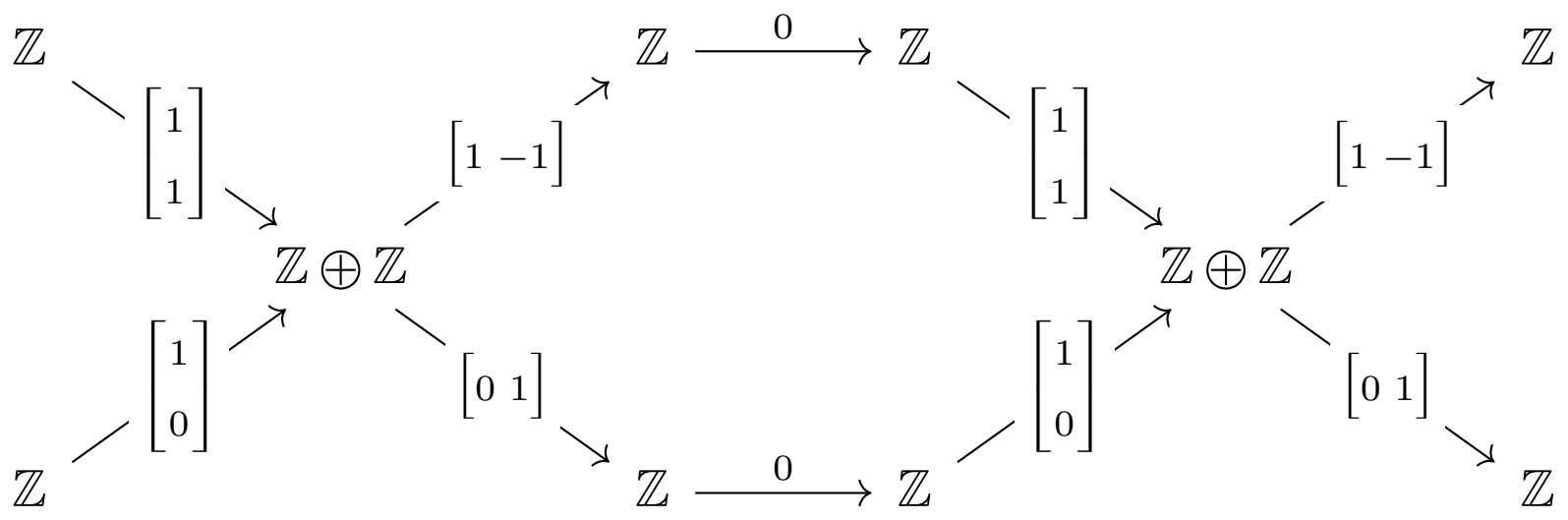
\[\begin{tikzcd}[row sep=large, column sep =1.2cm] 
\Z \arrow[rd, "\begin{bmatrix}1\\1\end{bmatrix}" description] &                                                                                                                                     & \Z \arrow[r, "0"  ] & \Z \arrow[rd, "\begin{bmatrix}1\\1\end{bmatrix}" description] &                                                                                                                                     & \Z \\
                                                               & \Z\oplus \Z \arrow[ru, "{\begin{bmatrix}1~-1\end{bmatrix}}" description] \arrow[rd, "{\begin{bmatrix}0~1\end{bmatrix}}" description] &                               &                                                                & \Z\oplus \Z \arrow[ru, "{\begin{bmatrix}1~-1\end{bmatrix}}" description] \arrow[rd, "{\begin{bmatrix}0~1\end{bmatrix}}" description] &    \\
\Z \arrow[ru, "\begin{bmatrix}1\\0\end{bmatrix}" description]  &                                                                                                                                     & \Z \arrow[r, "0" ]             & \Z \arrow[ru, "\begin{bmatrix}1\\0\end{bmatrix}" description]  &                                                                                                                                     & \Z
\end{tikzcd}\]
The differential from $\HC(M_{2'})$ to $\HC(M_2)$ is given by the composition $jk$. Taking the direct sum over the grading, with the summands ordered by decreasing degree, we obtain:
\[
\dff^1_{22'} \sim \dff_{\epsilon\gamma\alpha} =
\begin{bmatrix}0&1&0&0\\0&0&0&1\end{bmatrix}\cdot
\begin{bmatrix}1&0\\1&0\\0&1\\0&1\end{bmatrix}  =
\begin{bmatrix}1&0\\0&1\end{bmatrix},
\quad
\dff^1_{12'}=\begin{bmatrix}1&0\\0&1\end{bmatrix},\quad \dff^1_{23}=\begin{bmatrix}1\\0\end{bmatrix}.
\]
Page-$1$ of the (fundamental) spectral net can be displayed as follows:
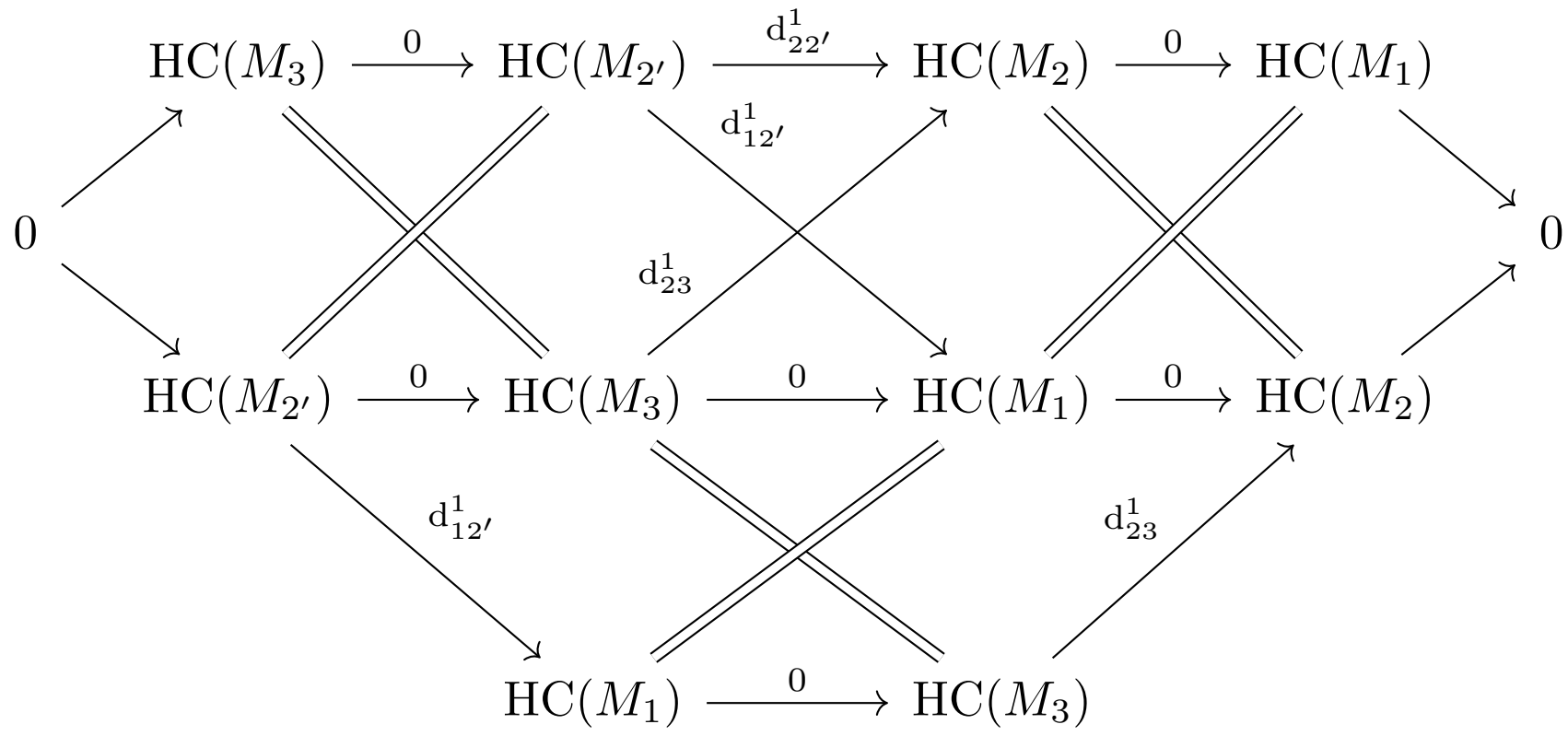
\begin{equation}
\begin{tikzcd}[column sep=small, row sep=1.8em]
\label{diamond4ever}
                        & \HC(M_3) \arrow[rr, "0"] \arrow[rrdd, no head, equals] &  & \HC(M_{2'}) \arrow[rr, "\dff^1_{22'}"] \arrow[rrdd, "\dff^1_{12'}"{pos=0.20}] \arrow[lldd, no head, equals] &              & \HC(M_2) \arrow[rr, "0"] \arrow[rrdd, no head, equals] &  & \HC(M_1) \arrow[rd] \arrow[lldd, no head, equals] &   \\
0 \arrow[ru] \arrow[rd] &                                                            &  &                                                                                                       &              &                                                            &  &                                                       & 0 \\
                        & \HC(M_{2'}) \arrow[rr, "0"] \arrow[rrdd, "\dff^1_{12'}"]   &  & \HC(M_3) \arrow[rr, "0"] \arrow[rrdd, no head, equals] \arrow[rruu, "\dff^1_{23}"{pos=0.20}]                &              & \HC(M_1) \arrow[rr, "0"] \arrow[lldd, no head, equals] &  & \HC(M_2) \arrow[ru]                                   &   \\
                        &                                                            &  &                                                                                                       & \hspace{1mm} &                                                            &  &                                                       &   \\
                        &                                                            &  & \HC(M_1) \arrow[rr, "0"]                                                                              &              & \HC(M_3) \arrow[rruu, "\dff^1_{23}"]                       &  &                                                       &
\end{tikzcd}
\end{equation}
One can follow $5$ distinct paths in the above diagram that provide a page-$1$ of a spectral sequence associated with the $5$ maximal chains in $\sO(\sM)$.
These are all possible linear Morse representations in this setting.
We will not consider perspective pairs in page-$1$ in this example since they do not provide additional information about differentials. In principle the first order differentials yield connecting orbits when considering different paths: $C(M_{2'},M_1)\neq \emptyset$, $C(M_{2'},M_2)\neq \emptyset$, and $C(M_3,M_2)\neq \emptyset$. We want to illustrate that higher order differentials can also fulfill this role.

Let us start with the path
\[
2\prec 3 \prec 1 \prec 2',
\]
which corresponds to the maximal chain:
$\bt\subsetneq \beta \subsetneq \delta \subsetneq \zeta\subsetneq \tp$.
From  Diagram \eqref{diamond4ever} (bottom row) we obtain page-$1$:
\[\begin{tikzcd}[column sep =large]
\cdots \arrow[r] & \Z\oplus \Z \arrow[r, "\id"  ] & \Z\oplus\Z \arrow[r, "0"  ] & \Z \arrow[r, "\begin{bmatrix}1\\0\end{bmatrix}" description] & \Z\oplus\Z \arrow[r] & \cdots
\end{tikzcd}\]
Computing homology yields page-$2$:
 
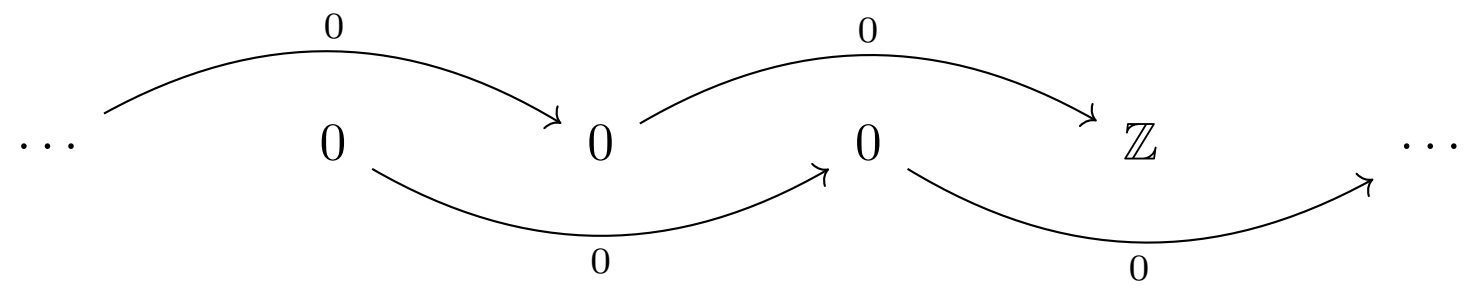
\[
\begin{tikzcd}[column sep=large]
\cdots \arrow[rr, "0", bend left] & 0 \arrow[rr, "0"', bend right] & 0 \arrow[rr, "0", bend left] & 0 \arrow[rr, "0"', bend right] & \Z & \cdots
\end{tikzcd}
\]
This maximal chain
  stabilizes at page-$2$.
  The higher order differentials do not play a role for this linear Morse representation.
There are two attractor-repeller pairs with $M_{2'},M_1$ and $M_3,M_2$, and since the differentials are non-trivial, connections are found by Corollary \ref{conn1}.
No additional  information about connections may be obtained from this Morse representation. However, the spectral net allows for more paths.

Consider the linear order
\[
2\prec 1\prec 3 \prec 2',
\]
which is associated with the maximal chain $\bt\subsetneq \beta\subsetneq \gamma\subsetneq \zeta \subsetneq \tp$. Diagram \eqref{diamond4ever} (middle row) yields the page-$1$:
 \[\begin{tikzcd}[column sep =large]
\cdots \arrow[r] & \Z\oplus \Z \arrow[r, "0"  ] & \Z \arrow[r, "0"  ] & \Z\oplus\Z \arrow[r, "0"  ] & \Z\oplus\Z \arrow[r] & \cdots
\end{tikzcd}\]
Computing homology and  establishing
the non-trivial second order differentials $\dff^2_{12'}$ and $\dff^2_{23}$ from the earlier data via $ji^{-1} k$, yields:
 
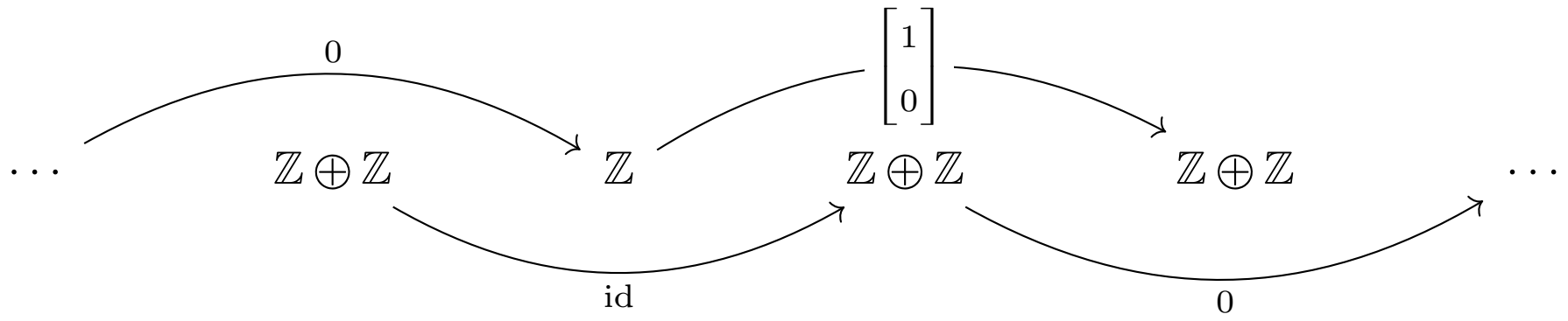
\[
\begin{tikzcd}[column sep=large]
\cdots \arrow[rr, "0", bend left] & \Z\oplus\Z \arrow[rr, "\id"', bend right] & \Z \arrow[rr, "\begin{bmatrix}1\\0\end{bmatrix}" description, bend left] & \Z\oplus\Z \arrow[rr, "0"', bend right] & \Z\oplus \Z & \cdots
\end{tikzcd}
\]
Computing homology reveals that the spectral sequence stabilizes at page-$3$. Here we have the opposite scenario, where all differentials on page-$1$ are trivial, and some differentials on page-$2$ are non-trivial. By Theorem \ref{connviahigherdiff},
since $\dff\colon S^{\smtp\smtp}_{\gamma\zeta}\rightarrow S^{\gamma\zeta}_{\smbt\beta}$ and
  $\dff\colon S^{\zeta\smtp}_{\beta\gamma}\rightarrow S^{\beta\gamma}_{\smbt\smbt}$
are not trivial, we obtain connections.
For the former
from $M_{2'}$ to either $M_1,M_3\in \zeta\setminus\beta$ and reciprocally from $M_{2'}, M_3 \in \tp\setminus\gamma$ to $M_1$. Since $M_{2'}\Vert  M_3$, we conclude  that  $C(M_{2'},M_{1})\neq \emptyset$, and since $M_1\Vert M_3$, we conclude again that $C(M_{2'},M_1)\neq \emptyset$. With  analogous reasoning we obtain from the second non-trivial differential that $C(M_3,M_2)\neq \emptyset$.

Finally, for the top row in Diagram \eqref{diamond4ever} we have the linear order
\[
1\prec 2\prec 2'\prec 3,
\]
which corresponds to the maximal chain $\bt\subsetneq \alpha\subsetneq \gamma\subsetneq \epsilon \subsetneq\tp$. This gives the page-$1$ spectral sequence
 \[\begin{tikzcd}[column sep =large]
\cdots \arrow[r] & \Z \arrow[r, "0"  ] & \Z\oplus \Z \arrow[r, "\id"  ] & \Z\oplus\Z \arrow[r, "0"  ] & \Z\oplus\Z \arrow[r] & \cdots
\end{tikzcd}\]
The first order differential $\dff^1_{22'}$ yields a connection from $M_{2'}$ to $M_2.$
Computing homology gives page-$2$, where all differentials are trivial:
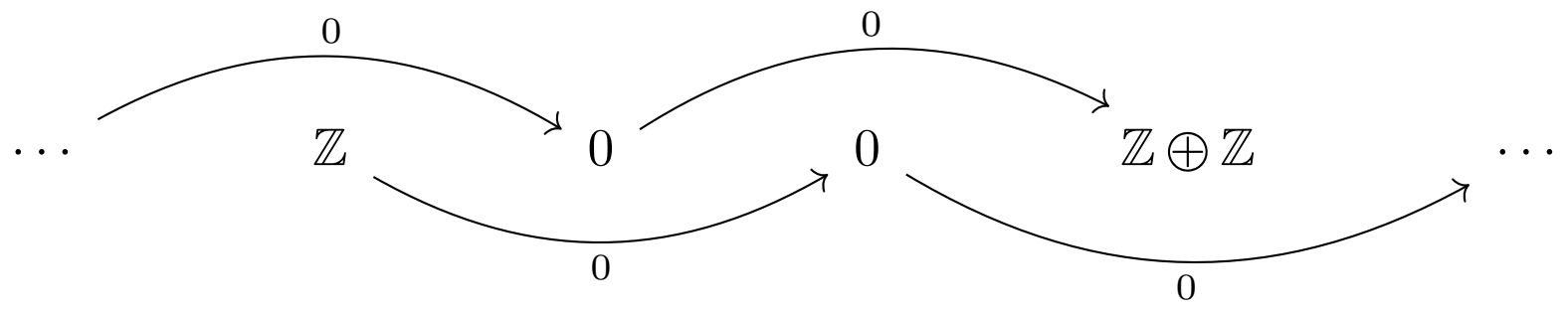
\[
\begin{tikzcd}[column sep=large]
\cdots \arrow[rr, "0", bend left] & \Z \arrow[rr, "0"', bend right] & 0 \arrow[rr, "0", bend left] & 0 \arrow[rr, "0"', bend right] & \Z\oplus \Z & \cdots
\end{tikzcd}
\]
Page-$3$ has a non-trivial differential since the spectral sequence has to converge at page-$4$. The differential $\dff^3_{13}$ need not be computed; it is non-trivial for the above reason.

  
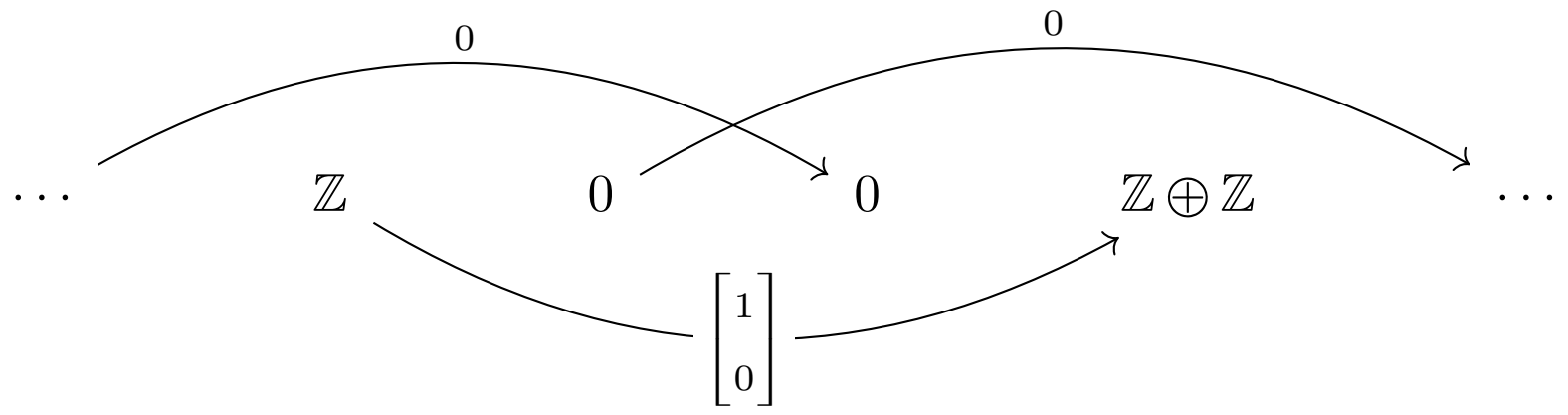
\[\begin{tikzcd}[column sep=large]
\cdots \arrow[rrr, "0", bend left] & \Z \arrow[rrr, "\begin{bmatrix}1\\0\end{bmatrix}"' description, bend right] & 0 \arrow[rrr, "0", bend left] & 0 & \Z\oplus \Z & \cdots
\end{tikzcd}\]
The differential
\[\dff\colon S^{\smtp\smtp}_{\alpha\epsilon}\rightarrow S^{\alpha\epsilon}_{\smbt\smbt}\]
is non-trivial. From Theorem \ref{conn2} we conclude that
the following sets $C(M_{2'},M_1)$, $C(M_{2'},M_2)$, and $C(M_3,M_2)$ are non-empty. Indeed, the sequence in Theorem~\ref{conn2} goes from $M_3$ to $M_1$ through comparable consecutive elements, and because $M_1\Vert M_2$, $M_1\Vert M_3$ and $M_{2'}\Vert M_3$, it passes through the pairs $\{M_1,M_{2'}\}$, $\{M_{2'},M_2\}$ and $\{M_2,M_3\}$.
Theorem~\ref{connviahigherdiff} gives two of these connections directly. There is a connection from $M_3$ to one of $M_1$, $M_2$, $M_{2'}\in \epsilon$, and reciprocally, from one of $M_3$, $M_{2'}$, $M_2 \in \tp\setminus\alpha$ to $M_{1}$. Given that $M_3\Vert  M_{2'}$, $M_3\Vert  M_{1}$ and $M_2\Vert  M_{1}$, that implies that  $C(M_3,M_2)\neq \emptyset$ and $C(M_{2'},M_1)\neq\emptyset$.
With the remaining two linear orders we obtain the same results. The above results have shown that $(\sM,\le)$ as given in Figure \ref{theflowonD3}[middle]
is a flow-order.
\end{example}

\begin{example}
\label{franexample}
In \cite{fran} the following system of differential equations is considered:
\[
\dot x=y,\quad \dot y = -\theta y - x(x-1/3)(1-x),\quad \dot\theta = \varepsilon(\theta'-\theta)(\theta''-\theta),
\]
with $x,y,\theta\in \R$ and $\varepsilon \ll 1$, and $0<\theta'\ll 1$, $\theta''\gg 1$.
The detailed description can be found in \cite[Ex.\ 6.2]{fran}, v.\ Fig.\ \ref{orbitsex1234}.
At $\theta =\theta'$ and   at $\theta = \theta''$ the system is known.
The Conley indices at the six rest points are isomorphic to $\Z$.
We consider Morse representations $(\sM,\le)$, where the Morse sets are the six fixed points as indicated in Figure \ref{orbitsex1234}, and we describe $\sM$ via a poset $(\sP,\le)$.
The
orbits given in Figure \ref{orbitsex1234} are known to exist which force order relations on $(\sM,\le)\cong (\sP,\le)$. From this point on we will describe $\sM$ via $\sP$.
From the given orbits the following sets are down-sets in $\sP$, and form a perspective pair in $\sO(\sP)$:
\[
  \alpha = \{1\}, \quad \beta = \{1,2\}, \quad \gamma = \{1,4\}, \quad \delta = \{1,2,4\}.
\]
We do not know whether a connection between $p=2$
and $p=6$ exists. As tail elements we take $\alpha_{-} = \emptyset$ and $\delta_{+} = \{1,2,4,6\}$. The order relations and the two possible sublattices are shown in Fig.~\ref{perspectivelattice}.

 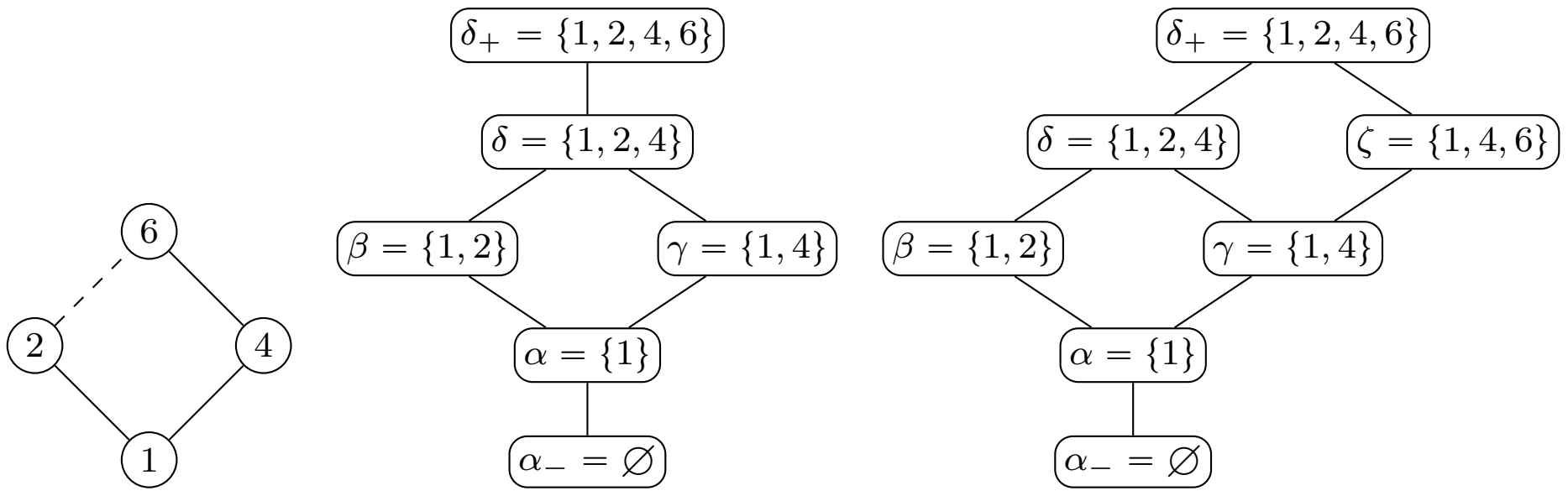
\begin{figure}[h] 
  \centering
\begin{tikzpicture}[scale=0.6, every node/.style={circle, draw, minimum size=2mm, inner sep=2pt, font=\small}]
  \node (1) at (0,0) {{\small $1$}};         
  \node (2) at (-1.5,1.5) {{\small $2$}};
  \node (4) at (1.5,1.5) {{\small $4$}};
  \node (6) at (0,3) {{\small $6$}};

  \draw (1) -- (4);
  \draw (1) -- (2);
  \draw[dashed] (2) -- (6);
  \draw (4) -- (6);
\end{tikzpicture}
\quad
\begin{tikzpicture}[
  scale=0.7,
  every node/.style={draw, rounded corners, fill=white, inner sep=2pt, font=\small}
]

  \node (D1) at (0,0) {$\alpha_{-} =\emptyset$};

  \node (D2) at (0,1.2) {$\alpha=\{1\}$};
  \node (D3) at (1.8,2.4) {$\gamma=\{1,4\}$};

  \node (D4) at (-1.8,2.4) {$\beta=\{1,2\}$};

  \node (D6) at (0.0,3.6) {$\delta=\{1,2,4\}$};

  \node (D8) at (0,4.8) {$\delta_{+}=\{1,2,4,6\}$};

  \draw (D1) -- (D2);
  \draw (D2) -- (D3);
  \draw (D2) -- (D4);
  \draw (D3) -- (D6);
  \draw (D4) -- (D6);
  \draw (D6) -- (D8);

\end{tikzpicture}
\quad
\begin{tikzpicture}[
  scale=0.7,
  every node/.style={draw, rounded corners, fill=white, inner sep=2pt, font=\small}
]

  \node (D1) at (0,0) {$\alpha_-=\emptyset$};

  \node (D2) at (0,1.2) {$\alpha=\{1\}$};
  \node (D3) at (1.8,2.4) {$\gamma=\{1,4\}$};

  \node (D4) at (-1.8,2.4) {$\beta=\{1,2\}$};

  \node (D6) at (0.0,3.6) {$\delta=\{1,2,4\}$};
  \node (D7) at (3.6,3.6) {$\zeta=\{1,4,6\}$};

  \node (D8) at (1.8,4.8) {$\delta_{+}=\{1,2,4,6\}$};

  \draw (D1) -- (D2);
  \draw (D2) -- (D3);
  \draw (D2) -- (D4);
  \draw (D3) -- (D7);
  \draw (D3) -- (D6);
  \draw (D4) -- (D6);
  \draw (D6) -- (D8);
  \draw (D7) -- (D8);

\end{tikzpicture}
\caption{The order relations in $\sP$ between the elements $1, 2, 4$ and $6$ [left].
The elements $2$ and $6$ may, or may not be ordered.
The two possible sublattices of $\sO(\sP)$, the first one supposing $2\leq 6$ [middle] and the second one for $2\Vert 6$ [right].
For both scenarios we have the square in the diagram indicating a perspective pair.}
 \label{perspectivelattice}
\end{figure}
 The Conley indices (with $R=\Z$) are used to define $G_pC \cong \Z$ for all $p\in \sP$. For the differentials between $E$-terms we have:
\[
\begin{aligned}
    d_{\beta\alpha\alpha_{-}} = \dff^1_{12}&\colon G_2 C \to G_1 C,\qquad
  \ell^{-1} d_{\delta_+ \delta \gamma} = \dff^1_{26} \colon G_6 C \to G_2 C;\\
  \dff_{\gamma\alpha\alpha_{-}} = \dff^1_{14} &\colon G_4 C \to G_1 C,\qquad
  \ell'^{-1} d_{\delta_+ \delta \beta} = \dff^1_{46} \colon G_6 C \to G_4 C.
  \end{aligned}
\]
By Theorem \ref{firstantocomm} we conclude that
\[
  \dff^1_{12}\, \dff^1_{26} + \dff^1_{14}\, \dff^1_{46} = 0.
\]
A direct analysis of the system shows that $E^\beta_{\alpha_-}=E^\gamma_{\alpha_-}=E^{\delta_+}_\beta=0$, and the exact sequences of the triples $\alpha_-\subset\alpha\subset\beta$, $\alpha_-\subset\alpha\subset\gamma$ and $\beta\subset\delta\subset\delta_+$ show that $d_{\beta\alpha\alpha_-}$, $d_{\gamma\alpha\alpha_-}$ and $d_{\delta_+\delta\beta}$ are isomorphisms, i.e.
 $\dff^1_{12}$, $\dff^1_{14}$ and $\dff^1_{46}$ are isomorphisms which implies that
$\dff^1_{26} \neq 0$. Since $\dff_{\delta_+\delta\gamma}=\ell\,\dff^1_{26}\neq0$, Corollary~\ref{conn1}, applied to $\gamma\subset\delta\subset\delta_+$ with $\delta\setminus\gamma=\{2\}$ and $\delta_+\setminus\delta=\{6\}$, implies that
there exists a connecting orbit from $M_6$ to $M_2$. Therefore,
  $2<6$ and hence the sublattice in Figure \ref{perspectivelattice}[middle] is the correct sublattice.
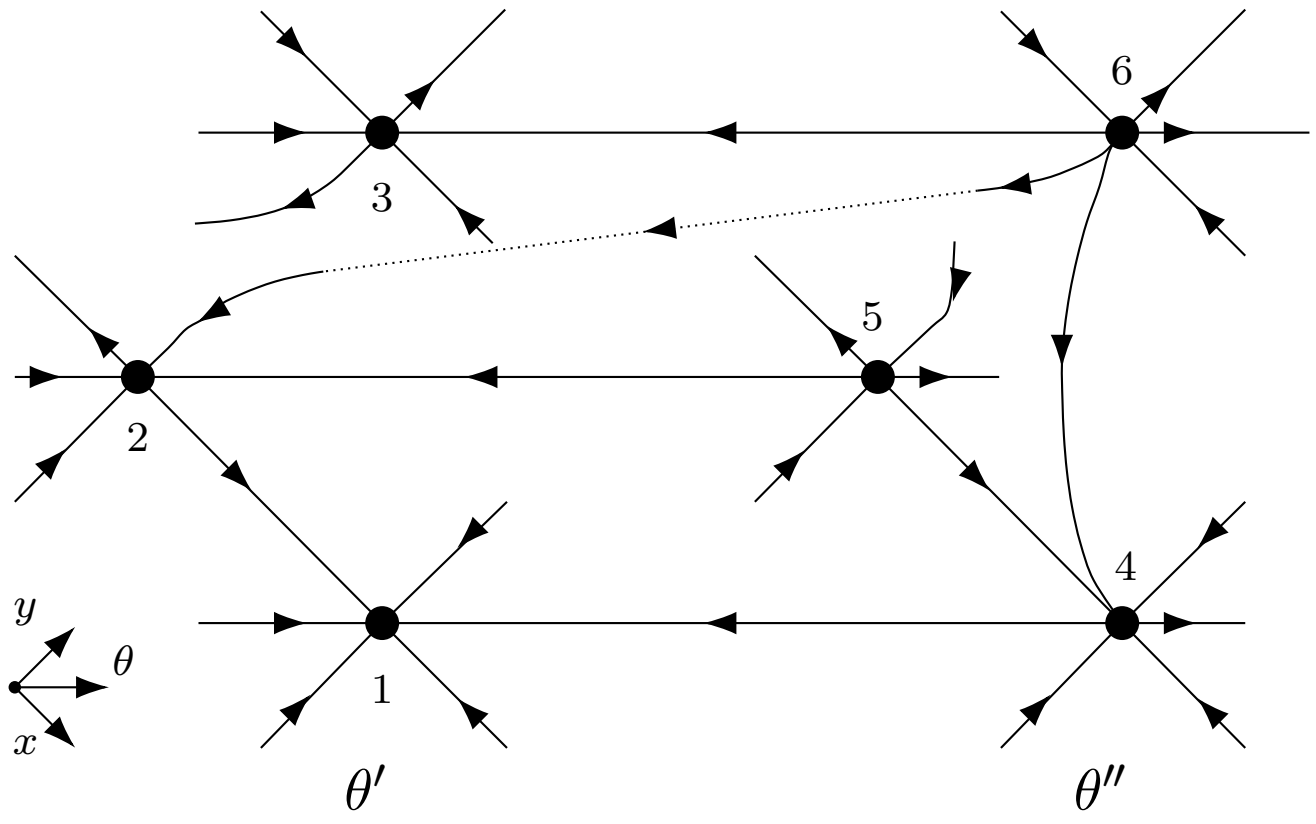
\begin{figure}[h]
\centering
  \begin{tikzpicture}[x=0.0118cm,y=-0.0118cm,thin,
    tip/.style={-{Latex[length=2.2mm,width=1.5mm]}}]
    \draw (122,85) -- (745,85);
    \draw (19,222) -- (571,222);
    \draw (122,360) -- (709,360);
    \draw[tip] (180,85) -- (183,85);
    \draw[tip] (408,85) -- (405,85);
    \draw[tip] (679,85) -- (682,85);
    \draw[tip] (43,222) -- (46,222);
    \draw[tip] (274,222) -- (271,222);
    \draw[tip] (542,222) -- (545,222);
    \draw[tip] (180,360) -- (183,360);
    \draw[tip] (408,360) -- (405,360);
    \draw[tip] (679,360) -- (682,360);
    \draw (157,17) -- (225,85) -- (294,16);
    \draw (225,85) -- (287,147);
    \draw (572,17) -- (640,85) -- (709,16);
    \draw (640,85) -- (709,154);
    \draw (19,154) -- (88,222) -- (225,360);
    \draw (19,292) -- (88,222);
    \draw (434,154) -- (503,222) -- (640,360);
    \draw (434,292) -- (503,222);
    \draw (157,430) -- (225,360) -- (295,292);
    \draw (225,360) -- (295,430);
    \draw (572,430) -- (640,360) -- (709,292);
    \draw (640,360) -- (709,430);
    \draw[tip] (180.879,40.879) -- (183,43);
    \draw[tip] (252.879,57.121) -- (255,55);
    \draw[tip] (267.121,127.121) -- (265,125);
    \draw[tip] (596.879,41.879) -- (599,44);
    \draw[tip] (665.879,59.121) -- (668,57);
    \draw[tip] (62.121,196.121) -- (60,194);
    \draw[tip] (45.879,264.121) -- (48,262);
    \draw[tip] (149.879,283.879) -- (152,286);
    \draw[tip] (476.121,195.121) -- (474,193);
    \draw[tip] (459.879,265.121) -- (462,263);
    \draw[tip] (563.879,283.879) -- (566,286);
    \draw[tip] (269.121,315.879) -- (267,318);
    \draw[tip] (182.879,402.121) -- (185,400);
    \draw[tip] (269.121,404.121) -- (267,402);
    \draw[tip] (684.121,315.879) -- (682,318);
    \draw[tip] (597.879,402.121) -- (600,400);
    \draw[tip] (684.121,404.121) -- (682,402);
    \draw[] plot[smooth, tension=0.6] coordinates {(120,136) (140,134.2) (160,130.4) (178,124) (196,113) (211,99) (225,85)};
    \draw[tip] (170.827,126.995) -- (168,128);
    \draw[tip] (682.121,127.121) -- (680,125);
    \draw[] plot[smooth, tension=0.6] coordinates {(546,146) (545.4,158) (544.5,173) (541.5,186) (533,194) (517,209) (503,222)};
    \draw[tip] (545.728,177.09) -- (545,180);
    \draw[] plot[smooth, tension=0.6] coordinates {(640,85) (633,95) (627.3,115) (620.2,135) (614.8,155) (610.6,175) (607.7,195) (606.2,215) (606.5,235) (607.3,255) (609,275) (611.9,295) (616.5,315) (624,336) (640,360)};
    \draw[tip] (606.2,214) -- (606.2,217);
    \draw[] plot[smooth, tension=0.6] coordinates {(640,85) (630,96) (620,101.7) (610,106) (600,109.8) (590,112.3) (580,114.4) (570,116.2) (559,117.5)};
    \draw[tip] (572.969,115.77) -- (570,116.2);
    \draw[densely dotted] (559,117.5) -- (192,162.9);
    \draw[tip] (373.977,140.232) -- (371,140.6);
    \draw[] plot[smooth, tension=0.6] coordinates {(192,162.9) (185,164) (175,166) (165,168.5) (155,171.9) (145,176) (135,181) (125,189.2) (115,195.2) (105,205.8) (88,222)};
    \draw[tip] (123.572,190.057) -- (121,191.6);
    \foreach \p in {(225,85),(640,85),(88,222),(503,222),(225,360),(640,360)} \fill \p circle (9.5);
    \node at (225,121) {$3$};
    \node at (640,50) {$6$};
    \node at (88,256) {$2$};
    \node at (500,188) {$5$};
    \node at (225,397) {$1$};
    \node at (642,328) {$4$};
    \node at (216,452) {\large$\theta'$};
    \node at (628,452) {\large$\theta''$};
    \fill (19,396) circle (3.5);
    \draw[tip] (19,396) -- (53,362);
    \draw[tip] (19,396) -- (72,396);
    \draw[tip] (19,396) -- (53,430);
    \node at (25,354) {$y$};
    \node at (80,381) {$\theta$};
    \node at (25,428) {$x$};
  \end{tikzpicture}
  \caption{Some flow-lines for the dynamical system given above. The spectral methods imply that $\dff^1_{26}\neq 0$, which implies that a connecting orbit from $p=6$ to $p=2$ exists. We denote this with the dotted arrow.}
  \label{orbitsex1234}
\end{figure}
In \cite[Ex.\ 6.2]{fran} this fact is obtained via a connection
matrix argument. With this example we show
that the same conclusion can be obtained using
spectral nets.
\end{example}

\begin{example}
\label{franexample2}
\cite{fran} discusses the   case of a flow on the M\"obius strip $S\cong I^2/_\sim$. We now show how derived differentials recognize connecting orbits in the case we use integral homology.
Consider a flow $\phi$ on $S$ with an attractor-repeller pair decomposition as displayed in Figure \ref{mobiusstrip}.
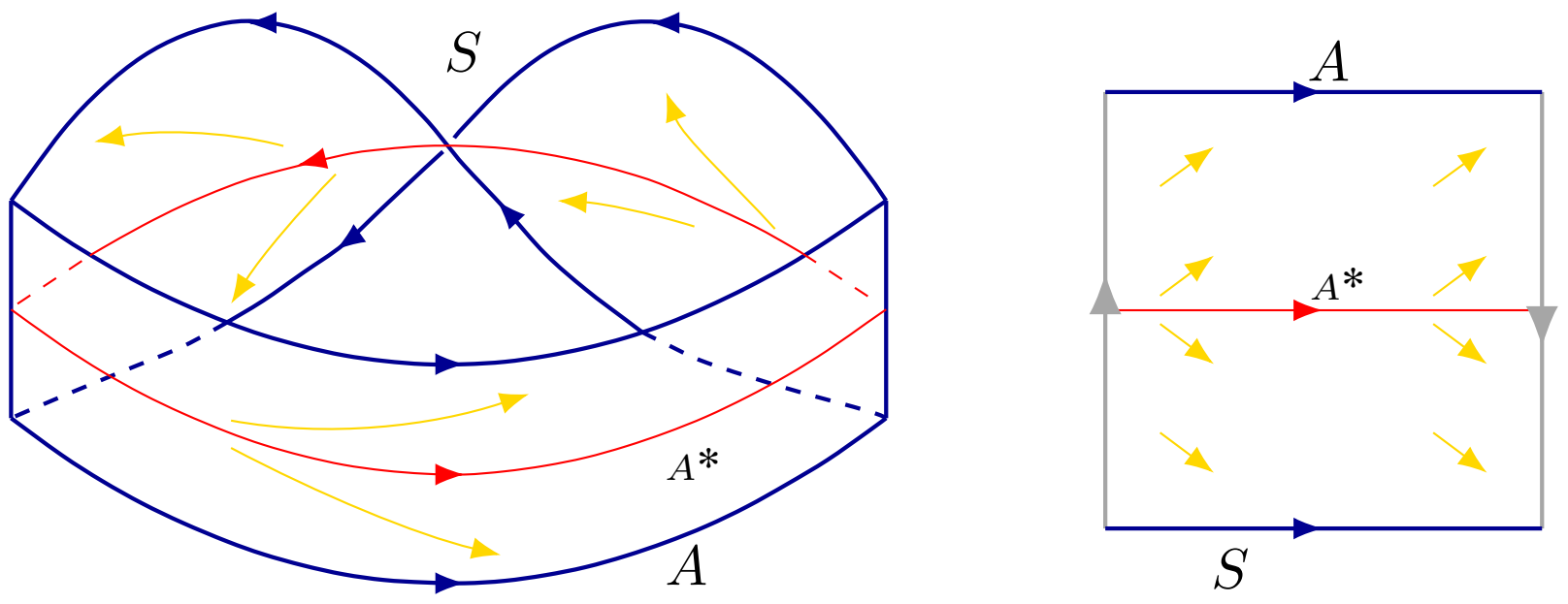
\begin{figure}[h]
\centering
  \begin{tikzpicture}[x=0.0055cm,y=-0.0055cm,
    btip/.style={sviiinavy,-{Latex[length=2mm,width=1.5mm]}},
    rtip/.style={red,-{Latex[length=2mm,width=1.5mm]}},
    ytip/.style={sviiigold,-{Latex[length=2mm,width=1.5mm]}},
    gtip/.style={gray!70,-{Latex[length=2.8mm,width=2.2mm]}}]
    \definecolor{sviiiblue}{RGB}{173,216,230}
    \definecolor{sviiigold}{RGB}{255,215,0}
    \definecolor{sviiipink}{RGB}{255,192,203}
    \definecolor{sviiigreen}{RGB}{144,238,144}
    \definecolor{sviiinavy}{RGB}{0,0,145}
    \draw[sviiinavy,thick] plot[smooth, tension=0.6] coordinates {(1123,265) (1100,232) (1050,169) (1000,119) (950,81) (900,56) (850,44) (800,44) (750,57) (700,82) (650,121) (600,172) (587,187)};
    \draw[sviiinavy,thick] plot[smooth, tension=0.6] coordinates {(573,204) (550,225) (500,272) (450,319) (400,354) (350,389) (305,416)};
    \draw[sviiinavy,thick] plot[smooth, tension=0.6] coordinates {(37,265) (100,179) (150,126) (200,86) (250,60) (305,44) (350,43) (400,54) (450,77) (500,113) (550,162) (580,198) (600,222) (650,275) (700,331) (750,375) (800,413) (820,428)};
    \draw[sviiinavy,thick] plot[smooth, tension=0.6] coordinates {(37,265) (100,310) (150,341) (200,369) (250,393) (305,416) (350,432) (400,446) (450,457) (500,464) (550,468) (600,468) (650,465) (700,458) (750,448) (800,435) (850,418) (900,397) (950,373) (1000,346) (1050,315) (1100,281) (1123,265)};
    \draw[sviiinavy,thick,dashed] plot[smooth, tension=0.6] coordinates {(305,416) (250,446) (175,477) (100,508) (37,535)};
    \draw[sviiinavy,thick,dashed] plot[smooth, tension=0.6] coordinates {(820,428) (850,443) (900,466) (1000,498) (1100,526) (1123,535)};
    \draw[sviiinavy,thick] plot[smooth, tension=0.6] coordinates {(37,535) (100,581) (150,613) (200,641) (250,665) (305,688) (350,704) (400,718) (450,729) (500,736) (550,739) (600,740) (650,737) (700,730) (750,720) (800,706) (850,689) (900,669) (950,645) (1000,617) (1050,587) (1100,552) (1123,535)};
    \draw[sviiinavy,thick] (37,265) -- (37,535) (1123,265) -- (1123,535);
    \draw[btip] (334,44) -- (330,44);
    \draw[btip] (834,44) -- (830,44);
    \draw[btip] (642.745,267.91) -- (640,265);
    \draw[btip] (443.277,323.706) -- (440,326);
    \draw[btip] (596,468) -- (600,468);
    \draw[btip] (596,740) -- (600,740);
    \draw[red] plot[smooth, tension=0.6] coordinates {(135,333) (150,324) (200,296) (250,271) (305,248) (350,233) (450,208) (500,201) (550,197) (600,197) (650,200) (700,207) (750,217) (800,230) (850,247) (950,292) (1000,319) (1020,331)};
    \draw[red,dashed] (45,393) -- (135,333) (1020,331) -- (1110,390);
    \draw[red] plot[smooth, tension=0.6] coordinates {(37,400) (100,445) (150,477) (200,505) (250,529) (305,552) (350,568) (400,582) (450,593) (500,600) (550,604) (600,605) (650,601) (700,594) (750,584) (800,570) (850,553) (900,533) (950,509) (1000,482) (1050,451) (1100,416) (1123,400)};
    \draw[rtip] (395.881,220.03) -- (392,221);
    \draw[rtip] (596,605) -- (600,605);
    \draw[ytip] (375,197) .. controls (300,177) and (210,177) .. (140,192);
    \draw[ytip] (440,232) .. controls (395,278) and (348,333) .. (310,393);
    \draw[ytip] (885,297) .. controls (830,280) and (770,268) .. (715,265);
    \draw[ytip] (985,300) .. controls (935,245) and (870,185) .. (850,130);
    \draw[ytip] (310,538) .. controls (430,560) and (560,545) .. (680,505);
    \draw[ytip] (310,572) .. controls (420,630) and (540,680) .. (645,705);
    \node at (598,80) {\large$S$};
    \node at (885,592) {\scriptsize$A^*$};
    \node at (875,718) {\large$A$};
    \draw[gray!70,thick] (1395,130) -- (1395,672) (1937,130) -- (1937,672);
    \draw[sviiinavy,thick] (1395,130) -- (1937,130) (1395,672) -- (1937,672);
    \draw[red] (1395,401) -- (1937,401);
    \draw[btip] (1661,130) -- (1665,130);
    \draw[btip] (1661,672) -- (1665,672);
    \draw[rtip] (1661,401) -- (1665,401);
    \draw[gtip] (1395,359) -- (1395,355);
    \draw[gtip] (1937,443) -- (1937,447);
    \draw[ytip,thin] (1463,247) -- (1530,198);
    \draw[ytip,thin] (1463,383) -- (1530,333);
    \draw[ytip,thin] (1463,418) -- (1530,468);
    \draw[ytip,thin] (1463,553) -- (1530,603);
    \draw[ytip,thin] (1802,247) -- (1869,198);
    \draw[ytip,thin] (1802,383) -- (1869,333);
    \draw[ytip,thin] (1802,418) -- (1869,468);
    \draw[ytip,thin] (1802,553) -- (1869,603);
    \node at (1672,92) {\large$A$};
    \node at (1550,722) {\large$S$};
    \node at (1685,368) {\scriptsize$A^*$};
  \end{tikzpicture}
  \caption{The M\"obius strip $S$ with the middle periodic orbit repelling (red) and the boundary is an attracting periodic orbit (blue).}
  \label{mobiusstrip}
\end{figure}
The attractor-repeller pair decomposition yields an $E$-system with poset $\sP=\{1,2\}$, $1<2$ and $M_1=A$ (blue) and $M_2=A^*$ (red). The lattice of down-sets is given by $\sO(\sP)= \{\bt,\alpha,\tp\}$, where $\bt=\emptyset$,
$\alpha=\{1\}$ and $\tp=\{1,2\}$. The $E$-terms are given by the Conley indices: $E^\alpha_\smbt=\HC(A)\cong H(S^1)$, $E^\smtp_\alpha = \HC(A^*)\cong H(S,\partial S)$ and $E^\smtp_\smbt = \HC(S)\cong H(S^1)$.
The exact triangle
\[
\begin{tikzcd}
 E^\alpha_\smbt \arrow[rr, "i"] &                                 & E^\smtp_\smbt \arrow[ld, "j"] \\
                                           & E^\smtp_\alpha \arrow[lu, "k"] &                                      \end{tikzcd}
\]
may be expressed as long exact sequence of graded Conley indices:
\tikzset{
  snake/.style={
    rounded corners=8pt,
    to path={
      let \p1=(\tikztostart.east), \p2=(\tikztotarget.west), \n1={(\y1+\y2)/2} in
        (\p1) -- (\x1+3.5ex, \y1)
              -- (\x1+3.5ex, \n1)
              -- node[midway, above, inner sep=3pt] {#1} (\x2-3.5ex, \n1)
              -- (\x2-3.5ex, \y2)
              -- (\p2)
    }
  }
}
\[
\begin{tikzcd}[row sep=2.6em, column sep=3.4em]
0 \arrow[r] & \HC_2(A) \arrow[r] & \HC_2(S) \arrow[r] & \HC_2(A^*) \arrow[dll, snake=$k_2$] & \\
& \HC_1(A) \arrow[r] & \HC_1(S) \arrow[r] & \HC_1(A^*) \arrow[dll, snake=$k_1$] & \\
& \HC_0(A) \arrow[r] & \HC_0(S) \arrow[r] & \HC_0(A^*) \arrow[r] & 0
\end{tikzcd}
\]
Let us start with $\bR=\Z$-coefficients. The Conley indices are given by:

\begin{center}
\begin{tabular}{cccc}
\toprule
degree & $\HC(A;\Z)$ & $\HC(A^*;\Z)$  & $\HC(S;\Z)$ (M\"obius strip)\\
\midrule
$0$      & $\Z$ & $0$      & $\Z$\\
$1$      & $\Z$ & $\Z_2$  & $\Z$\\
$\ge 2$  & $0$  & $0$      & $0$\\
\bottomrule
\end{tabular}
\end{center}

If we implement these modules in the long exact sequence we obtain:
\[
\begin{tikzcd}[row sep=2.6em, column sep=3.4em]
0 \arrow[r] & 0 \arrow[r] & 0 \arrow[r] & 0 \arrow[dll, snake=$0$] & \\
& \Z \arrow[r] & \Z \arrow[r] &\Z_2 \arrow[dll, snake=$0$] & \\
& \Z \arrow[r] & \Z \arrow[r] & 0 \arrow[r] & 0
\end{tikzcd}
\]
The connecting homomorphism $k_{1}=0$ since $\Hom(\Z_2,\Z)=0$, and $k_2=0$ since $\HC_2(A^*)=0$, and therefore $k=0$. The above exact triangle then becomes a short exact sequence
\[
\begin{tikzcd}
0 \arrow[r] & E^{\alpha}_{\smbt} \arrow[r, "i"] & E^{\smtp}_{\smbt} \arrow[r, "j"] & E^{\smtp}_{\alpha} \arrow[r] & 0,
\end{tikzcd}
\]
which is non-split. By Lemma \ref{critconniff} this implies that $C(A^*,A)\neq \emptyset$. The latter is based on the algebraic interpretation of Morse representations. However, purely in terms of algebra of spectral systems this is manifested in terms of the derived differential
$\kappa\colon E^\smtp_\alpha\to E^\alpha_\smbt[1]$. To compute $\kappa$ we use a projective resolution for $\Z_2$ given by $\cdots \to 0\to \Z\xrightarrow[\times 2]{\dff_1} \Z\xrightarrow[\mod 2]{\epsilon} \Z_2\to 0$.
In the diagram below we choose $f_0=\id$ and consequently $f_1=\id$:
\[
\begin{tikzcd}[column sep=large, row sep=large]
\cdots \arrow[r] & 0 \arrow[r, "0"] \arrow[d, "0"'] & \Z \arrow[r, "\times 2"] \arrow[d, "f_1", dotted] & \Z \arrow[d, "f_0", dotted] \arrow[r, "\mod 2"] & \Z_2 \arrow[r] \arrow[d,equal] & 0 \\
\cdots \arrow[r] & 0 \arrow[r]                      & \Z \arrow[r, "i"]                                       & \Z \arrow[r, "j"]                                     & \Z_2 \arrow[r]           & 0
\end{tikzcd}
\]
The diagram for $P_\bullet\to \Z[1]$ is given by
\[
\begin{tikzcd}[column sep=large, row sep=large]
\cdots \arrow[r] & 0 \arrow[r, "0"] \arrow[d, "0"'] & \Z \arrow[r, "\times 2"] \arrow[d, "\id"] & \Z \arrow[d, "0"] \\
\cdots \arrow[r] & 0 \arrow[r]                      & \Z \arrow[r]                         & 0
\end{tikzcd}
\]
which is a non-trivial representative of $\kappa$ in $\Hom_{\Kh(\sRmodZ)}(P_\bullet,\Z[1])$: a null-homotopy would be a homomorphism $h\colon\Z\to\Z$ with $h\circ(\times2)=\id$, which does not exist. This shows that the derived differential is non-trivial.
In \cite{fran} the equivalent of the derived differential $\kappa$ is a non-trivial entry in the $C$-connection matrix. This example illustrates that spectral systems also record the connecting orbits algebraically.

The same calculations can be carried out in different coefficient rings providing different conclusions. For example if we use the field coefficients in $\bR=\Z_2$ we obtain

\begin{center}
\begin{tabular}{cccc}
\toprule
degree & $\HC(A;\Z_2)$ & $\HC(A^*;\Z_2)$ & $\HC(S;\Z_2)$ (M\"obius strip)\\
\midrule
$0$     & $\Z_2$ & $0$     & $\Z_2$\\
$1$     & $\Z_2$ & $\Z_2$ & $\Z_2$\\
$2$     & $0$     & $\Z_2$ & $0$\\
$\ge 3$ & $0$     & $0$     & $0$\\
\bottomrule
\end{tabular}
\end{center}
In  this case the derived differential is necessarily trivial since all modules are free. One concludes in this case that $k\neq 0$, which by Theorem \ref{conn1forder} algebraically records the existence of a connecting orbit.
If we use $\bR=\Z_4$ the theory in \cite{fran} is not applicable since $\Z_4$ is not a PID. However, the analysis here still provides a non-trivial connecting homomorphism, since $k_2\colon\HC_2(A^*;\Z_4)\cong\Z_2\to\HC_1(A;\Z_4)\cong\Z_4$ is injective, as well as a non-trivial derived differential, and therefore an algebraic manifestation of a connecting orbit. If we consider $\bR=\R$, or for example $\bR=\Z_3$ the Conley index $\HC(A^*)$ is trivial and the connecting orbit cannot be detected algebraically.
The derived differential may also contain information about the connecting orbit in the Morse representation besides its existence. This is why considering different coefficient rings may be useful.
\end{example}

\section*{Disclosure of AI use}
This manuscript was written without using generative AI. Before submission of the article the content was carefully checked by the authors using
 Claude Opus 5.5 Extra. After using this tool, the authors have reviewed the suggested edits to the content. The authors take full responsibility for the content of the publication.

\begin{acknowledgements}
The  authors wish to thank Thomas Rot for pointing out the subject of spectral systems and for many inspiring discussions on the subject.
\end{acknowledgements}

\bibliographystyle{spmpsci}      
\bibliography{references}   


\end{document}